\documentclass[11pt,a4paper,onecolumn]{article}

\usepackage{a4wide}
\usepackage{microtype}

\usepackage[x11names,dvipsnames]{xcolor}
\colorlet{method6color}{magenta}
\definecolor{method5color}{RGB}{220,20,60}

\usepackage{amsmath,amsfonts,bm}

\def\1{\bm{1}}

\DeclareMathAlphabet{\mathsfit}{\encodingdefault}{\sfdefault}{m}{sl}
\SetMathAlphabet{\mathsfit}{bold}{\encodingdefault}{\sfdefault}{bx}{n}

\def\gC{{\mathcal{C}}}

\def\gF{{\mathcal{F}}}

\def\gK{{\mathcal{K}}}
\def\gL{{\mathcal{L}}}

\def\gN{{\mathcal{N}}}

\def\gP{{\mathcal{P}}}

\newcommand\N{{\mathbb{N}}}
\newcommand\bN{{\mathbb{N}}}

\newcommand{\E}{\mathbb{E}}

\newcommand{\R}{\mathbb{R}}

\newcommand{\cF}{\mathcal{F}}

\newcommand{\nrm}[1]{\|#1\|}

\DeclareMathOperator*{\argmin}{arg\,min}

\newcommand{\cz}{c(z^{k+1})}

\DeclareMathOperator{\gph}{gph}
\DeclareMathOperator{\dist}{dist}
\DeclareMathOperator{\diag}{diag}

\renewcommand{\Pr}{\mathbb{P}}
\newcommand{\comp}[1]{#1^{\mathrm{c}}}
\renewcommand{\epsilon}{\varepsilon}
\newcommand{\avgnoise}{\bar{\varepsilon}}

\usepackage{amsmath,amssymb}
\usepackage{amsthm,mathtools}
\usepackage{aliascnt}
\usepackage{aligned-overset}

\usepackage{graphicx}
\usepackage{subcaption}
\usepackage{tikz}
\usepackage{pgfplots}
\usepgfplotslibrary{groupplots,fillbetween}
\pgfplotsset{compat=1.18}

\usepackage{booktabs}
\usepackage{multirow}

\usepackage{algorithm}
\usepackage[noend]{algpseudocode}
\algrenewcommand{\algorithmicrequire}{\textbf{Input:}}

\usepackage{natbib}
\setcitestyle{authoryear,round,citesep={;},aysep={,},yysep={;}}

\usepackage{url}
\usepackage{hyperref}
\usepackage[capitalize,noabbrev,nameinlink]{cleveref}

\hypersetup{
    colorlinks=true,
    linkcolor=Blue4,
    citecolor=Green4,
    pdftitle={Stochastic Extragradient Methods Beyond Lipschitz Continuity}
}

\theoremstyle{plain}
\newtheorem{theorem}{Theorem}[section]

\newaliascnt{lemma}{theorem}
\newtheorem{lemma}[lemma]{Lemma}
\aliascntresetthe{lemma}

\newaliascnt{proposition}{theorem}
\newtheorem{proposition}[proposition]{Proposition}
\aliascntresetthe{proposition}

\newaliascnt{corollary}{theorem}
\newtheorem{corollary}[corollary]{Corollary}
\aliascntresetthe{corollary}

\theoremstyle{definition}

\newaliascnt{definition}{theorem}
\newtheorem{definition}[definition]{Definition}
\aliascntresetthe{definition}

\newaliascnt{assumption}{theorem}
\newtheorem{assumption}[assumption]{Assumption}
\aliascntresetthe{assumption}

\newaliascnt{example}{theorem}
\newtheorem{example}[example]{Example}
\aliascntresetthe{example}

\theoremstyle{remark}

\newaliascnt{remark}{theorem}
\newtheorem{remark}[remark]{Remark}
\aliascntresetthe{remark}

\crefname{assumption}{assumption}{assumptions}
\Crefname{assumption}{Assumption}{Assumptions}

\title{Stochastic Extragradient Methods Beyond Lipschitz Continuity}
\author{%
    Konstantinos Oikonomidis%
    \thanks{Equal contribution. \\
        KU Leuven, Department of Electrical Engineering,
        STADIUS Center for Dynamical Systems, Signal Processing,
        and Data Analytics.\newline \emph{E-mails:}
        \texttt{%
            \{%
            \href{mailto:konstantinos.oikonomidis@esat.kuleuven.be}%
                 {konstantinos.oikonomidis},%
            \href{mailto:jan.quan@esat.kuleuven.be}%
                 {jan.quan},%
            \href{mailto:mathijs.schuurmans@esat.kuleuven.be}%
                 {mathijs.schuurmans},%
            \href{mailto:panos.patrinos@esat.kuleuven.be}%
                 {panos.patrinos}%
            \}\href{mailto:konstantinos.oikonomidis@esat.kuleuven.be,jan.quan@esat.kuleuven.be,mathijs.schuurmans@esat.kuleuven.be,panos.patrinos@esat.kuleuven.be}{@esat.kuleuven.be}%
        }%
    }%
    \qquad
    Jan Quan\footnotemark[1]\\[0.25em]
    Mathijs Schuurmans\footnotemark[1]%
    \qquad
    Panagiotis Patrinos
}

\pgfplotsset{compat=1.18}
\begin{document}

\maketitle

\begin{abstract}
    Stochastic variational inequality problems have received significant attention due to their ability to model large-scale machine learning tasks. Typically, the operators involved are assumed to be globally Lipschitz, with only a handful of works extending beyond the standard setting to the so-called \((L_0,L_1)\)-Lipschitz setting. In this work, we propose a novel stochastic past extragradient algorithm for variational inequality problems under a broad regularity assumption that is even less restrictive than \((L_0,L_1)\)-Lipschitz continuity and provide convergence guarantees under bounded variance of the stochastic oracle. More precisely, we show a.s. convergence with suitable stepsize schedules and obtain convergence rates in expectation and in probability for fixed horizons. Finally, in numerical simulations we compare our proposed algorithm against standard methods. 
\end{abstract}

\section{Introduction}

Variational inequalities (VI) and root-finding problems \citep{facchinei2003finite} 
provide a common framework to address modern-day learning problems that go beyond classical 
unconstrained minimization, such as constrained minimization, multiplayer games, and minimax optimization. Typical applications include adversarial network training, distributionally robust optimization and optimal transport. In their simplest (unconstrained) form, VI problems reduce to finding the zeroes of an operator \(F\), i.e., finding \(x^\star\) such that \(F(x^\star) = 0\).

Some of the most well-known methods to tackle VIs are the extragradient (EG) \citep{korpelevich1977extragradient} and the past extragradient (PEG) method \citep{popov1980modification}. 
Despite their prevalence, these methods are typically analyzed under monotonicity and Lipschitz continuity of \(F\), which can be restrictive in modern learning problems.
For that reason, a large part of the recent literature has been devoted to lifting these assumptions. 

A prominent direction for relaxing monotonicity is based on \emph{Minty-type} conditions.
The Minty variational inequality (MVI) assumes the existence of a solution \(x^\star\) such that
\begin{equation} \label{eq:minty_vi} 
\langle F(x),x-x^\star\rangle \geq 0
\qquad \text{for all }x,
\end{equation}
whereas the traditional monotonicity assumption requires the same inequality between all pairs of points. This assumption was relaxed to the \emph{weak Minty} condition in \citet{diakonikolas2021efficient}, which allows the right-hand side of \eqref{eq:minty_vi} to be proportional to \(\|F(x)\|^2\),
and it is shown that a relaxed extragradient method, EG+, still converges under this substantially weaker assumption. This was then further developed by \citet{pethick2023escaping} for constrained and regularized problems; and by \citet{pethick2023solving} in the stochastic setting. Multi-step exploration was introduced in \citet{fan2024weaker}, and refined fixed-point analyses and schemes interpolating between extragradient and proximal-point methods in the weak Minty framework have been considered in \citet{alacaoglu2024revisiting} and \citet{pethick2025efficient}, respectively.

A complementary line of work has sought to relax the global Lipschitz continuity assumption instead. In particular, the $(L_0,L_1)$-smoothness assumption \citep{zhang2019gradient,zhang2020improved} was introduced to better describe the smoothness conditions of neural networks.
\citet{vankov2024generalized} extended this framework to variational inequalities and established convergence guarantees for extragradient-type methods under generalized smoothness.
Most recently, \citet{choudhury2025extragradient,quan2026dual} provided analyses of extragradient methods for \((L_0,L_1)\)-Lipschitz root-finding problems in the deterministic setting.
Nevertheless, much less attention has been devoted to the stochastic setting. In \citet{vankov2025generalized}, $(L_0,L_1)$-Lipschitz problems have been studied in the stochastic setting under a $p$-quasi-sharpness condition, which is sufficiently strong to yield convergence guarantees even for suitably clipped projection/gradient descent-ascent (GDA) schemes, whereas the analysis in \citet{quan2026dual} requires increasing batch sizes.

In this paper we aim at addressing this gap by proposing and analyzing a nonlinearly preconditioned stochastic past extragradient method under a regularity assumption that is in fact even less restrictive than $(L_0,L_1)$-Lipschitz continuity. Concretely, our contributions are the following:

\begin{itemize}
    \item We propose a stochastic past extragradient method for solving VI problems involving operators that are not Lipschitz continuous. The class of operators we consider, described in \Cref{def:prec_lipschitz}, is larger than that of the recently introduced class of \((L_0,L_1)\)-Lipschitz continuous operators \citep{vankov2025generalized,choudhury2025extragradient}. Assuming bounded variance of the stochastic oracle, we show almost sure convergence to solutions of the VI in \Cref{thm:main_as} and obtain convergence rates in expectation for a natural preconditioned residual in \Cref{cor:main_rate}. Moreover, we get convergence guarantees for the standard residual in probability in \Cref{thm:convergence-in-prob,thm:convergence-in-prob-unconstrained}, complementing our in expectation results.  
    \item We present various examples of operators satisfying our regularity assumptions, describing both theoretical constructions and potential applications. Furthermore, we compare our proposed \Cref{alg:polyak-sep} against standard solvers such as the two-stepsize versions of the stochastic extragradient and stochastic optimistic gradient methods in various numerical experiments ranging from simple quartic games to neural optimal transport and GAN training.
\end{itemize}

\subsection{Notation}
Let \(\N = \{0, 1, 2, \dots \}\) denote the naturals including \(0\). We denote by $\langle\cdot,\cdot\rangle$ the standard Euclidean inner product on $\mathbb{R}^n$ and by $\|\cdot\|$ the standard Euclidean norm on $\mathbb{R}^n$ as well as the spectral norm for matrices. For a closed, convex set \(C\subseteq \mathbb{R}^n\) and \(\bar x \in C\) we denote its normal cone \(N_C(\bar x) = \{y \in \mathbb{R}^n \mid \langle y,x-\bar x \rangle \leq 0: \forall x \in C\}\). The squared distance to a set $E$ will be denoted as $\dist^2(x, E) = \inf_{\bar x\in E} \|x - \bar x\|^2$. For \(T: \mathbb{R}^n \rightrightarrows \mathbb{R}^n\) a set-valued mapping, we define its graph as \(\gph T\coloneqq \{(x,y) \in \mathbb{R}^n \times \mathbb{R}^n\mid y \in T(x)\}\) and its inverse as \(T^{-1}(y) \coloneqq \{x \in \mathbb R^n \mid (x,y) \in \gph T\}\).
\section{Related work}

\textbf{Generalized smoothness.} Our work is closely related to works on generalized smoothness conditions in optimization. Starting from \citet{zhang2019gradient} that introduced the notion of \((L_0,L_1)\)-smoothness to better describe functions where the Hessian norm grows with the gradient norm, there has been a substantial body of work devoted to this topic. In \citet{wang2023convergence,faw2023beyond,koloskova2023revisiting,pethick2025generalizedgradientnormclipping,riabinin2025gluon}, the authors focused on analyzing the convergence properties of existing algorithms under this new condition, while \citet{gorbunov2025methods,vankov2025optimizing} suggest new algorithms tailored to this smoothness condition. Our proposed method \Cref{alg:polyak-sep} utilizes independent sampling for the update direction and the normalization, which has been used in the context of generalized smoothness in \citet{gaash2026convergence} and for variational inequality problems in \citet{vankov2025generalized}. Methods for \((L_0,L_1)\)-Lipschitz variational inequalities were also analyzed recently in \cite{greenberg2026extragradient}.

\textbf{Nonlinear (dual space) preconditioning.} The generalized Lipschitz continuity inequality we utilize in \Cref{def:prec_lipschitz} is influenced by the literature on dual space or nonlinear preconditioning for gradient descent and the corresponding anisotropic descent inequality. In \cite{maddison2021dual} dual space preconditioning for convex optimization was introduced based on a notion of generalized cocoercivity. The method was extended to the nonconvex additive nonsmooth setting in \cite{laude2025anisotropic} and was analyzed under the notion of anisotropic smoothness which was itself first introduced in \citet{laude2023dualities}. The connection of this new descent condition with \((L_0,L_1)\)-smoothness that we build upon in this work was first established in \cite{oikonomidis2025nonlinearly}. \cite{oikonomidis2025nonlinearlymomentum} analyzed the method with Polyak momentum and under simple noise assumptions, while \cite{bodard2025escaping} studied the saddle point escape properties of nonlinearly preconditioned gradient methods and provided various new examples of \((L_0,L_1)\)-smooth functions. A nonlinearly preconditioned (anisotropic) proximal point algorithm for finding zeroes of monotone operators was extensively studied in \citet{laude2023anisotropic}. The work closest to ours is \cite{quan2026dual} that studies root-finding and VI problems under generalized cocoercivity and preconditioned Lipschitz continuity conditions. We stress that the stochastic results of this paper require increasing batch sizes, something that we avoid in our current work.

\textbf{Stochastic methods for VI problems.} In \citet{juditsky2011solving} the Stochastic Mirror-Prox was analyzed and finite-time guarantees for monotone variational inequalities using stochastic prediction–correction steps were obtained. Subsequently, stochastic extragradient type methods have been studied under increasingly general conditions, including pseudomonotonicity, unbounded feasible sets, and nonuniform stochastic noise \citep{iusem2017extragradient}. In \citet{hsieh2020explore} the convergence of a two-timescale version of EG was established and in \citet{bohm2022two} the convergence of the stochastic forward-backward-forward method for convex-concave minimax problems was studied. Moreover, an extension of the stochastic EG method without large batch sizes to the weak Minty setting was proposed in \citet{pethick2023solving}, while \citet{gorbunov2022clipped,gorbunov2023high} established the convergence of clipped versions of EG under heavy-tailed noise. A major line of work, that is not directly related to our paper, focuses on variance reduction techniques for finite-sum VIs \citep{palaniappan2016stochastic,chavdarova2019reducing,alacaoglu2022stochastic,emmanouilidis2024stochastic}. A significant part of the related literature originates from learning in games with a particular emphasis on training GANs \citep{daskalakis2017training,mertikopoulos2018optimistic,gidel2018variational} and reinforcement learning \citep{daskalakis2020independent,kalogiannis2026policy}.

\section{Problem statement and algorithm}
\subsection{Problem statement}
We study stochastic methods for tackling VI problems of the form: find \(x^\star \in C\) such that
\begin{equation} \tag{VIP} \label{eq:vi}
    \langle F(x^\star),x-x^\star \rangle \geq 0 \qquad \forall x \in C,
\end{equation}
where \(C\subseteq \mathbb{R}^n\) is a closed, convex set. We assume that such an \(x^\star\) always exists and denote by \(S\) the set of all such points. Recall that \eqref{eq:vi} is also equivalent to finding the zeroes of the operator $F + N_C$, so that $S = (F + N_C)^{-1}(0)$.

Throughout the text, the jointly measurable map $G: \mathbb{R}^n \times \Xi \to \mathbb{R}^n$, 
which satisfies \(\E[\|G(x,\xi)\|] < +\infty\), $\forall x \in \mathbb{R}^n$,
denotes a stochastic oracle for the population operator $F(x) = \E[ G(x, \xi) ]$.
Regarding the regularity of \(F\) we do not follow the standard definition of \((L_0, L_1)\)-Lipschitz continuity of \citet[Assumption 2.1]{vankov2024generalized} but rather the preconditioned Lipschitz continuity inequality of \citet[Assumption 3.2]{quan2026dual}, which is an extension of \citet[Assumption 2.5]{oikonomidis2025nonlinearlymomentum}.
\begin{definition} \label{def:prec_lipschitz}
    We say that \(F\) is \(\gK\)-Lipschitz continuous with \(\gK(y) = \frac{y}{L_0 + L_1\|y\|}, L_0 >0, L_1 \geq 0\) if the following inequality holds for all \(x, \bar x \in \mathbb{R}^n\):
    \begin{equation} \label{eq:prec_Lipschitz}
        \|\gK(F(x)) - \gK(F(\bar x))\| \leq \|x-\bar x\|.
    \end{equation}
\end{definition}
We denote \(T(x)\coloneqq \gK(F(x))\) and note that \eqref{eq:prec_Lipschitz} is equivalent to the \(1\)-Lipschitz continuity of \(T\). This condition lies at the center of our analysis, since the Lipschitz continuity of \(T\), together with the fact that \(T\) and \(F\) share zeroes, allows us to design algorithms acting directly on \(T\). Nevertheless, in the stochastic setting we consider, this is not trivial: applying the preconditioner \(\gK(G(x,\xi))\) leads to a biased estimate of \(T(x)\), i.e., \(\E[\gK(G(x,\xi))] \neq T(x)\). We explore this issue in more detail in the following subsection.

At the same time, \Cref{def:prec_lipschitz} serves another purpose: it defines an operator class that is larger than the class of \((L_0,L_1)\)-Lipschitz operators. In the case where \(F \in \gC^1(\mathbb{R}^n)\), this can be seen by the differential characterizations. Indeed, for \((L_0,L_1)\)-Lipschitz continuity, the necessary and sufficient condition is that
\begin{equation} \label{eq:l0l1_snd_order}
    \|\nabla F(x)\| \leq L_0 + L_1\|F(x)\| \qquad \forall x \in \mathbb{R}^n,
\end{equation}
while for \Cref{def:prec_lipschitz}, \(\|\nabla T(x)\| \leq 1\), i.e., 
\begin{equation} \label{eq:prec_lipschitz_snd_order}
    \frac{1}{L_0+L_1\|F(x)\|}\left \| \left (I-\frac{L_1}{L_0+L_1\|F(x)\|}\frac{F(x)F(x)^\top}{\|F(x)\|} \right )\nabla F(x) \right\| \leq 1,
\end{equation}
for \(x \in \mathbb R^n \setminus F^{-1}(0)\), while by continuity \(\nabla T(x) = \nabla F(x)/ L_0\) for \(x \in F^{-1}(0)\). It is straightforward to check that \eqref{eq:l0l1_snd_order} implies \eqref{eq:prec_lipschitz_snd_order}. The reverse implication does not hold, however, with a standard example being from \cite{bodard2025escaping}: considering \(F(x) = 2x\exp(x^2)\), there cannot be any \(L_0,L_1\) such that \eqref{eq:l0l1_snd_order} holds globally. On the other hand, for \(L_0=2,L_1=1\), \(\|\nabla T(x)\| \leq 1\) and thus \Cref{def:prec_lipschitz} holds. This relation between the two classes holds without requiring \(F \in \gC^1(\mathbb{R}^n)\), as is detailed in \Cref{app:prec_lipschitz}.
 \begin{figure}[htbp]
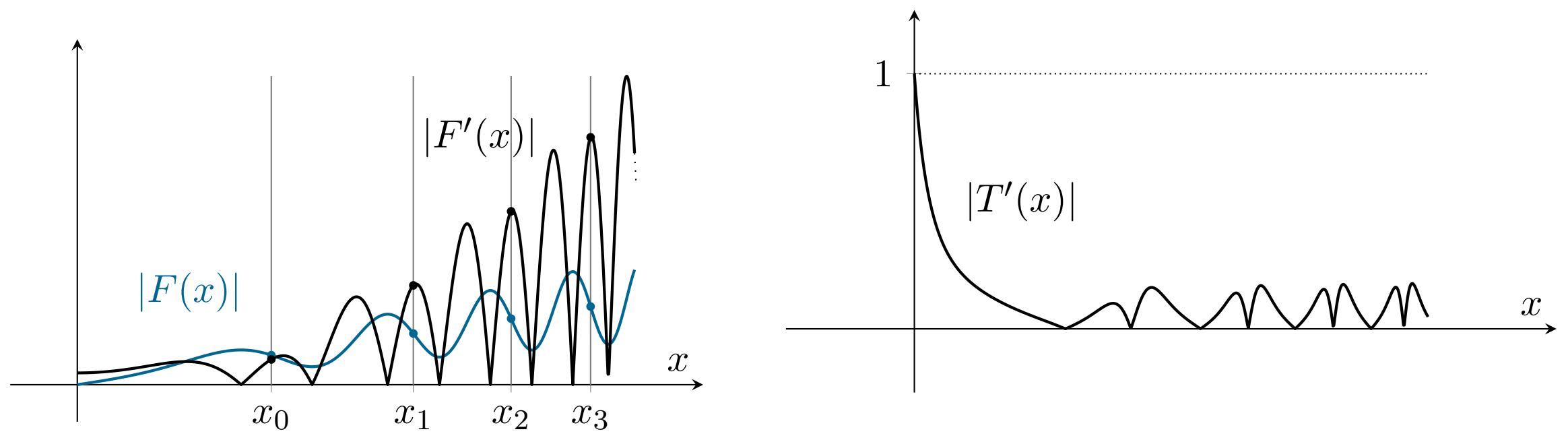

  \centering

  \begin{minipage}[b]{0.49\textwidth}
    \centering
    \input{figures/ex1_gradient_curvature.tex}
  \end{minipage}
  \hfill
  \begin{minipage}[b]{0.49\textwidth}
    \centering
    \input{figures/ex1_T.tex}
  \end{minipage}
  \caption{
    Comparison between preconditioned Lipschitz continuity and \((L_0,L_1)\)-Lipschitz continuity for \cref{ex:sin}.
    \textit{(Left)} \((L_0,L_1)\)-Lipschitz continuity fails although $\| F(x) \|$ grows at most linearly.
    \textit{(Right)} The gradient norm of the preconditioned operator $T$ remains bounded by \(1\).
  }
  \label{fig:l0l1failure}
\end{figure}
Examples of preconditioned Lipschitz operators can be found in the related literature, with \cite{maddison2021dual,vankov2025optimizing,gorbunov2025methods,oikonomidis2025nonlinearly,bodard2025escaping} covering the functional case \(F = \nabla f\) for \(f:\mathbb{R}^n \to \mathbb R\) and \cite{vankov2024generalized,choudhury2025extragradient} covering the root-finding case. We provide some more interesting examples and applications in what follows, with more detailed derivations and further examples in \Cref{app:prec_lipschitz}. We start with an example that satisfies \Cref{def:prec_lipschitz} but is not \((L_0,L_1)\)-Lipschitz although it has at most linear growth, in contrast to \(2x\exp(x^2)\) that we considered before.
\begin{example}[Failure of \((L_0,L_1)\)-Lipschitz with at most linear growth] \label{ex:sin}
    Let \(F\) be the gradient of the function \(f(x)=\tfrac{x^2}{2}+\tfrac{1-\cos(x^2/2)}{2}\), i.e.,
    \[
      F(x)=x\big(1+ \tfrac{1}{2}\sin(x^2/2)\big).
    \]
    \(F\) satisfies \Cref{def:prec_lipschitz} with constants \(L_0=1\) and \(L_1 =2\).
    Indeed,
    note that $|F(x)| \geq \frac{1}{2}|x|$,
    and since $F'(x) = 1 + \tfrac{1}{2}\sin(x^2 / 2) + \tfrac{x^2}{2} \cos(x^2 / 2)$,
    using \(|\sin x| \leq | x |\), we have \(| F'(x) | \leq 1 + \tfrac{3}{4} x^2 \leq 1 + x^2\). Therefore,
    \(
      |T'(x)| = \frac{|F'(x)|}{(1+2|F(x)|)^2} \leq \frac{1 + x^2}{(1 + |x|)^2} \leq 1.
    \)
    On the other hand, \(F\) is not \((L_0,L_1)\)-Lipschitz for any finite combination of \(L_0\) and \(L_1\).
    This can be seen by taking a sequence \((x_{\nu})_{\nu \in \bN}\) of points \(x_{\nu}= \sqrt{2(2\nu+1)\pi}\).
    Then, \(\sin(x_\nu^2 / 2) = 0\) and \( \cos (x^2_\nu / 2) = -1 \), and thus, 
    \(F(x_\nu) = x_\nu \) and \(|F'(x_{\nu})| = \frac{x_\nu^2}{2} - 1\). Therefore, 
    \(\frac{|F'(x_\nu)|}{|F(x_\nu)|} \to \infty\) and \eqref{eq:l0l1_snd_order} cannot hold for finite 
    \(L_0>0\) and \(L_1\geq0\).
    A visualization of this discrepancy can be found in \Cref{fig:l0l1failure}.
\end{example}

Another interesting example is a saddle operator where the coupling between the two variables 
is also nonlinear.
\begin{example}[Quartic game with nonlinear cross-player interaction] \label{example:quartic_nonlinear}
    Let $x\in\mathbb{R}^n,y\in\mathbb{R}^m$ and $B\in\mathbb{R}^{n\times m}$. Consider the saddle function \(L(x,y) = \tfrac{1}{4}\|x\|^4 - \tfrac{1}{4}\|y\|^4 + \tfrac{\lambda}{2}\|x\|^2\|y\|^2 + \langle x,By \rangle\) with \(\lambda \in \mathbb R\). In this setting, the two players interact not only via the linear coupling but also via
    \(\tfrac{\lambda}{2}\|x\|^2\|y\|^2\),
    leading to higher-order growth in the shared cost,
    rather than only in the player-specific terms. 
    The corresponding saddle operator,
    \[
        F(x, y) = \begin{pmatrix} (\|x\|^2 + \lambda \|y\|^2)x + B y \\ (\|y\|^2 - \lambda \|x\|^2)y - B^\top x \end{pmatrix}
    \]
    satisfies \cref{def:prec_lipschitz} with \(L_0=\|B\|+3+2|\lambda|\) and \(L_1=6+4|\lambda|\).
\end{example}

As is now standard in the literature,
our analysis assumes the Minty variational inequality,
an extension of monotonicity requiring that the involved inequality
only holds between general points and the solution set, rather than all pairs of points.
\begin{definition} \label{def:minty}
    We say that \(F\) satisfies the Minty variational inequality if for all \(x^\star \in S\), the following inequality holds:
    \begin{equation} \label{eq:minty}
        \langle y,x-x^\star \rangle \geq 0 \qquad \forall (x,y) \in \gph (F + N_C).
    \end{equation}
\end{definition}
Specializing to the case where \(C= \mathbb{R}^n\), it is clear that if \(F\) satisfies \Cref{def:minty}, then so does \(T\). 
This observation motivates the design of a stochastic EG method 
by applying standard techniques to \(T\) instead of \(F\). In general, this estimator is biased, that is, \(\E[\gK(G(x,\xi))] \neq T(x)\), so directly applying an analysis that requires an unbiased oracle for $T$ does not guarantee convergence to a solution of the original problem. Our algorithm design therefore utilizes a more careful approach, which we detail 
in the next section.
\begin{remark}
    In the related literature, the Minty variational inequality typically requires only one \(x^\star \in S\) in \eqref{eq:minty} and not the whole set \(S\) as in \Cref{def:minty}. Most of our analysis goes through with this more relaxed assumption, except for the a.s.\ full sequence convergence result in \Cref{thm:main_as}. We thus choose to keep this more restrictive assumption to make the presentation clear and concise and avoid different conditions between convergence theorems.
\end{remark}
\begin{algorithm}[t]
\caption{Stochastic Past Extra Preconditioned Gradient (SPEPG)}
\label{alg:polyak-sep}
\begin{algorithmic}[1]

\Require Stochastic oracle $G$.
\Require Initial point \(x^0 \in C\); $L_0>0$, $L_1\ge0$.
\Require Deterministic $0<\alpha_k\le1$, $\gamma,\omega_k>0$, with $\alpha_0=1$.

\State \(m^{-1}\gets 0\), \(z^0 \gets x^0\)
\State Draw \(\xi^0,\zeta^0\) i.i.d.\ and compute \(\widehat T(z^0;\xi^0,\zeta^0)=\frac{G(z^0,\xi^0)}{L_0 + L_1\|G(z^0,\zeta^0)\|}\).

\For{\(k=0,1,\ldots,\)}
    
        \State \(m^k \gets (1-\alpha_k)m^{k-1} + \alpha_k \widehat T(z^k;\xi^k,\zeta^k)\)
        \State \(z^{k+1} \gets P_C(x^k - \gamma m^k)\)
        \State Draw \(\xi^{k+1},\zeta^{k+1}\) i.i.d.\ and compute \(\widehat T(z^{k+1};\xi^{k+1},\zeta^{k+1})=\frac{G(z^{k+1},\xi^{k+1})}{L_0 + L_1\|G(z^{k+1},\zeta^{k+1})\|}\)
        \State \(x^{k+1} \gets P_C(x^k - \omega_k \widehat T(z^{k+1};\xi^{k+1},\zeta^{k+1}))\)

\EndFor

\end{algorithmic}
\end{algorithm}

\subsection{Algorithm}
Our proposed extension of the past extragradient method is presented in \Cref{alg:polyak-sep}. The intuition behind our construction is relatively simple: assuming for simplicity that \(C=\mathbb{R}^n\) and
\(F\) is Lipschitz and monotone, the two-stepsize past extragradient (PEG)
iteration takes the form:
\begin{align*}
    z^{k+1} &= x^k - \gamma F(z^k)
    \\
    x^{k+1} &= x^k - \omega F(z^{k+1}).
\end{align*}
Since in our case, \(T\) is Lipschitz and satisfies the Minty VI, the next step is to replace \(F\) with \(T\) in the recursion displayed above. However, as discussed in the previous section,
since in our stochastic setting we only have access to unbiased stochastic oracles 
\(G\) of \(F\) and not \(T\), directly utilizing \(\gK \circ G\) as the oracle 
would not lead to convergence. 
Instead, we construct an estimator of \(T(z^k)\) using two independently 
drawn samples \(\xi, \zeta\): given any \(z \in \mathbb{R}^n\), we define
\begin{equation} \label{eq:def_hat_T}
    \widehat T(z; \xi, \zeta)  \coloneqq  \frac{G(z,\xi)}{L_0 + L_1\|G(z,\zeta)\|}. 
\end{equation}
Despite the separate samples, in general \(\E[\widehat T(z; \xi, \zeta)] \neq T(z)\).
Nevertheless, due to the independence of the samples, it is easy to show that
\(\E[\widehat T(z; \xi, \zeta)] = \E\Big[\frac{L_0+L_1\|F(z)\|}{L_0+L_1\|G(z, \zeta)\|}\Big]T(z)\).
Denoting \(c(z) \coloneqq \E\Big[\frac{L_0+L_1\|F(z)\|}{L_0+L_1\|G(z, \zeta)\|}\Big]\), we have that \(c(z) > 0\) which implies that in expectation the direction of the estimator is correct.
In fact, under bounded variance (\Cref{assum:bounded_variance}),
we have uniform lower and upper bounds for \(c\) as we state formally in the following section.
Finally, we use Polyak momentum in the update of \(z^{k+1}\) instead of using
\(\widehat T(z^k;\xi^k,\zeta^k)\) directly in order to make the most of its variance reduction properties.

\begin{remark} \label{rem:eg_variant}
    \Cref{alg:polyak-sep} can be extended to a standard extragradient type method where in the update of \(z^{k+1}\) we would use two more samples to compute \(\widehat T(x^k; \xi_x^k, \zeta_x^k)\). Nevertheless, since this approach adds significant computational cost to the method, we decided not to include it in the current paper. In the unconstrained case, and with the same stepsize \(\omega = \gamma\), it is well known that PEG and the optimistic gradient methods have equivalent forms \citep[Equation (OG)]{gorbunov2022last}. Since we are interested in the general VI problem and we utilize the two-sequence form of PEG in our analysis, we choose the name stochastic past extragradient for our proposed method.
\end{remark}
\begin{remark}
    There are multiple ways in which $L_0,L_1$ in \cref{alg:polyak-sep} can be chosen as hyperparameters. First of all, if $F$ is $(L_0,L_1)$-Lipschitz, then it is also $(\tilde{L}_0, L_1)$-Lipschitz for any $\tilde{L}_0 \geq L_0$ as is clear from \eqref{eq:l0l1_snd_order} and \eqref{eq:l0l1}. Similarly, if $F$ satisfies \cref{def:prec_lipschitz} for some $L_0>0$ and $L_1\geq 0$ and we have $\tilde{L}_0 \geq L_0$ and $\tilde{L}_1^2 L_0 \geq L_1^2 \tilde{L}_0$, then $F$ also satisfies \cref{def:prec_lipschitz} with $\tilde{L}_0$ and $\tilde{L}_1$. In general, $L_0$ and $L_1$ are connected through some nonlinear relation which may differ depending on the operator, see for example \cite[Example 3]{choudhury2025extragradient} in the $(L_0,L_1)$-Lipschitz case.
\end{remark}

\section{Convergence analysis} 
\label{sec:convergence_analysis}

In this section we present our main theoretical results regarding the convergence properties of our proposed scheme. The proofs can be found in \Cref{app:analysis}.
We next formulate the main assumptions. We begin with the bounded variance assumption for the stochastic oracle, a standard assumption in the related literature,
as well as the aforementioned regularity conditions on the operator \(F\).

\begin{assumption} \label{assum:bounded_variance}
    For every $x \in \mathbb{R}^n$, we have \(F(x) = \E[G(x,\xi)]\) and
        \[
            \E[\|G(x, \xi) - F(x)\|^2] \leq \sigma^2.
        \]
\end{assumption}
\begin{assumption} \label{assum:prec_lip_minty}
    \(F\) is continuous and satisfies \Cref{def:prec_lipschitz} and \Cref{def:minty}.
\end{assumption}

Under the bounded variance assumption we can formalize our results on the estimation of the preconditioned operator:
\begin{lemma} \label{thm:stoch_bounds}
    Let $\xi, \zeta$ be i.i.d. random variables, 
    and suppose that \Cref{assum:bounded_variance} holds. Define
    \begin{equation*} 
        c(x) \coloneqq \E \left[ \frac{L_0 + L_1 \|F(x)\|}{L_0 + L_1 \|G(x, \zeta)\|} \right], \quad x\in\mathbb{R}^n
    \end{equation*}
    and $\delta \coloneqq \frac{L_1\sigma}{L_0}$. Then, for all $x \in \mathbb{R}^n$,
    \begin{equation} \label{eq:c_bound}
          \E[\widehat T(x; \xi, \zeta)] = c(x)T(x), \qquad \frac{1}{1+\delta} \leq c(x) \leq 1+\delta,
    \end{equation}
    \begin{equation} \label{eq:cond_bound_var}
        \E\left[\|\widehat T(x;\xi, \zeta)-c(x) T(x)\|^2\right] \leq \frac{2\sigma^2}{L_0^2}, \qquad \E \left[\|\widehat T(x; \xi, \zeta)\|^2 \right] \leq c(x)^2 \|T(x)\|^2 + \frac{2\sigma^2}{L_0^2}.
    \end{equation}
\end{lemma}
Finally, we assume the following step-size conditions, which follow directly 
from our analysis.
\begin{assumption} \label{assum:param_choice}
    The deterministic parameters satisfy 
    \[
          0 < \gamma \leq \frac{1}{4(1+\delta)}, \quad 0 < \theta\leq \frac{\gamma}{3(1+\delta)^2},\quad \omega_k = \theta\alpha_{k+1}, \quad 0 < \alpha_{k+1} \leq 1, \quad k \in \N
    \]
    where $\delta = \frac{L_1\sigma}{L_0}$ with $\sigma$ as in \Cref{assum:bounded_variance} and $L_0,L_1$ as in \Cref{assum:prec_lip_minty}.
\end{assumption}

Under these assumptions, we obtain the following main result.

\begin{theorem} \label{thm:main}
    Let \Cref{assum:bounded_variance,assum:prec_lip_minty,assum:param_choice} hold true. 
    Let, moreover, $(z^k)_{k\in\mathbb{N}}$ be generated by \Cref{alg:polyak-sep}.
    Then, for every $K\geq 1$ and every $x^\star\in S$,
    \[
        \frac{1}{S_K} \sum_{k=0}^{K-1} \alpha_{k+1} \E[\dist^2 (0, (T + N_C)(z^{k+1}))] \leq \frac{18}{7}\frac{ (1+\delta)^2 }{\theta\gamma S_K}\overline {V}_0(x^\star) + \frac{159}{14}(1+\delta)^4 \frac{2\sigma^2}{L_0^2} \frac{Q_K}{S_K},
    \]
    where $\overline{V}_0(x^\star)$ is some constant that depends on $x^\star$. Here, $S_K \coloneqq \sum_{k=0}^{K-1}\alpha_{k+1}$ and $Q_K\coloneqq \sum_{k=0}^{K-1}\alpha_{k+1}^2$.
\end{theorem}

\begin{corollary} \label{cor:main_rate}
  Let \Cref{assum:bounded_variance,assum:prec_lip_minty,assum:param_choice} hold true. Let, moreover, $(z^k)_{k = 1}^{K}$ be the sequence of iterates 
  generated by running \Cref{alg:polyak-sep} for $K\geq 1$ iterations with $\alpha_{k+1} = K^{-1/2}$ for all $k \in \{0,1,\dots,K-1\}$. Then,
    \[
        \frac{1}{K}\sum_{k=0}^{K-1}\E[\dist^2(0,(T+N_C)(z^{k+1}))] = O(K^{-1/2}).
    \]
\end{corollary}
\begin{remark}\label{rmk:hardness-expectation}
    Some remarks are in order regarding \Cref{cor:main_rate}. To begin with, our analysis leads to obtaining a rate on the mixed residual involving \(T\) instead of the usual \(F\). This is a natural outcome since our regularity condition is the Lipschitzness of \(T\) which we combine with the fact that it is Minty. Moreover, \(\dist^2(0,(T+N_C)(z)) = 0\) for some \(z \in \mathbb{R}^n\) implies \(-T(z) \in N_C(z)\) and thus that \(-F(z) \in N_C(z)\) due to the cone properties of \(N_C(z)\). In fact, it is not possible to transfer this rate to $F$, even when $C=\mathbb{R}^n$. To see this, consider $F : \mathbb{R} \to \mathbb{R} : x \mapsto x^3$, which satisfies \cref{def:prec_lipschitz} with $L_0=L_1=1$. Then, for $x^0 = 0$, and the stochastic oracle $G(x,\xi) = x^3 + \xi$ where $\xi$ has zero mean, finite variance and $\E[|\xi|^3] = \infty$ (e.g., when $\xi$ follows a Student-$t_3$ distribution), we have $z_1 = -\gamma \xi^0/(1 + |\zeta^0|)$, hence $\E[|F(z^1)|] = +\infty$. Nevertheless, we complement our results by showing a.s.\ full sequence convergence under standard momentum schedules and obtaining convergence rates in probability for the standard metrics.
\end{remark}

\begin{theorem} \label{thm:main_as}
     Let \Cref{assum:bounded_variance,assum:prec_lip_minty,assum:param_choice} hold true. 
     Let, moreover, $(x^k)_{k\in\mathbb{N}}$ be the sequence of iterates generated by \Cref{alg:polyak-sep} with
    \[
        \sum_{k=0}^\infty \alpha_{k+1} = +\infty, \qquad \sum_{k=0}^\infty \alpha_{k+1}^2 < +\infty.
    \]
    Then, the sequence \((x^k)_{k \in \mathbb N}\) converges almost surely to some possibly random \(x^\star \in S\).
\end{theorem}

\begin{theorem} \label{thm:convergence-in-prob}
  Let \Cref{assum:bounded_variance,assum:prec_lip_minty,assum:param_choice} hold true. Let, moreover, $(z^k)_{k = 1}^{K}$ be the sequence of iterates 
  generated by running \Cref{alg:polyak-sep} for $K\geq 1$ iterations with $\alpha_{k+1} = K^{-1/2}$ for all $k \in \{0,1,\dots,K-1\}$.
  Let $J$ be uniformly distributed on $\{0,\ldots,K-1\}$ and independent
  of the iterates.
  Then, for every $\eta\in(0,1)$, there exists a constant
  $C_\eta>0$, independent of $K$, such that, for every $K\geq1$,
  \[
    \Pr\left\{
      \dist^2\left(
        0,(F+N_C)(z^{J+1})
      \right)
      >
      C_\eta K^{-1/2}
    \right\}
    \leq \eta.
  \]
\end{theorem}

\begin{remark}
From the proof in \cref{sec:proof-convergence-in-prob},
it can be observed that although the constant $C_{\eta}$ is independent from $K$, it may grow rapidly as $\eta$ approaches 0,
depending on the growth of $\| F(z) \|$ away from the solution.
This reveals that for some choices of $F$ a small finite-$K$ error may be difficult to obtain for very small violation probabilities $\eta$.
We emphasize, however, that this is not a consequence of our analysis.
To see this, consider again the example from \cref{rmk:hardness-expectation} with $\xi \sim t_3$.
Since $z^1=-\gamma \xi^0/(1+|\zeta^0|)$, we have (conditioned on $\{|\zeta^0| \leq 1\}$) that $|F(z^1)|^2\ge (\gamma^6/64)|\xi^0|^6$.
By independence between $\xi^0$ and $\zeta^0$,
we thus have for any $\beta > 0$, 
\begin{equation}\label{eq:prob-bound-ex}
\mathbb P \{ |F(z^1)|^2 > \beta \} \geq \mathbb P(|\zeta^0|\le 1)
\mathbb P \left\{ |\xi^0|>\frac{2}{\gamma}\beta^{1/6}\right\}.
\end{equation}
Now since $\xi^0 \sim t_3$, there exists $c_1 > 0$, such that $\mathbb P\{ |\xi^0| > t \} \geq c_1 t^{-3}$ for sufficiently large $t$. Thus, from \eqref{eq:prob-bound-ex}, 
there exists $c > 0$ such that
$\mathbb P \{ |F(z^1)|^2>\beta \} \geq c \beta^{-1/2}$
for sufficiently large $\beta$.
Consequently, if $\mathbb P\{ |F(z^1)|^2 > C_\eta \} \leq \eta$, 
then $C_{\eta} \geq c^2 \eta^{-2}$ for all sufficiently small $\eta$.
This construction can be repeated with operators of increasingly fast growth (e.g., taking $F(x) = x^p$ with $p$ odd) and noise distributions with increasingly heavy--yet finite-variance--tails (e.g., $\xi \sim t_\nu$, with $\nu > 2$).
This yields arbitrarily severe dependence of the constant $C_\eta$ on $\eta$.
Hence, under the present assumptions (which even cover super-exponential operators $F$, see \cref{ex:prec_greater}), no uniform ($F$-independent) quantitative dependence of $C_\eta$ on $\eta$ can be expected without additional restrictions controlling either the growth of $F$ or the tails of the stochastic oracle.
\end{remark}

In the unconstrained case (i.e., $C = \mathbb{R}^n$), we can strengthen the conclusions 
of \cref{thm:convergence-in-prob} and additionally 
obtain an $O(K^{-1/2})$ bound on the failure probability given a tolerance $\epsilon > 0$.

\begin{theorem} \label{thm:convergence-in-prob-unconstrained}
  Let the hypotheses of \cref{thm:convergence-in-prob} hold true.
  If, additionally, $C = \mathbb{R}^n$, then,
  there exists a constant $C_T$, independent of $K$, such that for all $\epsilon > 0$,
  \[
  \Pr\{ \nrm{F(z^{J+1})} > \epsilon \} \leq \frac{C_T (L_0 + L_1 \epsilon)^2}{\sqrt{K}\epsilon^2}.
  \]
\end{theorem}

\begin{remark}
    The constant $C_T$ in \cref{thm:convergence-in-prob-unconstrained} may be taken as $$C_T \coloneq \frac{18(1+\delta)^2}{7\theta \gamma} \Big(\overline{V}_0(x^\star) + \frac{53}{12}\frac{\theta\gamma (1+\delta)^22\sigma^2}{L_0^2}\Big)$$ where $\overline{V}_0$ is as in \eqref{eq:bound-EV0}. This constant is independent of $K$ and depends on the initialization, noise level, normalization parameters, and chosen stepsizes. Its polynomial dependence on $\delta$ reflects the deterioration of the bound when the noise is large relative to the normalization constant.
\end{remark}

To the best of our knowledge, our results are the first under \((L_0,L_1)\)-Lipschitz continuity, and in fact under the less restrictive condition \Cref{def:prec_lipschitz}, that merely require the Minty variational inequality \Cref{def:minty} to hold. The work \citet{vankov2025generalized}
that also tackles stochastic VI problems under generalized continuity assumptions requires a strong growth assumption in the form of \citet[Assumption 2.4]{vankov2025generalized}. 

\section{Numerical experiments}

In this section we compare our proposed \Cref{alg:polyak-sep} 
with standard algorithms like stochastic extragradient (SEG) and stochastic optimistic gradient (SOG),
where for SEG and SOG we use two timescales in the stepsizes to allow for more flexibility.
In addition, we compare against the clipped SEG method analyzed in \cite{vankov2025generalized}.
Details and further experimental results including GAN training on MNIST and neural optimal transport are presented in \Cref{app:additional_experiments}. Having established the theoretical guarantees of our proposed algorithm, in this section we aim to test its performance in a variety of problems and compare it against standard methods. To that aim we do not use the theoretically obtained stepsizes, but rather tune them, as is the case in practice. Moreover, we consider problems of practical interest where the assumptions we consider potentially do not hold. We provide some synthetic experiments in \cref{app:synth} that satisfy our theoretical assumptions.

\subsection{Quartic game with nonlinear interaction} \label{sec:quartic_nonlinear}
We first consider the setting of \Cref{example:quartic_nonlinear}. We define \(F(x,y) = \tfrac{1}{N}\sum_{i=1}^N F_i(x,y)\), \(x,y \in \mathbb R^{d}\) with 
\[
    F_i(x,y) = \begin{pmatrix} (\|x\|^2 + \lambda \|y\|^2)x + B_i y \\ (\|y\|^2 - \lambda \|x\|^2)y - B_i^\top x \end{pmatrix},
\]
where each \(B_i = B + \sigma U D_i V^\top\) with \(D_i = \diag(\varepsilon_{i,1},\dots ,\varepsilon_{i,d})\) 
and \(\varepsilon_{i,j} \in \{-1,1\}\) such that \(\tfrac{1}{N}\sum_{i=1}^N B_i = B\) and 
\(\sigma = 0.5\). \(B\) is generated by taking singular values \(\sigma_i \in[1/\kappa_B,1]\) and setting \(B = U \diag(\sigma_i)_{i=1}^d V^\top\).
We set $d=30$, \(\lambda = L_1 = 1\) and tune the step-sizes using a coarse logarithmic grid
followed by a refinement over a linear grid around the selected parameters. More information on the tuning procedure is detailed in \Cref{app:details}.

\begin{figure}[htbp]
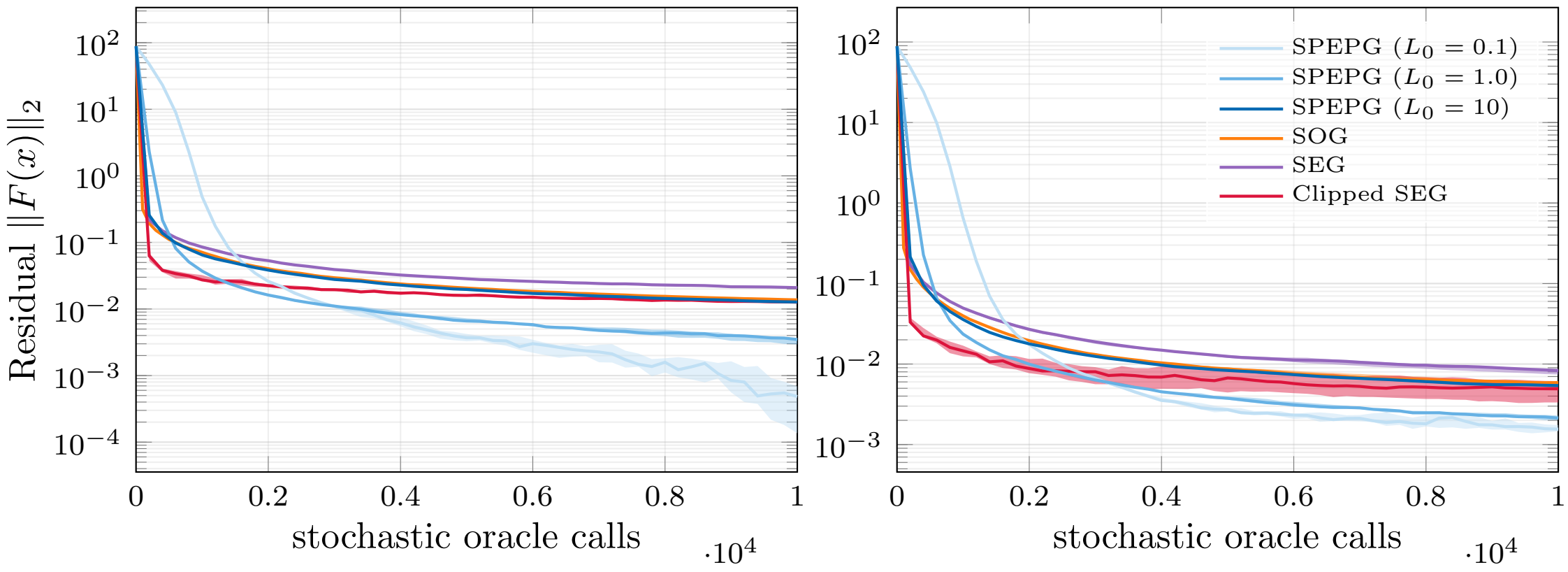

    \centering
    \begin{minipage}[t]{0.48\textwidth}
        \vspace{0pt}
        \centering
        \input{figures/tex/paper_kb10/plot_paper_kb10.tex}
    \end{minipage}
    \begin{minipage}[t]{0.48\textwidth}
        \vspace{0pt}
        \centering
          \input{figures/tex/paper_kb100/plot_paper_kb100.tex}
    \end{minipage}
    \caption{
    Convergence plots for the example in \cref{sec:quartic_nonlinear}.
    The lines and shaded areas represent the median and interquartile range, respectively.
    (\textit{left}) $\kappa(B) = 10$; 
    (\textit{right}) $\kappa(B) = 100$.
    }
    \label{fig:quartic_nonlinear}
\end{figure}

We ran ten instances of the experiment. The results can be found in \Cref{fig:quartic_nonlinear} for \(\kappa_B \in \{10, 100\}\).
In both cases, we observe a similar trade-off as \(L_0\) increases: 
For smaller values of \(L_0\), the method progresses more slowly, but tends to achieve higher accuracy.
As \(L_0\) increases, the effect of the normalization diminishes so the methods behave more similarly to the baselines without preconditioning.

\subsection{Cubic Game} \label{sec:cubic-game}

We consider the two-player game with saddle point operator 
\[
F(x,y)
=
\begin{pmatrix}
\sqrt{\langle x, A x \rangle}\,A x + B y \\
\sqrt{\langle y, C y \rangle}\,C y - B^\top x
\end{pmatrix},
\]
where $A, C \in \R^{d \times d}$ are positive definite matrices,
as described in \cite{choudhury2025extragradient}, and $B \in \R^{d \times d}$ 
is constructed similarly to the previous experiment.
We consider two constructions for the stochastic oracle. Details are provided in \cref{app:details}.

\textbf{Finite sum case.}
Similarly to the example in \cref{sec:quartic_nonlinear},
let
\[
    F_i(x,y)
    =
    \begin{pmatrix}
        \alpha_i \sqrt{x^\top A x}\,Ax + B_i y \\
        \delta_i \sqrt{y^\top C y}\,Cy - B_i^\top x
    \end{pmatrix},
\]
\(i = 1, \dots, N\),
where $\alpha_i, \delta_i$ and $B_i$ are constructed 
randomly to ensure $F(x,y) = \frac{1}{N} \sum_{i=1}^{N} F_i(x,y)$.

\paragraph{Student-\(t\) noise.}
To investigate the effect of heavy-tailed noise, we keep the nonlinear terms
deterministic and perturb only the coupling matrix $B$ for simplicity.
At each oracle evaluation, we sample a matrix \(H\) with independent
standard Student-\(t\) entries,
\(
    H_{jk}\sim t_\nu,
\) 
and replace $B$ by
\(
    \widehat B = B + \frac{\sigma_B}{\sqrt d}\,H,
\)
to obtain the oracle.

\begin{figure}[htb!]
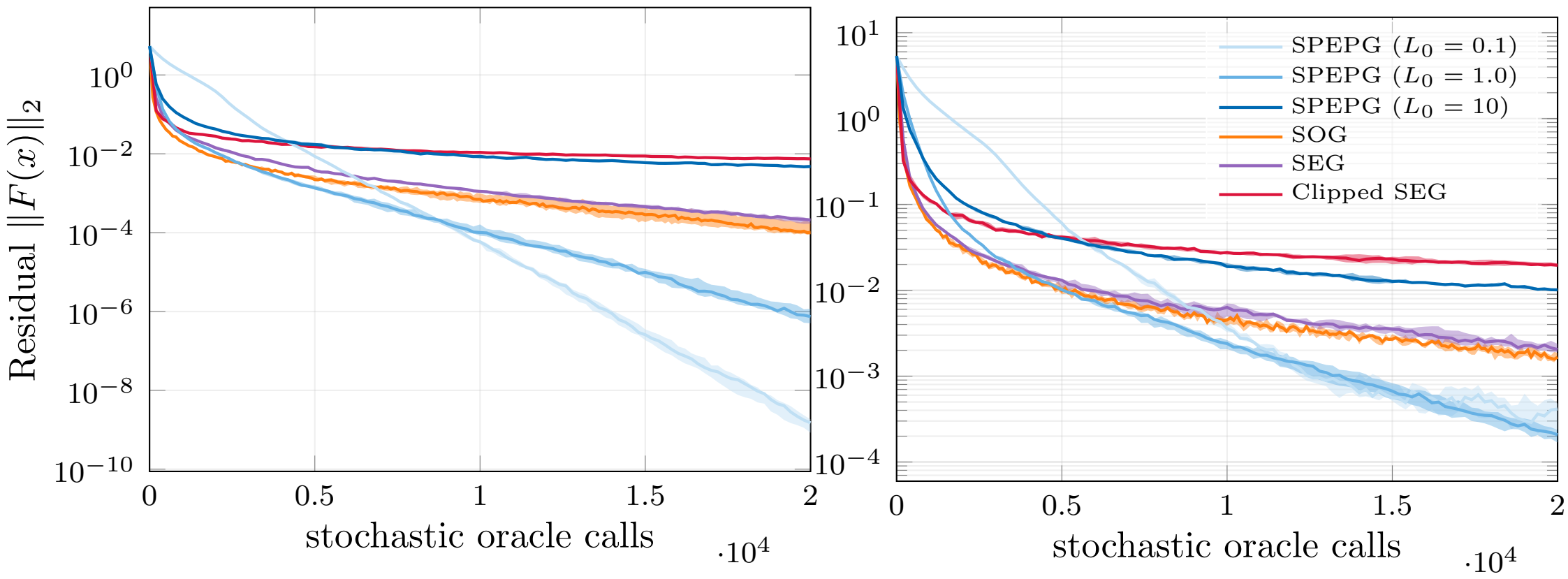

  \centering
  \begin{minipage}[t]{0.48\textwidth}
    \vspace{0pt}
    \centering
    \input{figures/tex/polygame_fs/plot_paper.tex}
  \end{minipage}
  \begin{minipage}[t]{0.48\textwidth}
    \vspace{0pt}
    \centering
    \input{figures/tex/paper_student_t/plot_paper_student_t.tex}
  \end{minipage}
    \caption{
      Convergence plots for the example in \cref{sec:cubic-game}.
      The lines and shaded areas represent the median and interquartile range, respectively.
      (\textit{left}) Finite-sum case;
      (\textit{right}) Student-$t$ noise.
    }
    \label{fig:polygame}
\end{figure}

The results are presented in \cref{fig:polygame}.
As in the previous experiment, we observe a similar trade-off as $L_0$ is 
varied. In particular, it is notable that for smaller values of $L_0$,
our method is able to achieve far greater accuracy than the baselines, even 
in the presence of heavy-tailed noise.

\subsection{Synthetic WGAN}
We evaluate the solvers on a two-dimensional Wasserstein generative adversarial network with gradient penalty (WGAN-GP) \citep{gulrajani2017improved}. The target distribution is a two-dimensional anisotropic Gaussian, \(\gN([2,-1], \Sigma_\star)\) where \(\Sigma_{\star}\) has eigenvalues \(2, 0.2\) and is rotated by \(45\) degrees. The latent variable is two-dimensional and sampled according to \(z \sim \gN(0, 1.5^2 I_2)\). The generator \(G_\phi\) and critic \(D_\psi\) are chosen as two-hidden-layer MLPs with width \(8\) and Leaky-ReLU activations, and the critic initialization scale is chosen in \(s_D \in \{1,4\}\), i.e. the weights are sampled as \(W_{ij} \sim \gN(0,s_D^2/d_{\textrm in})\) with \(d_{\textrm in}\) the input dimension of each linear layer. The form of the optimization problem we are solving is the following:
\[
\min_{\phi}\max_{\psi}\;
L(\phi,\psi),
\]
where
\begin{multline*}
L(\phi,\psi)
=
\mathbb{E}_{x\sim p_{\mathrm{data}}}
\left[D_{\psi}(x)\right]
-
\mathbb{E}_{z\sim p_z}
\left[D_{\psi}(G_{\phi}(z))\right]
\\
-
\lambda_{\mathrm{gp}}
\mathbb{E}_{x,z,\varepsilon}
\left[
\left(
\left\|
\nabla_{\hat{x}}D_{\psi}(\hat{x})
\right\|_2-1
\right)^2
\right]
-
\frac{\lambda_D}{2}\|\psi\|_2^2
+
\frac{\lambda_G}{2}\|\phi\|_2^2.
\end{multline*}
The regularization parameters are set to \(\lambda_{\mathrm{gp}} = 5\), \(\lambda_D = 10^{-4}\) and \(\lambda_G = 10^{-6}\). Moreover, \(\hat{x}=\varepsilon x+(1-\varepsilon)\operatorname{sg}(G_{\phi}(z)), \varepsilon\sim\mathcal{U}[0,1]\) where the symbol \(\operatorname{sg}\) denotes stop-gradient, so the gradient-penalty term contributes only to the critic update and not to the generator update.

We set the budget at \(8000\) batch operator evaluations with batch size \(64\). The results of this experiment can be found in \Cref{tab:gaussian}. We can see that our proposed method outperforms the rest of the methods in the SWD metric, especially in the case \(s_D = 4\) where the difference is more pronounced, while it also achieves lower mean and covariance errors. Both SEG and SOG seem to have high variance between runs, while SPEPG is more stable and more robust to the initialization. The tuning procedure can be found in \Cref{app:details} and the chosen hyperparameters in \Cref{tab:wgangp_hyperparameters}.

\begin{table}[hbpt]
    \centering
    \small
    \caption{
        Results on the 2D Gaussian WGAN-GP problem over 10 runs.
        Lower is better. It can be seen that our proposed method is more stable across runs and more robust towards initialization compared to SEG and SOG.
    }
    \label{tab:gaussian_wgan_results}
    \begin{tabular}{@{}cc ccc@{}}
        \toprule
        $s_D$ & Method
        & Last-5 SWD $\downarrow$
        & Mean error $\downarrow$
        & Covariance error $\downarrow$ \\
        \midrule

        \multirow{3}{*}{$1$}
        & SEG  & $0.3437 \pm 0.2719$ & $0.3662$ & $0.6591$ \\
        & SOG  & $0.3035 \pm 0.1714$ & $0.2772$ & $0.9186$ \\
        & SPEPG & $\mathbf{0.1673 \pm 0.0599}$ & $\mathbf{0.1461}$ & $\mathbf{0.2865}$ \\
        & Clipped SEG & $0.3241 \pm 0.1727$ & $0.5055$ &  $0.8558$ \\
        \midrule

        \multirow{3}{*}{$4$}
        & SEG  & $0.8033 \pm 0.1738$ & $0.9025$ & $4.2742$ \\
        & SOG  & $0.7350 \pm 0.4496$ & $1.1089$ & $2.5993$ \\
        & SPEPG & $\mathbf{0.2676 \pm 0.1481}$ & $\mathbf{0.2354}$ & $\mathbf{0.4026}$ \\
        & Clipped SEG & $0.5481 \pm 0.0533$ & $0.4283$ & $2.4634$ \\
        \bottomrule
    \end{tabular}
    \label{tab:gaussian}
\end{table}

\section{Conclusion}
In this paper we introduced a stochastic independently normalized past extragradient algorithm for problems beyond Lipschitz continuity, assuming merely the Minty variational inequality. We obtained convergence guarantees in expectation and in probability and we showed that the proposed method converges a.s., under bounded variance assumptions for the stochastic oracle. 

There are several open directions for future work. To begin with, analyzing the convergence of the method under more extreme noise assumptions such as heavy-tailed noise is an important research direction that is tightly connected to the normalization-like preconditioning the proposed algorithm performs. Moreover, extending the approach to include coordinate-wise preconditioning for stochastic VI problems is an important research avenue that is still largely underexplored.

\subsection*{Acknowledgments}
This work was supported by the Research Foundation Flanders (FWO) PhD grant 11A8T26N and Postdoc grant 1248626N; Research Council KUL grant C14/24/103. This work was supported in part by a grant of access to OpenAI models through the ChatGPT for Academic Researchers program. The resources and services used in this work were provided by the VSC (Flemish Supercomputer Center), funded by the Research Foundation - Flanders (FWO) and the Flemish Government.

In this work, we used generative AI to aid in implementing and testing the code for the experiments. All implementations were tested and checked for correctness by at least 2 authors.
It was also used to verify and improve mathematical proofs, in particular, to finetune the arguments in order to obtain better constants. 
Finally, it was used to proofread the manuscript and check it for errors and typos.
Generative AI was \textbf{not} used to write or edit the text, generate ideas or develop new theoretical arguments.

\bibliography{arxiv_references}

@inproceedings{alacaoglu2022stochastic,
  title={Stochastic variance reduction for variational inequality methods},
  author={Alacaoglu, Ahmet and Malitsky, Yura},
  booktitle={Conference on Learning Theory},
  pages={778--816},
  year={2022},
}

@InProceedings{alacaoglu2024revisiting,
  title = 	 {Revisiting Inexact Fixed-Point Iterations for Min-Max Problems: Stochasticity and Structured Nonconvexity},
  author =       {Alacaoglu, Ahmet and Kim, Donghwan and Wright, Stephen},
  booktitle = 	 {International Conference on Machine Learning},
  pages = 	 {840--878},
  year = 	 {2024},
}

@article{bertsekas2011incremental,
  title={Incremental proximal methods for large scale convex optimization},
  author={Bertsekas, Dimitri P},
  journal={Mathematical Programming},
  volume={129},
  number={2},
  pages={163--195},
  year={2011},
  publisher={Springer}
}

@inproceedings{bodard2025escaping,
 author = {Bodard, Alexander and Patrinos, Panagiotis},
 booktitle = {Advances in Neural Information Processing Systems},
 pages = {124173--124208},
 title = {Escaping saddle points without {L}ipschitz smoothness: the power of nonlinear preconditioning},
 year = {2025}
}

@article{bohm2022two,
  title={Two Steps at a Time---Taking {GAN} Training in Stride with {T}seng's Method},
  author={Bohm, Axel and Sedlmayer, Michael and Csetnek, Erno Robert and Bot, Radu Ioan},
  journal={SIAM Journal on Mathematics of Data Science},
  volume={4},
  number={2},
  pages={750--771},
  year={2022},
  publisher={SIAM}
}

@incollection{bayraksan_DataDrivenStochasticProgramming_2015,
  title={Data-driven stochastic programming using phi-divergences},
  author={Bayraksan, G{\"u}zin and Love, David K},
  booktitle={The operations research revolution},
  pages={1--19},
  year={2015},
  publisher={Informs}
}

@inproceedings{chavdarova2019reducing,
 author = {Chavdarova, Tatjana and Gidel, Gauthier and Fleuret, Fran\c{c}ois and Lacoste-Julien, Simon},
 booktitle = {Advances in Neural Information Processing Systems},
 title = {Reducing Noise in {GAN} Training with Variance Reduced Extragradient},
 year = {2019}
}

@inproceedings{choudhury2025extragradient,
 author = {Choudhury, Sayantan and Loizou, Nicolas},
 booktitle = {Advances in Neural Information Processing Systems},
 pages = {70738--70781},
 title = {Extragradient Method for {$(L_0, L_1)$}-{L}ipschitz Root-finding Problems},
 year = {2025}
}

@inproceedings{
    daskalakis2017training,
    title={Training {GAN}s with Optimism},
    author={Constantinos Daskalakis and Andrew Ilyas and Vasilis Syrgkanis and Haoyang Zeng},
    booktitle={International Conference on Learning Representations},
    year={2018},
}

@inproceedings{daskalakis2020independent,
 author = {Daskalakis, Constantinos and Foster, Dylan J and Golowich, Noah},
 booktitle = {Advances in Neural Information Processing Systems},
 pages = {5527--5540},
 title = {Independent Policy Gradient Methods for Competitive Reinforcement Learning},
 year = {2020}
}

@article{deng2012mnist,
  title={The {MNIST} database of handwritten digit images for machine learning research [best of the web]},
  author={Deng, Li},
  journal={IEEE Signal Processing Magazine},
  volume={29},
  number={6},
  pages={141--142},
  year={2012},
  publisher={IEEE}
}

@inproceedings{diakonikolas2021efficient,
  title={Efficient methods for structured nonconvex-nonconcave min-max optimization},
  author={Diakonikolas, Jelena and Daskalakis, Constantinos and Jordan, Michael I},
  booktitle={International Conference on Artificial Intelligence and Statistics},
  pages={2746--2754},
  year={2021},
}

@inproceedings{emmanouilidis2024stochastic,
  title={Stochastic extragradient with random reshuffling: Improved convergence for variational inequalities},
  author={Emmanouilidis, Konstantinos and Vidal, Ren{\'e} and Loizou, Nicolas},
  booktitle={International Conference on Artificial Intelligence and Statistics},
  pages={3682--3690},
  year={2024},
}

@book{facchinei2003finite,
  title={Finite-dimensional variational inequalities and complementarity problems},
  author={Facchinei, Francisco and Pang, Jong-Shi},
  year={2003},
  publisher={Springer}
}

@inproceedings{fan2024weaker,
  title={Weaker {MVI} condition: Extragradient methods with multi-step exploration},
  author={Fan, Yifeng and Li, Yongqiang and Chen, Bo},
  booktitle={International Conference on Learning Representations},
  pages={9353--9378},
  year={2024}
}

@inproceedings{faw2023beyond,
  title={Beyond uniform smoothness: A stopped analysis of adaptive {SGD}},
  author={Faw, Matthew and Rout, Litu and Caramanis, Constantine and Shakkottai, Sanjay},
  booktitle={Conference on Learning Theory},
  pages={89--160},
  year={2023},
}

@inproceedings{gaash2026convergence,
 author = {Gaash, Ofir and Levy, Kfir Y. and Carmon, Yair},
 booktitle = {Advances in Neural Information Processing Systems},
 pages = {120410--120442},
 title = {Convergence of Clipped {SGD} on Convex {$(L_0,L_1)$}-Smooth Functions},
 year = {2025}
}

@article{gangbo1996geometry,
  author  = {Gangbo, Wilfrid and McCann, Robert J.},
  title   = {The geometry of optimal transportation},
  journal = {Acta Mathematica},
  volume  = {177},
  number  = {2},
  pages   = {113--161},
  year    = {1996},
}

@inproceedings{
    gidel2018variational,
    title={A Variational Inequality Perspective on Generative Adversarial Networks},
    author={Gauthier Gidel and Hugo Berard and Gaëtan Vignoud and Pascal Vincent and Simon Lacoste-Julien},
    booktitle={International Conference on Learning Representations},
    year={2019},
}

@inproceedings{gorbunov2022last,
 author = {Gorbunov, Eduard and Taylor, Adrien and Gidel, Gauthier},
 booktitle = {Advances in Neural Information Processing Systems},
 pages = {21858--21870},
 title = {Last-Iterate Convergence of Optimistic Gradient Method for Monotone Variational Inequalities},
 year = {2022}
}

@article{goodfellow2020generative,
  title={Generative adversarial networks},
  author={Goodfellow, Ian and Pouget-Abadie, Jean and Mirza, Mehdi and Xu, Bing and Warde-Farley, David and Ozair, Sherjil and Courville, Aaron and Bengio, Yoshua},
  journal={Communications of the ACM},
  volume={63},
  number={11},
  pages={139--144},
  year={2020},
  publisher={ACM New York, NY, USA}
}

@inproceedings{gorbunov2022clipped,
 author = {Gorbunov, Eduard and Danilova, Marina and Dobre, David and Dvurechenskii, Pavel and Gasnikov, Alexander and Gidel, Gauthier},
 booktitle = {Advances in Neural Information Processing Systems},
 pages = {31319--31332},
 title = {Clipped Stochastic Methods for Variational Inequalities with Heavy-Tailed Noise},
 year = {2022}
}

@InProceedings{gorbunov2023high,
  title = 	 {High-Probability Convergence for Composite and Distributed Stochastic Minimization and Variational Inequalities with Heavy-Tailed Noise},
  author =       {Gorbunov, Eduard and Sadiev, Abdurakhmon and Danilova, Marina and Horv\'{a}th, Samuel and Gidel, Gauthier and Dvurechensky, Pavel and Gasnikov, Alexander and Richt\'{a}rik, Peter},
  booktitle = 	 {International Conference on Machine Learning},
  pages = 	 {15951--16070},
  year = 	 {2024},
}

@inproceedings{gorbunov2025methods,
  title={Methods for convex {$(L_0, L_1)$}-smooth optimization: Clipping, acceleration, and adaptivity},
  author={Gorbunov, Eduard and Tupitsa, Nazarii and Choudhury, Sayantan and Aliev, Alen and Richt{\'a}rik, Peter and Horv{\'a}th, Samuel and Tak{\'a}{\v{c}}, Martin},
  booktitle={International Conference on Learning Representations},
  pages={97261--97305},
  year={2025}
}

@inproceedings{gulrajani2017improved,
 author = {Gulrajani, Ishaan and Ahmed, Faruk and Arjovsky, Martin and Dumoulin, Vincent and Courville, Aaron},
 booktitle = {Advances in Neural Information Processing Systems},
 title = {Improved Training of {W}asserstein {GAN}s},
 year = {2017}
}

@article{greenberg2026extragradient,
  title={Extragradient and Past Extragradient Methods for $(L\_0, L\_1) $-Lipschitz Variational Inequalities},
  author={Greenberg, Ezra and Loizou, Nicolas},
  journal={arXiv preprint arXiv:2610.16886},
  year={2026}
}

@article{juditsky2011solving,
  title={Solving variational inequalities with {S}tochastic {M}irror-{P}rox algorithm},
  author={Juditsky, Anatoli and Nemirovski, Arkadi and Tauvel, Claire},
  journal={Stochastic Systems},
  volume={1},
  number={1},
  pages={17--58},
  year={2011},
  publisher={INFORMS}
}

@inproceedings{hsieh2020explore,
 author = {Hsieh, Yu-Guan and Iutzeler, Franck and Malick, J\'{e}r\^{o}me and Mertikopoulos, Panayotis},
 booktitle = {Advances in Neural Information Processing Systems},
 pages = {16223--16234},
 title = {Explore Aggressively, Update Conservatively: Stochastic Extragradient Methods with Variable Stepsize Scaling},
 year = {2020}
}

@article{iusem2017extragradient,
  title={Extragradient method with variance reduction for stochastic variational inequalities},
  author={Iusem, Alfredo N and Jofr{\'e}, Alejandro and Oliveira, Roberto Imbuzeiro and Thompson, Philip},
  journal={SIAM Journal on Optimization},
  volume={27},
  number={2},
  pages={686--724},
  year={2017},
  publisher={SIAM}
}

@inproceedings{kalogiannis2026policy,
 author = {Kalogiannis, Fivos and Farina, Gabriele},
 booktitle = {Advances in Neural Information Processing Systems},
 pages = {29162--29224},
 title = {Policy Gradient Methods Converge Globally in Imperfect-Information Extensive-Form Games},
 year = {2025}
}

@inproceedings{koloskova2023revisiting,
  title={Revisiting Gradient Clipping: Stochastic bias and tight convergence guarantees},
  author={Koloskova, Anastasia and Hendrikx, Hadrien and Stich, Sebastian U},
  booktitle={International Conference on Machine Learning},
  pages={17343--17363},
  year={2023},
}

@inproceedings{
    korotin2022neural,
    title={Neural Optimal Transport},
    author={Alexander Korotin and Daniil Selikhanovych and Evgeny Burnaev},
    booktitle={International Conference on Learning Representations},
    year={2023},
}

@article{korpelevich1977extragradient,
  title={Extragradient method for finding saddle points and other problems},
  author={Korpelevich, Galina M},
  journal={Matekon},
  volume={13},
  number={4},
  pages={35--49},
  year={1977},
  publisher={ME SHARPE INC 80 BUSINESS PARK DR, ARMONK, NY 10504}
}

@article{laude2023dualities,
  title={Dualities for non-{E}uclidean smoothness and strong convexity under the light of generalized conjugacy},
  author={Laude, Emanuel and Themelis, Andreas and Patrinos, Panagiotis},
  journal={SIAM Journal on Optimization},
  volume={33},
  number={4},
  pages={2721--2749},
  year={2023},
  publisher={SIAM}
}

@article{laude2023anisotropic,
  title={Anisotropic proximal point algorithm},
  author={Laude, Emanuel and Patrinos, Panagiotis},
  journal={arXiv preprint arXiv:2312.09834},
  year={2023}
}

@article{laude2025anisotropic,
  author={Laude, Emanuel and Patrinos, Panagiotis},
  title     = {Anisotropic proximal gradient},
  journal   = {Mathematical Programming},
  volume    = {214},
  pages     = {801--845},
  year      = {2025},
  publisher = {Springer},
}

@article{maddison2021dual,
	author = {Maddison, Chris J. and Paulin, Daniel and Teh, Yee Whye and Doucet, Arnaud},
	journal = {SIAM Journal on Optimization},
	number = {1},
	pages = {991--1016},
	publisher = {SIAM},
	title = {Dual space preconditioning for gradient descent},
	volume = {31},
	year = {2021}}

@article{malitsky_golden_2020,
	title = {Golden ratio algorithms for variational inequalities},
	volume = {184},
	language = {en},
	number = {1},
	journal = {Mathematical Programming},
	author = {Malitsky, Yura},
	month = nov,
	year = {2020},
	pages = {383--410},
}

@inproceedings{
    mertikopoulos2018optimistic,
    title={Optimistic mirror descent in saddle-point problems: Going the extra(-gradient) mile},
    author={Panayotis Mertikopoulos and Bruno Lecouat and Houssam Zenati and Chuan-Sheng Foo and Vijay Chandrasekhar and Georgios Piliouras},
    booktitle={International Conference on Learning Representations},
    year={2019},
}

@article{nie2024nash,
  title={Nash equilibrium problems of polynomials},
  author={Nie, Jiawang and Tang, Xindong},
  journal={Mathematics of Operations Research},
  volume={49},
  number={2},
  pages={1065--1090},
  year={2024},
  publisher={informs}
}

@article{namkoong2016stochastic,
  title={Stochastic gradient methods for distributionally robust optimization with f-divergences},
  author={Namkoong, Hongseok and Duchi, John C},
  journal={Advances in neural information processing systems},
  volume={29},
  year={2016}
}

@inproceedings{oikonomidis2025nonlinearlymomentum,
 author = {Oikonomidis, Konstantinos and Quan, Jan and Patrinos, Panagiotis},
 booktitle = {Advances in Neural Information Processing Systems},
 pages = {38957--38988},
 title = {Nonlinearly Preconditioned Gradient Methods: Momentum and Stochastic Analysis},
 year = {2025}
}

@InProceedings{oikonomidis2025nonlinearly,
  title = 	 {Nonlinearly Preconditioned Gradient Methods under Generalized Smoothness},
  author =       {Oikonomidis, Konstantinos and Quan, Jan and Laude, Emanuel and Patrinos, Panagiotis},
  booktitle = 	 {International Conference on Machine Learning},
  pages = 	 {47132--47154},
  year = 	 {2025},
}

@inproceedings{palaniappan2016stochastic,
 author = {Palaniappan, Balamurugan and Bach, Francis},
 booktitle = {Advances in Neural Information Processing Systems},
 title = {Stochastic Variance Reduction Methods for Saddle-Point Problems},
 year = {2016}
}

@inproceedings{
    pethick2023escaping,
    title={Escaping limit cycles: Global convergence for constrained nonconvex-nonconcave minimax problems},
    author={Thomas Pethick and Puya Latafat and Panagiotis Patrinos and Olivier Fercoq and Volkan Cevher},
    booktitle={International Conference on Learning Representations},
    year={2022},
}

@inproceedings{
    pethick2023solving,
    title={Solving stochastic weak Minty variational inequalities without increasing batch size},
    author={Thomas Pethick and Olivier Fercoq and Puya Latafat and Panagiotis Patrinos and Volkan Cevher},
    booktitle={International Conference on Learning Representations},
    year={2023},
}

@inproceedings{pethick2025efficient,
  title={Efficient interpolation between extragradient and proximal methods for weak {MVI}s},
  author={Pethick, Thomas and Mavrothalassitis, Ioannis and Cevher, Volkan},
  booktitle={International Conference on Learning Representations},
  pages={50896--50923},
  year={2025}
}

@inproceedings{pethick2025generalizedgradientnormclipping,
 author = {Pethick, Thomas and Xie, Wanyun and Erdogan, Mete and Antonakopoulos, Kimon and Silveti-Falls, Antonio and Cevher, Volkan},
 booktitle = {Advances in Neural Information Processing Systems},
 pages = {21170--21208},
 title = {Generalized Gradient Norm Clipping \& Non-{E}uclidean {$(L_0,L_1)$}-Smoothness},
 year = {2025}
}

@article{popov1980modification,
  title={A modification of the {A}rrow--{H}urwicz method for search of saddle points},
  author={Popov, Leonid Denisovich},
  journal={Mathematical notes of the Academy of Sciences of the USSR},
  volume={28},
  number={5},
  pages={845--848},
  year={1980},
  publisher={Springer}
}

@article{quan2026dual,
  title={Dual-Space Preconditioning for Variational Inequalities and Root-Finding Problems},
  author={Quan, Jan and Oikonomidis, Konstantinos and Bodard, Alexander and Patrinos, Panagiotis},
  journal={arXiv preprint arXiv:2611.13916},
  year={2026}
}

@InProceedings{riabinin2025gluon,
  title = 	 {From Muon to Gluon: Bridging Theory and Practice of {LMO}-based Optimizers for {LLM}s},
  author =       {Riabinin, Artem and Shulgin, Egor and Gruntkowska, Kaja and Richt\'{a}rik, Peter},
  booktitle = 	 {International Conference on Machine Learning},
  pages = 	 {105039--105082},
  year = 	 {2026},
}

@InProceedings{vankov2024generalized,
  title = 	 {Generalized Smooth Variational Inequalities: Methods with Adaptive Stepsizes},
  author =       {Vankov, Daniil and Nedich, Angelia and Sankar, Lalitha},
  booktitle = 	 {International Conference on Machine Learning},
  pages = 	 {49137--49170},
  year = 	 {2024},
}

@article{
    vankov2025generalized,
    title={Generalized Smooth Stochastic Variational Inequalities:  Almost Sure Convergence and Convergence Rates},
    author={Daniil Vankov and Angelia Nedich and Lalitha Sankar},
    journal={Transactions on Machine Learning Research},
    issn={2835-8856},
    year={2025},
}

@inproceedings{vankov2025optimizing,
  title={Optimizing {$(L_0, L_1)$}-Smooth Functions by Gradient Methods},
  author={Vankov, Daniil and Rodomanov, Anton and Nedich, Angelia and Sankar, Lalitha and Stich, Sebastian},
  booktitle={International Conference on Learning Representations},
  pages={15953--15979},
  year={2025}
}

@inproceedings{wang2023convergence,
  title={Convergence of {A}dagrad for non-convex objectives: Simple proofs and relaxed assumptions},
  author={Wang, Bohan and Zhang, Huishuai and Ma, Zhiming and Chen, Wei},
  booktitle={Conference on Learning Theory},
  pages={161--190},
  year={2023},
}

@book{williams1991probability,
  title={Probability with martingales},
  author={Williams, David},
  year={1991},
  publisher={Cambridge University Press}
}

@inproceedings{
    zhang2019gradient,
    title={Why Gradient Clipping Accelerates Training: A Theoretical Justification for Adaptivity},
    author={Jingzhao Zhang and Tianxing He and Suvrit Sra and Ali Jadbabaie},
    booktitle={International Conference on Learning Representations},
    year={2020},
}

@inproceedings{zhang2020improved,
 author = {Zhang, Bohang and Jin, Jikai and Fang, Cong and Wang, Liwei},
 booktitle = {Advances in Neural Information Processing Systems},
 pages = {15511--15521},
 title = {Improved Analysis of Clipping Algorithms for Non-convex Optimization},
 year = {2020}
}
\bibliographystyle{ref_style}

\appendix
\newpage
\section{Preconditioned Lipschitz continuity}
\label{app:prec_lipschitz}
\subsection{Preconditioned Lipschitz continuity and \texorpdfstring{\((L_0,L_1)\)-Lipschitz continuity}{(L0, L1)-Lipschitz continuity}}
The regularity condition we work with is related to \((L_0,L_1)\)-Lipschitz continuity but is in fact a less restrictive condition. In this section we explore the differences between the two notions and provide examples of operators that satisfy \Cref{def:prec_lipschitz}. To begin with, we state the definition of \((L_0,L_1)\)-Lipschitz continuity as adapted from \citet[Assumption 2.1]{vankov2024generalized} (with \(\alpha = 1\)), see also \citet[Assumption 1.1]{choudhury2025extragradient}. We focus on the case \(L_0 > 0\).
\begin{definition} \label{def:l0l1}
    We say that a continuous \(F:\mathbb{R}^n\to\mathbb{R}^n\) is \((L_0,L_1)\)-Lipschitz if for some \(L_0 > 0,L_1 \geq 0\), the following inequality holds:
    \begin{equation} \label{eq:l0l1}
        \|F(x)-F(\bar x)\| \leq \Bigl(L_0 + L_1\max_{\theta \in [0,1]}\|F(\theta x+(1-\theta)\bar x)\|\Bigr)\|x-\bar x\| \qquad \forall x,\bar x \in \mathbb{R}^n.
    \end{equation}
\end{definition}
On the other hand, the preconditioned Lipschitz condition we consider in \Cref{def:prec_lipschitz} imposes $1$-Lipschitz continuity of the operator $T(x) = \frac{F(x)}{L_0 + L_1\|F(x)\|}$, that is,
\[
    \left \| \frac{F(x)}{L_0 + L_1\|F(x)\|}-\frac{F(\bar x)}{L_0+L_1\|F(\bar x)\|}\right \| \le \|x-\bar x\| \qquad \forall x,\bar x \in \mathbb{R}^n.
\]
It is shown in \citet[Proposition 3.4]{quan2026dual} that every \((L_0,L_1)\)-Lipschitz operator is also preconditioned Lipschitz. We provide a self-contained proof here.
To this end, we first show a preliminary result that will be used in the proof.
\begin{lemma} \label{lem:prec_bound}
    Let $y,\bar y \in \mathbb{R}^n$, $L_0 > 0$, $L_1\geq0$. Then,
    \[
        \left\|\frac{y}{L_0 + L_1\|y\|} - \frac{\bar y}{L_0 + L_1\|\bar y\|}\right\| \leq \frac{\|y - \bar y\|}{\sqrt{(L_0 + L_1\|y\|)(L_0 + L_1\|\bar y\|) }}.
    \]  
\end{lemma}
\begin{proof}
     Denote $Y = L_0 + L_1\|y\|$ and $\bar Y = L_0 + L_1\|\bar y\|$. It suffices to show that $\|\bar Y y - Y\bar y\|^2 \leq Y\bar Y\|y - \bar y\|^2$. Indeed, we have
     \begin{align*}
         Y\bar Y \|y - \bar y\|^2 - \|\bar Y y - Y \bar y\|^2 &= (Y - \bar Y) (\bar Y \|y\|^2 - Y \|\bar y\|^2) \\
         &= L_1(\|y\| - \|\bar y\|)^2 (L_0(\|y\| + \|\bar y\|) + L_1\|y\|\|\bar y\|) \geq 0,
     \end{align*}
     hence $\|\bar Y y - Y \bar y\| \leq \sqrt{Y\bar Y} \|y - \bar y\|$, which yields the desired result after dividing by $Y\bar Y\neq 0$.
\end{proof}

We now move on to proving the claim.
\begin{proposition} \label{thm:prec_lip_greater}
    Let \(F:\mathbb{R}^n \to \mathbb{R}^n\) satisfy \Cref{def:l0l1}. Then \(F\) satisfies \Cref{def:prec_lipschitz}.
\end{proposition}
\begin{proof}
    Let $x,\bar x\in \mathbb{R}^n$ and denote $\gamma(t) = (1-t)x + t\bar x$ for $t\in[0,1]$. Since $F$ is continuous, $L_0 + L_1\|F(\gamma(t))\|$ is uniformly continuous on $[0,1]$. Now let $\varepsilon>0$ be arbitrary, and choose $N$ sufficiently large, so that on each interval $[i/N, (i+1)/N]$ for $i=0,\ldots,N-1$, we have
    \begin{align*}
        \max_{t\in[i/N, (i+1)/N]} (L_0 + L_1\|F(\gamma(t))\|) &\leq \min_{t\in[i/N,(i+1)/N]} (L_0 + L_1\|F(\gamma(t))\|) + \varepsilon L_0 \\
        &\leq (1+\varepsilon)\min_{t\in[i/N,(i+1)/N]} (L_0 + L_1\|F(\gamma(t))\|).
    \end{align*}
    Let $x_i = \gamma(i/N)$ denote the endpoints of the intervals. Then, by \Cref{lem:prec_bound}, we have
    \begin{align*}
        \|T(x_{i+1}) - T(x_i)\| &\leq \frac{\|F(x_{i+1}) - F(x_i)\|}{\sqrt{(L_0 + L_1 \|F(x_{i+1})\|)(L_0 + L_1\|F(x_i)\|)}} \\
        &\leq \frac{\max_{t\in[i/N, (i+1)/N]} (L_0 + L_1\|F(\gamma(t))\|)}{\sqrt{(L_0 + L_1 \|F(x_{i+1})\|)(L_0 + L_1\|F(x_i)\|)}} \|x_{i+1} - x_i\| \\
        &\leq \frac{\max_{t\in[i/N, (i+1)/N]} (L_0 + L_1\|F(\gamma(t))\|)}{\min_{t\in[i/N, (i+1)/N]} (L_0 + L_1\|F(\gamma(t))\|)} \|x_{i+1} - x_i\| \leq \frac{1+\varepsilon}{N}\|x - \bar x\|
    \end{align*}
    where we used \Cref{def:l0l1} in the second inequality.
    
    Summing over $i=0,1,\ldots,N-1$ and using the triangle inequality yields
    \[
        \|T(x) - T(\bar x)\| \leq \sum_{i=0}^{N-1} \|T(x_{i+1}) - T(x_i)\| \leq (1 + \varepsilon)\|x - \bar x\|.
    \]
    As $\varepsilon > 0$ is arbitrary, the desired result follows.
\end{proof}

Having established \Cref{thm:prec_lip_greater} the claim about the class of preconditioned Lipschitz operators being larger follows from \citet[Theorem A.3]{bodard2025escaping}, see also \Cref{ex:sin,ex:prec_greater}.

\subsection{Examples}
In this section we present examples of operators that satisfy \Cref{def:prec_lipschitz} but need not be globally Lipschitz continuous. A list of examples can also be found in \citet[Appendix C]{choudhury2025extragradient}. We start with a minimax problem where the interaction between the two players is higher-order nonlinear.
\begin{example}
    Let $\lambda\in\mathbb{R},x\in\mathbb{R}^n,y\in\mathbb{R}^m$ and $B\in\mathbb{R}^{n\times m}$. Consider \(\gL(x,y) = \tfrac{1}{4}\|x\|^4-\tfrac{1}{4}\|y\|^4 + \tfrac{\lambda }{2}\|x\|^2\|y\|^2 + \langle x,By \rangle\) which leads to the saddle operator:
    \[
        F(x, y) = \begin{pmatrix} (\|x\|^2 + \lambda \|y\|^2)x + B y \\ (\|y\|^2 - \lambda \|x\|^2)y - B^\top x \end{pmatrix}.
    \]
    \(F\) satisfies \Cref{def:prec_lipschitz} with constants \(L_0=\|B\|+3+2|\lambda|, L_1 = 6+4|\lambda|\).    
\end{example}
\begin{proof}
    Consider \((x,y) \neq (0,0)\). To begin with we lower bound \(\|F(x,y)\|\). From Cauchy--Schwarz,
    \[
        \|F(x,y)\| \geq \langle F(x,y),(x,y)\rangle/\|(x,y)\| = (\|x\|^4 + \|y\|^4) / \|(x, y)\| \geq (\|x\|^2+\|y\|^2)^{3/2}/2
    \]
    where we used \eqref{eq:young} with $a=\|x\|^2, b=\|y\|^2$ and $e=1$ in the last inequality. The inequality above holds also at the origin by inspection. At the same time, for the Jacobian of \(F\) we have:
    \[
        \nabla F(x, y) = \begin{pmatrix} (\|x\|^2+\lambda\|y\|^2)I+2xx^\top & 2\lambda xy^\top + B \\ -2\lambda yx^\top - B^\top  & (\|y\|^2-\lambda \|x\|^2)I + 2yy^\top \end{pmatrix}.
    \]
    We may decompose this matrix as
    \[
        \nabla F(x,y) = \begin{pmatrix}
            0 & B \\ -B^\top & 0
        \end{pmatrix} + \begin{pmatrix} \|x\|^2I+2xx^\top & 0 \\ 0  & \|y\|^2I + 2yy^\top \end{pmatrix} + \begin{pmatrix}
            \lambda \|y\|^2 I & 2\lambda xy^\top \\ -2\lambda yx^\top & -\lambda\|x\|^2I
        \end{pmatrix}
    \]
    where the spectral norms of the components can be shown to be bounded by $\|B\|$, $3(\|x\|^2+\|y\|^2)$ and $2|\lambda|(\|x\|^2+\|y\|^2)$ respectively. Therefore,
    \[
        \|\nabla F(x, y)\| \leq \|B\| + (3 + 2|\lambda|)(\|x\|^2 + \|y\|^2) \leq (\|B\| + 3 + 2|\lambda|) + (6 + 4|\lambda|)\|F(x,y)\|
    \]
    where we also used the standard bound \(t^2 \leq 1 + t^3\) for \(t = (\|x\|^2 + \|y\|^2)^{1/2} \ge 0\). The result follows from \eqref{eq:l0l1_snd_order} and \Cref{thm:prec_lip_greater}.
\end{proof}
The next example goes beyond polynomial growth of the operator, involving the hyperbolic sine.
\begin{example}
    Consider
    \[
        F(x) = B^\top \sinh(Bx - b),
    \]    
    where \(B \in \mathbb R^{m \times n}\) has full row rank, \(b \in \mathbb R^m, x \in \mathbb R^n\) and $\sinh$ is considered componentwise. This operator is the gradient of \(f(x) = \sum_{i=1}^m \cosh((Bx-b)_i)\), so $F$ is monotone as the gradient of a convex function. Then \(F\) satisfies \Cref{def:prec_lipschitz} with constants \(L_0=\|B\|^2\) and \(L_1 = \tfrac{\|B\|^2}{\sigma_{\rm min}(B)}\), where \(\sigma_{\rm min}(B) \coloneqq \min_{u\in\mathbb{R}^m, \|u\| = 1} \|B^\top u\| > 0\).
\end{example}
\begin{proof}
    We have \(\nabla F(x) = B^\top \diag(\cosh(Bx-b))B\). Using submultiplicativity of the spectral norm, we get
    \[  
        \|\nabla F(x)\| \leq \|B\|^2 \|\diag(\cosh(Bx-b))\| \leq \|B\|^2(1 + \|\sinh(Bx-b)\|) 
    \]
    where we also used $\cosh(t) = \sqrt{1 + \sinh^2(t)} \leq 1 + |\sinh(t)|$. By definition of \(\sigma_{\rm min}(B)\), \(\|F(x)\| \geq \sigma_{\rm min}(B)\|\sinh(Bx-b)\|\) and therefore
    \[
        \|\nabla F(x)\| \leq \|B\|^2 (1 + \|\sinh(Bx-b)\|) \leq \|B\|^2 + \frac{\|B\|^2}{\sigma_{\rm min}(B)} \|F(x)\|,
    \]
    leading to the claimed result in light of \eqref{eq:l0l1_snd_order} and \Cref{thm:prec_lip_greater}.
\end{proof}

Another setting where these problems arise naturally is in distributionally robust optimization (DRO).
In particular, for a DRO problem with a $\chi^2$-ambiguity set \citep{namkoong2016stochastic,bayraksan_DataDrivenStochasticProgramming_2015},
replacing the ambiguity set constraint with a soft penalty leads to a problem
of the following form.

\begin{example}[Higher-order distributionally robust optimization]
    Consider the following distributionally robust optimization problem
    \[
        \min_{x \in \mathbb R^n}\max_{y \in \Delta^m}\gL(x,y)=\sum_{i=1}^m y_i\frac{(\langle a_i,x \rangle - b_i)^{4}}{4} - \frac{\mu}{2}\|y\|^2,
    \]
    where $b\in\mathbb{R}^m$, \(a_i^\top\) are the rows of the data matrix \(A \in \mathbb R^{m \times n}\), \(\Delta^m\) is the probability simplex and \(\mu>0\) is a regularization parameter. The corresponding saddle operator is 
    \[
        F(x, y) = \begin{pmatrix} \sum_{i=1}^m y_i(\langle a_i,x \rangle - b_i)^{3}a_i \\ -l(x)+\mu y\end{pmatrix},
    \]
    where \(l:\mathbb R^n \to \mathbb R^m\) is given by \(l(x) = (\frac{(\langle a_i,x \rangle - b_i)^{4}}{4})_{i=1}^m\). Then, \(F\) satisfies \Cref{def:prec_lipschitz} on $\mathbb{R}^n \times \Delta^m$ with \(L_0 = \mu + 6(\max_i \|a_i\|)^2(2+\mu^{3/4}) + 4^{3/4}(\sum_{i=1}^m \|a_i\|^8)^{1/8}(1 + \mu^{3/4})\) and \(L_1 = 6(\max_i \|a_i\|)^2+4^{3/4}(\sum_{i=1}^m \|a_i\|^8)^{1/8}\).
\end{example}
\begin{proof}
    Denote $r(x) \coloneqq (\langle a_i,x \rangle - b_i)_{i=1}^m$. Then, 
    \[
        \nabla F(x,y) = \begin{pmatrix} H(x,y) & B(x) \\ -B(x)^\top & \mu I_m\end{pmatrix}
    \]
    with $H(x,y) = 3\sum_{i=1}^m y_i r_i(x)^2 a_ia_i^\top$ and $B(x) = [r_1(x)^3 a_1 \dots r_m(x)^3  a_m]$. We now proceed by bounding the operator norm of each submatrix in \(\nabla F(x,y)\). Since $y\in \Delta^m$, it follows that $\|H(x, y)\| \leq 3 \big( \max_i \|a_i\| \big)^2 \|r(x)\|_\infty^2$. Moreover, for \(v \in \mathbb R^n\) with \(\|v\|=1\), we have
    \[
        \|B(x)^\top v\|^2 = \sum_{i=1}^m|r_i(x)|^6|\langle a_i,v \rangle|^2 \le \left ( \sum_{i=1}^m|r_i(x)|^8\right )^{3/4}\left ( \sum_{i=1}^m |\langle a_i,v \rangle|^8 \right )^{1/4}
    \]
    where we used H\"older's inequality. Note that $\|l(x)\|^2 = \frac{1}{16} \sum_{i=1}^m |r_i(x)|^8$, hence taking square roots and then the supremum over $\|v\| = 1$, we obtain
    \[
        \|B(x)\| \leq 4^{3/4}\left ( \sum_{i=1}^m \|a_i\|^8 \right )^{1/8}\|l(x)\|^{3/4}.
    \]    
    From the triangle inequality on the decomposition of \(\nabla F\) we have 
    \[
        \|\nabla F(x,y)\| \le \max\{\|H(x,y)\|, \mu\} + \|B(x)\|.
    \]
    We now aim to associate the obtained bounds with \(\|F(x,y)\|\). Note that \(l(x) = \mu y - F_y(x,y)\) and thus \(\|l(x)\| \leq \mu + \|F(x,y)\|\) since \(\|y\| \le 1\) due to the simplex constraint. Moreover, \(|r_i(x)|^4 \leq 4(\mu + \|F(x,y)\|)\) for all \(i \in \{1,\dots, m\}\) by the same reasoning, leading to \(\|r(x)\|_{\infty}^2 \leq 2\sqrt{\|F(x,y)\|+\mu}\). We are thus led to
    \[
        \|\nabla F(x,y)\| \le \mu + 6\big(\max_{i} \|a_i\|\big)^2\sqrt{\|F(x,y)\| + \mu} + 4^{3/4}\Bigg (\sum_{i=1}^m \|a_i\|^8 \Bigg )^{1/8}(\|F(x,y)\|+\mu)^{3/4}.
    \]
    Using the standard bounds \(t^{1/2} \leq 1 + t^{3/4}\) and \((c+d)^{3/4} \le c^{3/4}+d^{3/4}\) for nonnegative \(c,d\) we get
    \begin{multline*}
        \|\nabla F(x,y)\| 
        \le \mu + 6\big(\max_i \|a_i\|\big)^2(1+\mu^{3/4}) + 4^{3/4}\Bigg(\sum_{i=1}^m \|a_i\|^8\Bigg)^{1/8}\mu^{3/4}
        \\
        + \left[6(\max_i \|a_i\|)^2+4^{3/4}\Bigg(\sum_{i=1}^m \|a_i\|^8\Bigg)^{1/8}\right]\|F(x,y)\|^{3/4}.
    \end{multline*}
    The claimed result now follows from \(t^{3/4} \leq 1+t\).
\end{proof}

In what follows, we focus mostly on multiplayer games. We denote by \(p\) the number of players, where each player \(i\) controls \(x_i \in \mathbb R^{n_i}\), and concatenate the actions in the variable \(x = (x_1, x_2,\dots, x_p) \in \mathbb R^n\) with \(n = \sum_{i=1}^p n_i\). Each player has cost \(f_i(x_i,x_{-i})\) where \(x_{-i}\) encodes the actions of the rest of the players. The game mapping is given by 
\[
    F(x) = \begin{pmatrix} \nabla_{x_1}f_1(x_1,x_{-1}) \\ \nabla_{x_2}f_2(x_2,x_{-2}) \\ \vdots \\ \nabla_{x_p}f_p(x_p,x_{-p})\end{pmatrix}
\]
and in the unconstrained setting under additional convexity assumptions, the first-order equilibrium conditions are \(F(x^\star) = 0\).

\begin{example}[Three-player polynomial game]
    Consider \citet[Example 5.8]{nie2024nash} where a polynomial game with three players is described and assume that all \(x_i \in \mathbb R^d\) such that \(n = 3d\) for simplicity. Each player has cost
    \[
        f_i(x) = \sum_{j=1}^d x_{i,j}^4 + \frac{1}{d^2}\sum_{0 \leq m\leq j \leq k \leq d}x_{i,m}x_{i,j}(x_{i,k} + x_{i+1,m} + x_{i+2,j}),
    \]
    where \(x_{i,0} = 1\) and we interpret player indices modulo three. Note that for \(1\leq l\leq d\), \(F_{i,l}(x) = 4x_{i,l}^3 + p_{i,l}(x)\) where \(p_{i,l}\) is a polynomial of degree at most \(2\). In this case, for any given $L_1>0$, there exists a finite $L_0>0$ such that  \(F\) satisfies \Cref{def:prec_lipschitz}.
\end{example}
\begin{proof}
    We use the proof technique of \cite{bodard2025escaping} and aim to bound the Jacobian by polynomials of degree at most two and show that $\|F(x)\|$ has cubic growth up to lower-order terms. Let \(H_3\) be the leading term in \(F\), i.e. \(H_3(x) = 4x^{[3]}\), where \(x^{[3]}\) denotes componentwise cubing. We have
    \[
        \|H_3(x)\| = 4\|x\|_6^3 \geq \frac{4}{n}\|x\|^3.
    \]
    Moreover, $\nabla H_3(x) = 12\diag(x_j^2)_{j=1}^n$ so that $\|\nabla H_3(x)\|\leq 12\|x\|^2$. Regarding the remainder $P_2(x) \coloneqq F(x) - H_3(x)$, we note that it is a polynomial of degree at most $2$ and thus \(\|P_2(x)\| \leq C_0 + C_1\|x\| + C_2\|x\|^2\) and \(\|\nabla P_2(x)\| \leq C_1 + C_2\|x\|\) for some \(C_0,C_1,C_2 \ge 0\). Therefore, $\|\nabla F(x)\| \leq 12\|x\|^2 + C_2\|x\| + C_1$ by the triangle inequality. At the same time, the reverse triangle inequality ensures that $\|F(x)\| \ge \frac{4}{n}\|x\|^3 - C_2\|x\|^2 - C_1\|x\| - C_0$. For fixed $L_1>0$, we are led to
    \[
        \|\nabla F(x)\| - L_1\|F(x)\|\leq -\frac{4L_1}{n}\|x\|^3 + (12+L_1C_2)\|x\|^2 + (C_2+L_1C_1)\|x\| + C_1 + L_1C_0.
    \]
    For \(\|x\| \to \infty\), the r.h.s.\ tends to \(-\infty\) and thus admits a finite supremum. Choosing \(L_0\) larger than this supremum concludes the proof.
\end{proof}

\begin{example}[A restriction of the Nash--Cournot game]
    Consider the Nash--Cournot oligopolistic equilibrium model as described in \citet{malitsky_golden_2020}. We explicitly restrict the nonnegative orthant constraints to the form \(C = \prod_{i=1}^p[l_i, +\infty)\) with \(l_i > 0\) so that the operator \(F\) is well-defined. We consider the equilibrium operator \(F(x) = (F_i(x))_{i=1}^p\) where each coordinate is given by 
    \[
        F_i(x) = f_i'(x_i) - \pi\Bigg(\sum_{j=1}^p x_j\Bigg) - x_i \pi'\Bigg(\sum_{j=1}^p x_j\Bigg), \qquad i=1,\ldots,p.
    \]
    The inverse demand function \(\pi\) and each firm's cost function \(f_i\) are given by 
    \[
        \pi(q) = a^{\rho}q^{-\rho} \quad \textnormal{and} \quad f_i(x_i) = c_i x_i + \frac{d_i}{\nu_i+1}x_i^{\nu_i+1},
    \]
    where \(a,\nu_i,d_i, \rho>0\) and \(c_i \ge 0\). In this setting, \(F\) satisfies \Cref{def:prec_lipschitz} on \(C\), with constants \(L_0 = \rho a^\rho (\sum_{i=1}^p l_i)^{-\rho-1} [p + 1 + (\rho + 1) \sqrt{p}] + a^\rho (\max_i \nu_i/l_i) (\sum_{i=1}^p l_i)^{-\rho}\) and \(L_1 = \max_i\nu_i/l_i\).
\end{example}
\begin{proof}
    We work in the constraint set \(C\) where all the functions involved are differentiable. We proceed in a similar manner to the previous proofs by bounding the higher-order terms in the Jacobian in terms of $\|F(x)\|$. To begin with, denoting $q = \sum_{j=1}^px_j$, we have $F_i(x) = c_i + d_i x_i^{\nu_i} - a^\rho q^{-\rho} + \rho a^{\rho}x_i q^{-\rho-1}$. Therefore, we can write $F = c + D + H$ where $D_i(x) = d_i x_i^{\nu_i}$ and $H_i(x) = -a^\rho q^{-\rho} (1 - \rho x_i / q)$.
    
    Since $x_i\geq l_i$, we have
    \[
        H_i(x) \geq -a^\rho \Bigg(\sum_{j=1}^p l_j\Bigg)^{-\rho}, \qquad D_i(x) \leq |F_i(x)| + a^\rho \Bigg(\sum_{j=1}^p l_j\Bigg)^{-\rho}  \leq \|F(x)\| + a^\rho \Bigg(\sum_{j=1}^p l_j\Bigg)^{-\rho},
    \]
    where we also used $c_i\geq0$. It follows that 
    \[
        \|\nabla D(x)\| = \max_i \frac{\nu_i}{x_i} D_i(x) \leq \left( \max_i \frac{\nu_i}{l_i}\right) \Bigg(\|F(x)\| + a^\rho \Bigg(\sum_{j=1}^p l_j\Bigg)^{-\rho} \Bigg).
    \]
    Moreover,
    \[
        \|\nabla H(x)\| \leq \rho a^\rho q^{-\rho - 1}\left[\|I + \mathbf{1}\mathbf{1}^\top\| + (\rho + 1)\frac{\|x\|}{q}\|\mathbf{1}\|\right] \leq \rho a^\rho \Bigg(\sum_{j=1}^p l_j\Bigg)^{-\rho - 1} [ p +1 + (\rho + 1)\sqrt{p}],
    \]
    where $\mathbf{1}$ is the vector of all ones and $I$ denotes the identity matrix. Here, we used that $\|I + \mathbf{1}\mathbf{1}^\top\| = p + 1$, $\|x\| \leq q$, $\|\mathbf{1}\| = \sqrt{p}$ and $q \geq \sum_{j=1}^p l_j$. Therefore, $\|\nabla F(x)\| \leq \|\nabla D(x)\| + \|\nabla H(x)\| \leq L_0 + L_1 \|F(x)\|$ with the stated $L_0,L_1$ yielding the claimed result.
\end{proof}

Finally, we move on to provide another example that satisfies \Cref{def:prec_lipschitz} but fails to be \((L_0,L_1)\)-Lipschitz. To this end, we first prove the following helper result.
\begin{lemma} \label{thm:radial}
    Let $J\in\mathbb{R}^{n\times n}$ satisfy $J^\top J = I$, and let $g : [0, \infty) \to [0, \infty)$ satisfy $g(0) = 0$ and $g\in \mathcal{C}^1((0, \infty))$. Define $F : \mathbb{R}^n \to \mathbb{R}^n$ by $F(x) = g(\|x\|) \frac{Jx}{\|x\|}$ for $x\neq 0$ and $0$ otherwise. Let $L_0>0$, $L_1\geq0$, and define $h(r) = \frac{g(r)}{L_0 + L_1g(r)}$ for $r>0$. If $|h'(r)|\leq1$ and $\frac{|h(r)|}{r}\leq1$ for all $r>0$, then $F$ satisfies \cref{def:prec_lipschitz} with constants $L_0,L_1$.
\end{lemma}
\begin{proof}
    For $x\neq0$, let $r=\|x\|$. Since $J^\top J = I$, $\|F(x)\| = g(r)$, and hence
    \[
        T(x) = h(r) \frac{Jx}{r}, \qquad \nabla T(x) = J\left[h'(r)\frac{xx^\top}{r^2} + \frac{h(r)}{r} \left(I - \frac{xx^\top}{r^2}\right) \right].
    \]
    The two matrices above are orthogonal projections and $J^\top J = I$, so we have that $\|\nabla T(x)\| \leq \max \left\{|h'(r)|, \frac{|h(r)|}{r}\right\} \leq 1$. Moreover, $\|T(x)\| = |h(\|x\|)| \leq \|x\|$, so $T$ is continuous at the origin. It follows that $T$ is $1$-Lipschitz continuous, from which the desired result follows.
\end{proof}

\begin{example} \label{ex:prec_greater}
    Let $J\in\mathbb{R}^{n\times n}$ satisfy $J^\top J = I$, and define the super-exponential operator $F : \mathbb{R}^n \to \mathbb{R}^n$ by 
    \[
        F(x) = \frac{\exp(\|x\|^2)-1}{\|x\|}Jx, \qquad x \neq 0
    \]
    and \(F(0) = 0\). Then \(F\) satisfies \Cref{def:prec_lipschitz} with constants \(L_0 = L_1 = 1\). Moreover, \(F\) is not \((L_0,L_1)\)-Lipschitz for any $L_0>0$ and \(L_1\geq0\).
\end{example}
\begin{proof}
    Let $g(r) = e^{r^2} - 1$. For $L_0 = L_1 = 1$, $h(r) = \frac{g(r)}{1 + g(r)} = 1 - e^{-r^2}$, and, using $1 - e^{-r^2} \leq \min\{r^2, 1\}$, we have $\frac{h(r)}{r} \leq \min\{r, 1/r\} \leq 1$ and $|h'(r)| = 2re^{-r^2} \leq \sqrt{\frac{2}{e}} < 1$. Hence, the first claim follows from \cref{thm:radial}.

    For the second claim, suppose that $F$ is $(L_0, L_1)$-Lipschitz for some $L_0 > 0$ and $L_1\geq 0$. Fix a unit vector $u$ and apply \Cref{def:l0l1} to $x = (r + t)u$ and $\bar x = ru$. Since $g$ is increasing, $\frac{g(r+t) - g(r)}{t} \leq L_0 + L_1g(r+t)$. Letting $t\downarrow0$ and using $g(r) = e^{r^2} - 1$ gives
    \[
        2re^{r^2} \leq L_0 + L_1 (e^{r^2} - 1) \implies 2r \leq L_0e^{-r^2} + L_1(1 - e^{-r^2}) \leq L_0 + L_1
    \]
    which is impossible as $r\to\infty$. Hence $F$ is not $(L_0,L_1)$-Lipschitz for any $L_0>0$ and $L_1\geq0$.
\end{proof}

\section{Analysis of stochastic past extra preconditioned gradient} \label{app:analysis}
Throughout this appendix, we analyze the past-extragradient method described in \Cref{alg:polyak-sep} under \Cref{assum:bounded_variance,assum:prec_lip_minty} and the parameter choices in \Cref{assum:param_choice}. The main difficulty is that the conditional mean of the oracle is a positive, iteration-dependent rescaling of $T$,
cf.\ \Cref{thm:stoch_bounds}, and the momentum estimator averages these oracle evaluations at past points.
To address this challenge, we introduce an expected mean tracker and bound the momentum tracking error by an averaged noise term and a history-dependent term in \Cref{thm:bound_error}. \Cref{prop:dist_est,thm:bound_rk} establish a one-step descent estimate for the distance to the solution, involving this momentum tracking error. \Cref{thm:bound_noise,thm:bound_hki,prop:bound_lyap_partial} then control the noise term and the history-dependent term, leading to the total Lyapunov recursion in \Cref{thm:Lyapunov_main}.
Using the resulting Lyapunov recursion, we then move on to prove 
the results stated in the main text.

\subsection{Setup}
For convenience, we repeat the algorithm here and establish some basic facts:
\begin{align}
    m^k &= (1 - \alpha_k)m^{k-1} + \alpha_k \widehat T(z^k; \xi^k, \zeta^k) \label{eq:mom_est} \\
    z^{k+1} &= P_C(x^k - \gamma m^k) \label{eq:pred} \\
    x^{k+1} &= P_C(x^k - \underbrace{\theta\alpha_{k+1}}_{\eqqcolon \omega_k} \widehat T(z^{k+1}; \xi^{k+1}, \zeta^{k+1})) \label{eq:corr}
\end{align}
with $x^0 \in C$ deterministic, $z^0 = x^0$, $m^{-1} = 0$, $\alpha_0 = 1$. We define the filtration \((\gF_{k})_{k \in \N}\) corresponding to our sampling scheme:
\[
  \gF_{-1} \coloneq \sigma(\emptyset), \qquad \gF_{k} = \gF_{k-1} \vee \sigma(\xi^{k}, \zeta^{k}), \qquad k \in \N.
\]
It is clear that $x^k, m^k, z^{k+1}$ are $\gF_k$-measurable. For ease of notation, we also define the error between the estimator $\widehat T$ and its conditional mean at $z^k$, see \Cref{thm:stoch_bounds}
\begin{equation} \label{eq:error_def}
    \varepsilon^{k} \coloneq \widehat T(z^{k}; \xi^k, \zeta^k) - c(z^{k})T(z^{k}).
\end{equation}
Note that since $z^{k+1}$ is $\gF_{k}$-measurable and $\xi^{k+1}, \zeta^{k+1}$ are independent 
of $\gF_{k}$,
by standard conditioning under independence \citep[\S 9.10]{williams1991probability},
\cref{thm:stoch_bounds} gives
\begin{equation} \label{eq:unb_eps}
    \E[\varepsilon^{k} \mid \gF_{k-1}] = 0, \qquad 
    \E[\|\varepsilon^{k}\|^2 \mid \gF_{k-1}] \leq \frac{2\sigma^2}{L_0^2}, \qquad k \in \N.
\end{equation}
In the predictor step, we use Polyak momentum to reduce variance and as such we will have to bound the difference between $T(x^k)$ and its momentum approximation $m^k$.
Since the conditional mean of $\widehat T(z^k; \xi^k, \zeta^k)$ is not $T(z^k)$ but rather $c(z^k)T(z^k)$ in light of \Cref{thm:stoch_bounds}, we build the expected mean tracker $\lambda_kT(x^k)$ where
\begin{equation} \label{eq:def_lamb}
   \lambda_0 \coloneqq c(x^0) = c(z^0),\qquad \lambda_{k+1} \coloneqq (1-\alpha_{k+1})\lambda_{k} + \alpha_{k+1} \cz,
   \qquad k \in \N
\end{equation}
and define the momentum tracking error
\begin{equation} \label{eq:mom_res}
    e^k \coloneqq m^k - \lambda_k T(x^k).
\end{equation}

For our analysis we will also require a scaled version of the subgradient that was chosen when computing \(z^{k+1}\) in \eqref{eq:pred}:
\begin{equation} \label{eq:def_n}
    n^{-1} \coloneqq 0, \qquad n^k \coloneqq \frac{x^k - \gamma m^k - z^{k+1}}{\gamma \lambda_k} \in N_C(z^{k+1}), \qquad k \in \N.
\end{equation}
The scaling is chosen to ensure cancellation in the proof of \Cref{thm:bound_rk}. Moreover, proceeding in the standard way of stochastic (past) extragradient methods leads to very large constants and suboptimal bounds for the stepsizes. For that reason we try to utilize \eqref{eq:unb_eps} as much as possible in order to cancel out cross-products, which leads to the following quantity showing up:
\begin{equation} \label{eq:def_r}
    r^k\coloneqq x^k - z^{k+1} - \omega_k\cz (T(z^{k+1}) + n^k).
\end{equation}

Before providing the full proof, we record some trivial inequalities due to \eqref{eq:def_lamb}, \eqref{eq:c_bound} and the parameter choices in \Cref{assum:param_choice}.
\begin{lemma}
    In the setting of \Cref{thm:main}, we have for every $k\geq0$ 
    that $\lambda_{k}$ is $\gF_{k-1}$-measurable. Moreover, almost surely,
    \begin{equation} \label{eq:c_lam_bound}
    \begin{aligned}
       \frac{1}{1+\delta}&\leq \cz  \leq  1+\delta, \quad & \frac{1}{1+\delta} &\leq \lambda_k\leq 1+\delta,
    \end{aligned}
    \end{equation}
    \begin{equation} \label{eq:param_bounds}
        \gamma\leq\gamma(1+\delta) \leq \frac{1}{4}, \qquad \theta(1+\delta) \leq \theta(1+\delta)^2 \leq \frac{\gamma}{3} \leq \frac{1}{12} 
    \end{equation}
    \begin{equation} \label{eq:alpha_bound}
        (1 - \alpha_k) \left(1 + \frac{3\alpha_k}{4} + \frac{\alpha_k^2}{2}\right) \leq 1 - \frac{\alpha_k}{4} - \frac{\alpha_k^2}{4}, \qquad (1 - \alpha_k)\left(1 + \frac{\alpha_k}{3}\right) \leq 1.
    \end{equation}
\end{lemma}

\subsection{Useful results}
We will make use of several standard results which we state here for convenience. For $a,b,d\in \mathbb{R}^n$ and $e > 0$, we have
\begin{align} 
    2\langle a, b\rangle &\leq e\|a\|^2 + \frac{1}{e} \|b\|^2, \label{eq:fenchel_young} \\
    \|a + b\|^2 &\leq (1+e)\|a\|^2 + (1+\tfrac1e) \|b\|^2, \label{eq:young} \\
    2\langle a, b\rangle &= \|a\|^2 + \|b\|^2 - \|a - b\|^2 \label{eq:pyth}, \\
    -\|a\|^2 - 2\langle d, a+ b\rangle &= -\|a + d\|^2 - \|b\|^2 + \|b - d\|^2. \label{eq:pyth2}
\end{align}

Let \(C\subseteq \mathbb{R}^n\) be a closed, convex, nonempty set and \(y\in \mathbb{R}^n\). Then, for \(P_C(y) \coloneqq \argmin_{x \in C}\|x-y\|^2\) the (necessarily unique) projection of \(y\) onto \(C\), the following standard inequality holds:
\begin{equation} \label{eq:proj_ineq}
    \|P_C(y)-x\|^2 \leq \|y-x\|^2 - \|y-P_C(y)\|^2 \qquad \forall x \in C.
\end{equation}

Let $(\Omega, \cF, \Pr)$ be a probability space. For any two events $A, B \in \cF$,
\begin{equation}\label{eq:prob-decomposition}
  \Pr(A) \leq \Pr(\comp{B}) + \Pr(A \cap B).
\end{equation}

\subsection{Oracle bounds}

We start by bounding the statistical error induced by 
the stochastic estimator $\widehat{T}$.

\subsection{Proof of \texorpdfstring{\Cref{thm:stoch_bounds}}{Lemma~\ref{thm:stoch_bounds}}}
\begin{proof}
    By definition of $\widehat T$ in \eqref{eq:def_hat_T} and the hypotheses:
    \[
        \E[\widehat T(x; \xi, \zeta) ] = F(x) \E \left[ \frac{1}{L_0 + L_1 \|G(x, \zeta)\|} \right] = T(x) \E\left[\frac{L_0 + L_1 \|F(x)\|}{L_0 + L_1 \|G(x, \zeta)\|} \right] = c(x)T(x).
    \]

    To recover the bounds on $c(x)$ in \eqref{eq:c_bound},
    note that the triangle inequality ensures \(\|F(x)\| \leq \|G(x, \zeta)\| + \|F(x) - G(x, \zeta)\|\). Therefore,
    \begin{align*}
        c(x) &= \E \left[ \frac{L_0 + L_1 \|F(x)\|}{L_0 + L_1 \|G(x, \zeta)\|}\right] \leq \E \left[ \frac{L_0 + L_1 \|G(x, \zeta)\| + L_1 \|F(x) - G(x,\zeta)\|}{L_0 + L_1 \|G(x, \zeta)\|} \right]\\
        &\leq 1 + \frac{L_1}{L_0}\E[\|F(x)-G(x, \zeta)\|] \leq 1 + \frac{L_1}{L_0} \sqrt{\E [\|F(x) - G(x, \zeta)\|^2]} \leq 1 + \delta,
    \end{align*}
    where we used $L_0 \leq L_0 + L_1\|G(x, \zeta)\|$ in the second inequality, Jensen's inequality in the third inequality and \Cref{assum:bounded_variance} in the last inequality. Regarding the lower bound on $c(x)$, it follows from the triangle inequality that \(L_0 + L_1\|G(x, \zeta)\| \leq L_0 + L_1\|F(x)\| + L_1\|G(x, \zeta)-F(x)\|\), which leads to 
    \[
        \E[L_0 + L_1\|G(x, \zeta)\|] \leq L_0 + L_1\|F(x)\| + L_1\sigma
    \]
    by a similar argument. Applying Jensen's inequality for the convex function \(t \mapsto 1/t\) on \((0, +\infty)\) then yields
    \[
        c(x) \geq \frac{L_0 + L_1\|F(x)\|}{\E[L_0 + L_1\|G(x, \zeta)\|]} \geq \frac{1}{1+\delta},
    \]
    concluding the proof for the bounds in \eqref{eq:c_bound}.
    
    We now move on to proving the bounded variance of \(\widehat T\) in \eqref{eq:cond_bound_var}.
    By the least-squares property of the mean and \eqref{eq:c_bound},
    we have
    \[
        \E[\|\widehat T(x; \xi, \zeta) - c(x) T(x)\|^2] \leq \E[\|\widehat T(x; \xi, \zeta) - T(x)\|^2],
    \]
    so it suffices to upper bound the right-hand side.
    Adding and subtracting $\frac{F(x)}{L_0 + L_1 \|G(x, \zeta)\|}$ yields
    \[
        \widehat T(x; \xi, \zeta) - T(x) = \frac{G(x, \xi) - F(x)}{L_0 + L_1\|G(x, \zeta)\|} + \left(\frac{F(x)}{L_0 + L_1 \|G(x, \zeta)\|} - \frac{F(x)}{L_0 + L_1\|F(x)\|}\right).
    \]
    It follows that
    \begin{multline*}
        \E[\|\widehat T(x; \xi, \zeta) - T(x)\|^2] = \E\left[\left\| \frac{G(x, \xi) - F(x)}{L_0 + L_1\|G(x, \zeta)\|}\right\|^2 \right] \\+ \E \left[\left\| \frac{F(x)}{L_0 + L_1 \|G(x, \zeta)\|} - \frac{F(x)}{L_0 + L_1\|F(x)\|}\right\|^2 \right]
    \end{multline*}
    where the cross-term vanishes by conditioning first on $\sigma(\zeta)$,
    and using \Cref{assum:bounded_variance} and the independence hypotheses.
    The first term on the right-hand side is upper bounded by 
    $\frac{1}{L_0^2} \E[\|G(x, \xi) - F(x)\|^2] \leq \frac{\sigma^2}{L_0^2}$ by \Cref{assum:bounded_variance}. For the second term, we calculate
    \[
        \left\| \frac{F(x)}{L_0 + L_1 \|G(x, \zeta)\|} - \frac{F(x)}{L_0 + L_1\|F(x)\|}\right\| = \frac{L_1\|F(x)\|  | \|F(x)\| - \|G(x, \zeta)\||}{(L_0 + L_1 \|G(x, \zeta)\|) (L_0 + L_1 \|F(x)\|)}.
    \]
    By the reverse triangle inequality, we have $| \|F(x)\| - \|G(x, \zeta)\|| \leq \|F(x) - G(x, \zeta)\|$. Further, since $\frac{L_1\|F(x)\|}{L_0 + L_1\|F(x)\|} \leq 1$ and $L_0 \leq L_0 + L_1 \|G(x, \zeta)\|$, it follows that
    \[
         \left\| \frac{F(x)}{L_0 + L_1 \|G(x, \zeta)\|} - \frac{F(x)}{L_0 + L_1\|F(x)\|}\right\|^2 \leq \frac{1}{L_0^2} \|F(x) - G(x, \zeta)\|^2.
    \]
    Taking expectation, and using \Cref{assum:bounded_variance} yields
    \[
        \E \left[\left\| \frac{F(x)}{L_0 + L_1\|G(x, \zeta)\|} - \frac{F(x)}{L_0+ L_1 \|F(x)\|}\right\|^2  \right] \leq \frac{1}{L_0^2} \E[\|F(x) - G(x, \zeta)\|^2 ] \leq \frac{\sigma^2}{L_0^2}.
    \]
    We conclude that 
    \[
    \E [\|\widehat T(x; \xi, \zeta) - c(x)T(x)\|^2] \leq \E[\|\widehat T(x; \xi, \zeta) - T(x)\|^2] \leq \frac{2\sigma^2}{L_0^2},
    \] showing the first part of \eqref{eq:cond_bound_var}.
    Finally, write $\widehat T(x; \xi, \zeta) = (\widehat T(x; \xi, \zeta) - c(x) T(x)) + c(x)T(x)$, so that
    \[
        \E[\|\widehat T(x; \xi, \zeta)\|^2] = \E[\|\widehat T(x; \xi, \zeta) - c(x) T(x)\|^2] + c(x)^2 \|T(x)\|^2
    \]
    where the cross-term vanishes due to \eqref{eq:c_bound}. The second part of \eqref{eq:cond_bound_var} now follows from the first part of \eqref{eq:cond_bound_var}.
\end{proof}

As a consequence of \cref{thm:stoch_bounds}, we next establish 
integrability of several quantities of interest, ensuring the existence 
of all expectations in the subsequent analysis.

\begin{lemma}\label{lem:integrability}
  Let $x^\star\in S$ be fixed.
  Under the standing assumptions and independent sampling scheme,
  for every $k\geq0$,
  \[
    \E\left[
      \nrm{x^k-x^\star}^2+\nrm{z^k-x^\star}^2
      +\nrm{m^k}^2
      +\nrm{\widehat T(z^k;\xi^k,\zeta^k)}^2
      +\nrm{n^k}^2
    \right]<\infty.
  \]
\end{lemma}
\begin{proof}
    Let $X$ be a random vector with $\E[\nrm{X-x^\star}^2]<\infty$ and let $\xi,\zeta$ be i.i.d.\ samples independent of $X$. It follows from \eqref{eq:cond_bound_var} and \Cref{assum:prec_lip_minty} that
    \begin{align*}
        \E[\nrm{\widehat T(X; \xi, \zeta)}^2] &\leq (1+\delta)^2 \E[\nrm{T(x)}^2] + \frac{2\sigma^2}{L_0^2} \\
        &\leq 2(1+\delta)^2 \left(\E[\nrm{X - x^\star}^2] + \nrm{T(x^\star)}^2\right) + \frac{2\sigma^2}{L_0^2} < \infty.
    \end{align*}
    Therefore, since $x^0=z^0$ and $m^{-1}=0$ are deterministic,
    the momentum recursion \eqref{eq:mom_est}
    and nonexpansiveness of $P_C$ in \eqref{eq:pred}--\eqref{eq:corr}
    yield finite second moments of the iterates, momentum vectors,
    and oracle estimates by induction. 
\end{proof}

By \cref{lem:integrability}, all quantities appearing inside
expectations in the analysis below are integrable. Indeed, Lipschitz
continuity of $T$ and the bounds \eqref{eq:c_lam_bound} imply that
$\varepsilon^k$, $e^k$, and $r^k$, which are defined in
\eqref{eq:error_def}, \eqref{eq:mom_res}, and \eqref{eq:def_r},
respectively, have finite second moments, as do
$T(z^{k+1})+n^k$ and $\avgnoise^k$
defined in \cref{thm:bound_error}. Therefore, also $h_{k,i}$, $H_k$,
and $V_k(x^\star)$, defined in \eqref{eq:def_h} and
\eqref{eq:def_lyap}, have finite expectations.
Further, all inner products appearing inside expectations are integrable in light of Cauchy--Schwarz.

\subsection{Lyapunov recursion}
We start by establishing a one-step descent estimate for the distance to a solution.
\begin{proposition} \label{prop:dist_est}
    Let $x^\star \in S$.  In the setting of \Cref{thm:main}, for $k\geq0$, we have
    \[
        \E[\|x^{k+1} - x^\star\|^2 \mid \gF_k] \leq \|x^k - x^\star\|^2 - \|x^k - z^{k+1}\|^2 + \|r^k\|^2  + \omega_k^2 \frac{2\sigma^2}{L_0^2}.
    \]
\end{proposition}
\begin{proof}
    From the projection inequality \eqref{eq:proj_ineq} with \(x=x^\star, y=x^k-\omega_k\widehat T(z^{k+1}; \xi^{k+1},\zeta^{k+1})\) and \(P_C(y) \overset{\eqref{eq:corr}}{=} x^{k+1}\) along with the Pythagorean theorem, we have:
    \begin{align*}
        \|x^{k+1}-x^\star\|^2 
        &\leq  \|x^k-x^\star\|^2 - 2\omega_k\langle \widehat T(z^{k+1};  \xi^{k+1},\zeta^{k+1}),x^{k+1}-x^\star\rangle - \|x^{k+1}-x^k\|^2
        \\
        \overset{\eqref{eq:error_def}}&{=} \|x^k-x^\star\|^2 - \|x^{k+1}-x^k\|^2 \underbrace{- 2\omega_k\langle \varepsilon^{k+1}+c(z^{k+1})T(z^{k+1}) , x^{k+1}-x^\star\rangle}_{A_k}.
    \end{align*}
    For the inner product term  
    \(A_k\), adding and subtracting \(z^{k+1}\), we get:
    \begin{align*}
        A_k \leq - 2\omega_k\langle \varepsilon^{k+1}+c(z^{k+1})T(z^{k+1}) , x^{k+1}-z^{k+1}\rangle - 2\omega_k\langle \varepsilon^{k+1}, z^{k+1}-x^\star\rangle,
    \end{align*}
    where we have also used the nonnegativity of \(c(z^{k+1})\langle T(z^{k+1}),z^{k+1}-x^\star \rangle\). Indeed, $\cz$ is trivially nonnegative. Further, \Cref{assum:prec_lip_minty} ensures that \Cref{def:minty} holds. Taking the zero normal in \Cref{def:minty} gives $\langle F(z^{k+1}), z^{k+1} - x^\star\rangle\geq 0$, hence $\langle T(z^{k+1}), z^{k+1} - x^\star\rangle\geq 0$ due to $\frac{1}{L_0 + L_1\|F(z^{k+1})\|}$ being positive.
    
    Further, by adding and subtracting $n^k$ as defined in \eqref{eq:def_n}, we obtain
    \begin{multline*}
        -2\omega_kc(z^{k+1})\langle T(z^{k+1}),x^{k+1}-z^{k+1} \rangle 
        = 2\omega_kc(z^{k+1})[-\langle T(z^{k+1})+n^k,x^{k+1}-z^{k+1} \rangle
        \\
        +\langle n^k,x^{k+1}-z^{k+1} \rangle].
    \end{multline*}
    Since \(n^k \in N_C(z^{k+1})\) by \eqref{eq:def_n} and \(x^{k+1} \in C\), it follows from the definition of the normal cone that the second inner product term on the r.h.s.~of the equality above is nonpositive. Therefore,
    \begin{multline*}
        \|x^{k+1}-x^\star\|^2 
        \leq \|x^k-x^\star\|^2 - \|x^{k+1}-x^k\|^2 - 2\omega_k\langle \varepsilon^{k+1}, z^{k+1}-x^\star\rangle 
        \\
        \underbrace{- 2\omega_k\langle c(z^{k+1})(T(z^{k+1})+n^k) + \varepsilon^{k+1},x^{k+1}-z^{k+1} \rangle}_{B_k}.
    \end{multline*}
    Now applying \eqref{eq:pyth2} with $a = x^{k+1} - x^k$, $b = x^k - z^{k+1}$ and $d = \omega_k (\cz(T(z^{k+1}) + n^k) + \varepsilon^{k+1})$, we have:
    \begin{align*}
        -\|x^{k+1}-x^k\|^2 + B_k \overset{\eqref{eq:pyth2},\eqref{eq:def_r}}&{=} -\|x^{k+1} - x^k + \omega_k(\cz(T(z^{k+1}) + n^k) + \varepsilon^{k+1})\|^2 \\
        &\mspace{200mu} - \|x^k - z^{k+1}\|^2 + \|r^k - \omega_k \varepsilon^{k+1}\|^2 \\
        &\leq -\|x^k-z^{k+1}\|^2 + \|r^k - \omega_k\varepsilon^{k+1}\|^2.
    \end{align*}
    Putting everything together, and using the Pythagorean theorem once more, we get:
    \begin{multline*}
        \|x^{k+1}-x^\star\|^2 
        \leq \|x^k-x^\star\|^2-\|x^k-z^{k+1}\|^2 + \|r^k\|^2 - 2\omega_k \langle r^k,\varepsilon^{k+1} \rangle
       \\ + \omega_k^2\|\varepsilon^{k+1}\|^2- 2\omega_k\langle \varepsilon^{k+1}, z^{k+1}-x^\star\rangle .
    \end{multline*}
    Now take conditional expectation with respect to \(\gF_k\) and note that since \(r^k\) and \(z^{k+1}-x^\star\) are $\gF_k$-measurable, 
    \[
        \E[\langle \varepsilon^{k+1}, z^{k+1}-x^\star\rangle \mid \gF_k] = \E[\langle r^k,\varepsilon^{k+1} \rangle \mid \gF_k] \overset{\eqref{eq:unb_eps}}{=} 0.
    \]
    Furthermore, in light of \eqref{eq:unb_eps}, \(\E[\|\varepsilon^{k+1}\|^2 \mid \gF_k] \leq \tfrac{2\sigma^2}{L_0^2}\). Combining these facts leads to the claimed result.
\end{proof}

\begin{proposition} \label{thm:bound_rk}
   In the setting of \Cref{thm:main}, for $k\geq0$, we have
    \begin{multline*}
        \|r^k\|^2 \leq \|x^k - z^{k+1}\|^2 - \frac{11}{16}\frac{\omega_k\cz}{\gamma\lambda_k} \|x^k - z^{k+1}\|^2 \\ - \frac{2}{3}\omega_k\cz \gamma\lambda_k\|T(z^{k+1}) + n^k\|^2 + \frac{5}{4} \omega_k\gamma \frac{\cz}{\lambda_k}\|e^k\|^2.
    \end{multline*}
\end{proposition}
\begin{proof}
    By the Pythagorean theorem, we have
    \[
        \|r^k\|^2 \overset{\eqref{eq:def_r}}{=} \|x^k - z^{k+1}\|^2 -2\omega_k \cz \langle x^k - z^{k+1}, T(z^{k+1}) + n^k \rangle + \omega_k^2 \cz^2 \|T(z^{k+1}) + n^k\|^2.
    \]
    We now proceed to bound $\langle x^k - z^{k+1}, T(z^{k+1}) + n^k \rangle$. By \eqref{eq:pyth} with $a=x^k - z^{k+1}$ and $b=\gamma\lambda_k(T(z^{k+1}) + n^k)$, we have
    \begin{multline*}
        2\gamma\lambda_k \langle x^k - z^{k+1}, T(z^{k+1}) + n^k\rangle = \|x^k - z^{k+1}\|^2 + \gamma^2\lambda_k^2\|T(z^{k+1}) + n^k\|^2 \\
         - \|x^k - z^{k+1} - \gamma\lambda_k (T(z^{k+1}) + n^k)\|^2.
    \end{multline*}
    Further, the definitions of $n^k$ and $e^k$ \eqref{eq:def_n} and \eqref{eq:mom_res} yield $x^k - z^{k+1} - \gamma\lambda_k(T(z^{k+1}) + n^k) = \gamma e^k - \gamma \lambda_k (T(z^{k+1}) - T(x^k))$. Hence, by Young's inequality \eqref{eq:young} with weight $\gamma\lambda_k$, and the Lipschitz continuity of $T$, we obtain
    \begin{align*}
        \|x^k - z^{k+1} - \gamma \lambda_k(T(z^{k+1}) + n^k)\|^2 &= \|\gamma e^k - \gamma \lambda_k (T(z^{k+1}) - T(x^k))\|^2 \\
        \overset{\eqref{eq:young},\eqref{eq:prec_Lipschitz}}&{\leq} \gamma^2(1+\gamma\lambda_k)\|e^k\|^2 + \gamma \lambda_k(1+\gamma\lambda_k)\|x^k - z^{k+1}\|^2.
    \end{align*}
    It follows that
    \begin{multline*}
        -2\gamma\lambda_k\langle x^k - z^{k+1}, T(z^{k+1}) + n^k\rangle \leq - (1 - \gamma\lambda_k - \gamma^2\lambda_k^2) \|x^k - z^{k+1}\|^2 \\ - \gamma^2\lambda_k^2 \|T(z^{k+1}) + n^k\|^2 + \gamma^2(1+\gamma\lambda_k) \|e^k\|^2.
    \end{multline*}
    Multiplying by $\frac{\omega_k\cz}{\gamma\lambda_k} > 0$ and combining this with our preceding decomposition of $\|r^k\|^2$ leads to
     \begin{multline*}
        \|r^k\|^2 \leq \|x^k - z^{k+1}\|^2 - \frac{\omega_k \cz}{\gamma\lambda_k} (1 - \gamma\lambda_k - \gamma^2\lambda_k^2) \|x^k - z^{k+1}\|^2 \\
        - \omega_k \cz (\gamma\lambda_k - \omega_k\cz) \|T(z^{k+1}) + n^k\|^2 +\omega_k \gamma\frac{\cz}{\lambda_k} (1+\gamma\lambda_k) \|e^k\|^2.
    \end{multline*}

    Moreover, \eqref{eq:c_lam_bound}, \eqref{eq:param_bounds}, $\omega_k = \theta\alpha_{k+1}$, and $0 < \alpha_{k+1} \leq 1$ ensure that $\gamma\lambda_k \leq \gamma (1+\delta) \leq \frac{1}{4}$ and $\frac{\omega_k\cz}{\gamma\lambda_k} \leq \frac{\alpha_{k+1}}{\gamma}\theta(1+\delta)^2 \leq \frac{\alpha_{k+1}}{3} \leq \frac{1}{3}$. Therefore, 
    \[
        1 - \gamma\lambda_k - \gamma^2\lambda_k^2 \geq 1 - \frac{1}{4} - \frac{1}{16} = \frac{11}{16},\quad 
        \gamma\lambda_k - \omega_k\cz = \gamma\lambda_k \left(1 - \frac{\omega_k\cz}{\gamma\lambda_k}\right) \geq \frac{2}{3}\gamma\lambda_k.
    \]
    Inserting these bounds in our preceding bound on $\|r^k\|^2$, and using $1+\gamma\lambda_k\leq\frac54$, yields the desired result.
\end{proof}
In the r.h.s.~of \Cref{thm:bound_rk} we have two controllable terms and the one involving \(\|e^k\|^2\). Using a recursion directly on this quantity does not seem to lead to convergence and we thus proceed to obtain upper bounds in \Cref{thm:bound_error}.

First, we unroll the momentum recursion. Define the standard momentum weights
\begin{equation} \label{eq:def_mom_weights}
    w_{0,0} \coloneqq 1, \qquad w_{k,0} \coloneqq \prod_{l=1}^k (1-\alpha_l),\qquad w_{k,i} \coloneqq \alpha_i \prod_{l=i+1}^k (1-\alpha_l), \quad 1\leq i\leq k,
\end{equation}
where we adopt the usual convention that the empty product equals one. Due to the convex combination, $\sum_{i=0}^k w_{k,i} = 1$ for each $k \in \N$. It follows from \eqref{eq:mom_est} that
\begin{align*}
    m^0 &= \widehat T(z^0; \xi^0,\zeta^0) = w_{0,0}\widehat T(z^0; \xi^0,\zeta^0)
    \\
    m^1 &= (1-\alpha_1)m^0 + \alpha_1 \widehat T(z^1; \xi^1, \zeta^1) = w_{1,0}\widehat T(z^0; \xi^0,\zeta^0) + w_{1,1}\widehat T(z^1; \xi^1, \zeta^1) \\
        &\mspace{10mu}\vdots \\
    m^k &= (1-\alpha_k)m^{k-1} + \alpha_k \widehat T(z^k; \xi^k, \zeta^k) =  \sum_{i=0}^k w_{k,i}\widehat T(z^i; \xi^i, \zeta^i).
\end{align*}
Similarly, we have $\lambda_k = \sum_{i=0}^k w_{k,i}c(z^i)$ by \eqref{eq:def_lamb} due to satisfying the same recursion. We are now ready to bound the momentum tracking error by an averaged noise term and a history-dependent term.

\begin{proposition}  \label{thm:bound_error} Define $\avgnoise^k \coloneqq \sum_{i=0}^k w_{k,i} \varepsilon^i$. In the setting of \Cref{thm:main}, for $k\geq0$, we have
    \begin{equation} \label{eq:ek_bound_young}
        \|e^k\|^2 \leq 2\|\avgnoise^k\|^2 + 2\lambda_k \sum_{i=0}^k w_{k,i}c(z^i)\|x^k-z^i\|^2.
    \end{equation}
\end{proposition}
\begin{proof}
     By our preceding recursions for $m^k$ and $\lambda_k$, and \eqref{eq:mom_res} and \eqref{eq:error_def}, we obtain
    \[
        e^k = \sum_{i=0}^k w_{k,i} (\varepsilon^i + c(z^i)T(z^i)) - \sum_{i=0}^k w_{k,i}c(z^i) T(x^k) = \avgnoise^k + \sum_{i=0}^k w_{k,i}c(z^i) (T(z^i) - T(x^k)).
    \]
    It follows that
    \[
        \|e^k\|^2 \overset{\eqref{eq:young}}{\leq} 2\|\avgnoise^k\|^2 + 2\lambda_k \sum_{i=0}^kw_{k,i}c(z^i)\|T(z^i)-T(x^k)\|^2 \overset{\eqref{eq:prec_Lipschitz}}{\leq}  2\|\avgnoise^k\|^2 + 2\lambda_k\sum_{i=0}^kw_{k,i}c(z^i)\|z^i-x^k\|^2,
    \]
    where we also used Jensen's inequality in the first inequality in light of \(\sum_{i=0}^k w_{k,i}c(z^i)/\lambda_k = 1\).
\end{proof}

\begin{proposition} \label{thm:bound_noise}
    In the setting of \Cref{thm:main}, for $k\geq0$, we have
    \[
        \E[\|\avgnoise^{k+1}\|^2 \mid \gF_k] \leq \left(1 - \alpha_{k+1} \right) \|\avgnoise^k\|^2 + \alpha_{k+1}^2 \frac{2\sigma^2}{L_0^2}.
    \]
\end{proposition}
\begin{proof}
    Note that \(\avgnoise^{k+1} = (1 - \alpha_{k+1})\avgnoise^k + \alpha_{k+1} \varepsilon^{k+1}\) by the definition of the momentum weights \(w_{k,i}\). Therefore, by the Pythagorean theorem
    \[
        \|\avgnoise^{k+1}\|^2 = (1-\alpha_{k+1})^2 \|\avgnoise^k\|^2 + 2\alpha_{k+1} (1-\alpha_{k+1}) \langle \avgnoise^k, \varepsilon^{k+1} \rangle + \alpha_{k+1}^2 \|\varepsilon^{k+1}\|^2.
    \]
    Taking conditional expectation with respect to \(\gF_k\), noting that \(\avgnoise^k\) is $\gF_k$-measurable and using \eqref{eq:unb_eps} and $(1-\alpha_{k+1})^2 \leq 1 - \alpha_{k+1}$ due to $\alpha_{k+1}\leq 1$ leads to the claimed result.
\end{proof}

In \Cref{thm:bound_noise} we obtained a recursive bound for \(\avgnoise^k\). We now move on to bounding \(\sum_{i=0}^k w_{k,i}c(z^i) (T(z^i) - T(x^k))\). Due to the Lipschitz continuity of \(T\), the important part is \(\|x^k-z^i\|^2\). We thus define our partial Lyapunov function in the following way:
\begin{equation} \label{eq:def_h}
    H_k \coloneq \sum_{i=0}^k w_{k,i}c(z^i)h_{k,i}, \quad \textnormal{with} \quad  h_{k,i} \coloneq \|x^{k} - z^i\|^2 + 4\theta^2(1+\delta)^2 \|T(z^i) + n^{i-1}\|^2.
\end{equation}

\begin{proposition} \label{thm:bound_hki}
   In the setting of \Cref{thm:main}, we have the following inequality for \(k \ge 0\) and $0\leq i\leq k$:
    \[
        (1-\alpha_{k+1})\E[h_{k+1,i} \mid \gF_k] 
        \leq  \bigl(1 - \tfrac{\alpha_{k+1}(1 + \alpha_{k+1})}{4}\bigr)h_{k,i}
        + \omega_k \cz \|x^k - z^{k+1}\|^2 + (1-\alpha_{k+1}) \omega_k^2 \frac{2\sigma^2}{L_0^2}.
    \]
\end{proposition}
\begin{proof}
    We begin by bounding \(\|x^{k+1}-z^i\|^2\) since the other term in \eqref{eq:def_h} does not depend on \(k\). Note that \(x^{k+1}  \overset{\eqref{eq:corr}}{=} P_C(x^k - \omega_k \widehat T(z^{k+1}, \xi^{k+1},\zeta^{k+1})) \) and since \(n^{i-1} \in N_C(z^i)\), we have \(z^i = P_C(z^i + \omega_k \cz n^{i-1})\). Therefore, using the nonexpansiveness of the projection,
    \begin{align*}
        \|x^{k+1}-z^i\|^2 &\leq \|x^k-z^i -\omega_k \widehat T(z^{k+1},\xi^{k+1},\zeta^{k+1}) -\omega_k \cz n^{i-1}\|^2
        \\
        \overset{\eqref{eq:error_def}}&{=}\|x^k - z^i - \omega_k \cz(T(z^{k+1}) + n^{i-1}) - \omega_k \varepsilon^{k+1}\|^2
        \\
        &= \|\hat r^{k,i}\|^2 +  \omega_k^2\| \varepsilon^{k+1}\|^2 - 2\langle \hat r^{k,i}, \omega_k \varepsilon^{k+1}\rangle
    \end{align*}
    where we defined $\hat r^{k,i} = x^k - z^i-\omega_k \cz(T(z^{k+1}) + n^{i-1})$ for ease of notation. Taking conditional expectation with respect to \(\gF_k\) and using \eqref{eq:unb_eps} we are led to the following bound:
    \begin{equation} \label{eq:corr_dist}
        \E[\|x^{k+1} - z^i\|^2 \mid \gF_k] \leq \|\hat{r}^{k,i}\|^2 + \omega_k^2 \frac{2\sigma^2}{L_0^2}.
    \end{equation}
    Moreover, it follows from the Pythagorean theorem that
    \begin{multline*}
        \|\hat{r}^{k,i}\|^2 = \|x^k - z^i \|^2    + \omega_k^2\cz^2 \|T(z^{k+1}) + n^{i-1}\|^2 \\ 
        - 2\omega_k \cz \langle x^k - z^i, T(z^{k+1}) - T(z^i)\rangle - 2\omega_k \cz \langle x^k - z^i, T(z^i) + n^{i-1}\rangle ,
    \end{multline*}
    where we have also added and subtracted \(T(z^i)\) in the inner product term. We now proceed by bounding all the terms in the r.h.s.~of the equality above. First,
     \begin{align*}
        \omega_k^2\cz^2 \|T(z^{k+1}) + n^{i-1}\|^2 &\overset{\eqref{eq:young}}{\leq} 2\omega_k^2 \cz^2[\|T(z^{k+1})-T(z^i)\|^2 + \|T(z^i)+n^{i-1}\|^2]
        \\
        \overset{\eqref{eq:young},\eqref{eq:prec_Lipschitz},\eqref{eq:c_bound}}&{\leq} 4\omega_k^2 \cz^2[\|z^{k+1}-x^k\|^2+\|x^k-z^i\|^2] 
        \\
        & \mspace{200mu} + 2\omega_k^2 (1+\delta)^2\|T(z^i)+n^{i-1}\|^2.
    \end{align*}
    Second,
    \begin{align*}
        -2\omega_k\cz\langle x^k - z^i, T(z^{k+1}) - T(z^i)\rangle &
        \leq 2\omega_k \cz [\|x^k - z^i\|\|x^k - z^{k+1}\| + \|x^k - z^i\|^2] \\
        &\overset{\eqref{eq:fenchel_young}}{\leq} \omega_k \cz \|x^k - z^{k+1}\|^2 + 3 \omega_k \cz\|x^k - z^i\|^2,
    \end{align*}
    where the first inequality follows by Cauchy--Schwarz, \eqref{eq:prec_Lipschitz} and the triangle inequality. Third, we use $\omega_k = \theta\alpha_{k+1}$, \eqref{eq:fenchel_young} with $e = \frac{\alpha_{k+1}}{2}$ and \eqref{eq:c_lam_bound} to obtain
    \[
        -2\omega_k \cz\langle x^k - z^i, T(z^i) + n^{i-1}\rangle
      \leq \frac{\alpha_{k+1}}{2} \|x^k - z^i\|^2  + 2\alpha_{k+1}\theta^2(1+\delta)^2 \|T(z^i) + n^{i-1}\|^2.
    \]
   
    Putting everything together and using \(\omega_k = \theta \alpha_{k+1}\) and $\omega_k\cz \leq \alpha_{k+1}\theta(1+\delta) \leq \frac{\alpha_{k+1}}{12}$, see \eqref{eq:param_bounds}, we obtain
    \begin{multline*}
        \E[h_{k+1,i} \mid \gF_k]\leq \mathcolor{Blue4}{\left(1 + \frac{3\alpha_{k+1}}{4} + \frac{\alpha_{k+1}^2}{36}\right)} \|x^k - z^i\|^2 + \omega_k\cz \left(1 + \frac{\alpha_{k+1}}{3}\right)\|x^k - z^{k+1}\|^2  \\
         + \mathcolor{Blue4}{\left( 1 + \frac{\alpha_{k+1}}{2} + \frac{\alpha_{k+1}^2}{2}\right)} 4 \theta^2 (1+\delta)^2 \|T(z^i) + n^{i-1}\|^2 + \omega_k^2 \frac{2\sigma^2}{L_0^2}.
    \end{multline*}
    We note that \textcolor{Blue4}{the coefficients of the first and third terms} can both be bounded by $1 + \frac{3\alpha_{k+1}}{4} + \frac{\alpha_{k+1}^2}{2}$, hence
    \[
        \E[h_{k+1,i} \mid \gF_k] \leq \left(1 + \frac{3\alpha_{k+1}}{4} + \frac{\alpha_{k+1}^2}{2}\right)h_{k,i} + \omega_k\cz \left( 1 + \frac{\alpha_{k+1}}{3}\right)\|x^k - z^{k+1}\|^2 + \omega_k^2 \frac{2\sigma^2}{L_0^2}.
    \]
    Using \eqref{eq:alpha_bound} yields the claimed result.
\end{proof}

\begin{proposition} \label{prop:bound_lyap_partial}
In the setting of \Cref{thm:main}, we have the following inequality for \(k \ge 0\):
    \begin{multline*}
        \E[H_{k+1} \mid \gF_k] \leq \left(1 - \tfrac{\alpha_{k+1}(1+\alpha_{k+1})}{4}\right) H_k + \omega_k\cz \lambda_k \|x^k - z^{k+1}\|^2 + \alpha_{k+1} \cz \|r^k\|^2 \\
         +4\omega_k\theta\cz (1+\delta)^2 \|T(z^{k+1}) + n^k\|^2 +\omega_k^2\lambda_{k+1} \frac{2\sigma^2}{L_0^2}.
    \end{multline*}
\end{proposition}
\begin{proof}
    Based on the definition of the momentum weights \eqref{eq:def_mom_weights}, we have for $0\leq i\leq k$ that $w_{k+1,i} = (1-\alpha_{k+1})w_{k,i}$ and $w_{k+1,k+1} =\alpha_{k+1}$, hence
    \[
        \E[H_{k+1} \mid \gF_k] = \sum_{i=0}^k c(z^i) w_{k,i}(1-\alpha_{k+1}) \E[h_{k+1,i}\mid \gF_k] + \alpha_{k+1}\cz \E[h_{k+1,k+1}\mid\gF_k].
    \]
    Applying \Cref{thm:bound_hki}, and using \eqref{eq:def_h} and $\lambda_k = \sum_{i=0}^k w_{k,i}c(z^i) $, we obtain
    \begin{multline*}
        \E[H_{k+1}\mid\gF_k] \leq \left(1 - \tfrac{\alpha_{k+1}(1+\alpha_{k+1})}{4}\right)H_k + \omega_k\cz\lambda_k\|x^k - z^{k+1}\|^2 \\ 
        + (1-\alpha_{k+1})\omega_k^2\lambda_k \frac{2\sigma^2}{L_0^2} + \alpha_{k+1}\cz \E[h_{k+1,k+1}\mid\gF_k].
    \end{multline*}
    Regarding the term $h_{k+1,k+1} = \|x^{k+1} - z^{k+1}\|^2 + 4\theta^2 (1+\delta)^2 \|T(z^{k+1}) + n^k\|^2$, we argue as in the derivation of \eqref{eq:corr_dist}, now with $i=k+1$, to obtain
    \[
        \E[\|x^{k+1} - z^{k+1}\|^2 \mid \gF_k] \leq \|x^k - z^{k+1} - \omega_k\cz (T(z^{k+1}) + n^{k})\|^2 + \omega_k^2 \frac{2\sigma^2}{L_0^2} \overset{\eqref{eq:def_r}}{=} \|r^k\|^2 + \omega_k^2 \frac{2\sigma^2}{L_0^2}.
    \]
    It follows that
    \begin{multline*}
        \alpha_{k+1}\cz \E[h_{k+1,k+1} \mid \gF_k] \leq \alpha_{k+1} \cz\|r^k\|^2 + \alpha_{k+1}\cz\omega_k^2 \frac{2\sigma^2}{L_0^2}
        \\
        + 4\omega_k\theta\cz (1+\delta)^2 \|T(z^{k+1}) + n^k\|^2.
    \end{multline*}
    Putting everything together and using the definition of $\lambda_{k+1}$ leads to the claimed result.
\end{proof}

Having sufficient bounds for all the terms showing up in our analysis, we can now define our Lyapunov function:
\begin{equation} \label{eq:def_lyap}
    V_k(x^\star) := \|x^k - x^\star\|^2 + 4\theta\gamma(1+\delta)^2 \|\avgnoise^k\|^2 + 10\theta\gamma(1+\delta)H_k,
\end{equation}
for some \(x^\star \in S\).

\begin{proposition} \label{thm:Lyapunov_main}
In the setting of \Cref{thm:main}, we have the following inequality for \(k \ge 0\) and for any $x^\star\in S$:
    \begin{multline*}
        \E[V_{k+1}(x^\star)\mid \gF_k] \leq V_k(x^\star) - \frac{\theta\alpha_{k+1}}{96\gamma(1+\delta)^2} \|x^k - z^{k+1}\|^2 \\ - \frac{7}{18}\frac{\theta\gamma\alpha_{k+1}}{(1+\delta)^2} \|T(z^{k+1}) + n^k\|^2 + \frac{53}{12}\theta\gamma\alpha_{k+1}^2 (1 + \delta)^2 \frac{2\sigma^2}{L_0^2}.
    \end{multline*}
\end{proposition}
\begin{proof}
   By the definition of $V_k$ in \eqref{eq:def_lyap}, we have
   \[
        \E[V_{k+1}(x^\star)\mid \gF_k] = \E[\|x^{k+1} - x^\star\|^2 + 4\theta\gamma(1+\delta)^2\|\avgnoise^{k+1}\|^2 + 10\theta\gamma(1+\delta) H_{k+1} \mid \gF_k].
   \]
   Applying linearity of expectation and \Cref{prop:dist_est,thm:bound_noise,prop:bound_lyap_partial} yields
    \begin{multline*}
        \E[V_{k+1}(x^\star)\mid\gF_k] \leq V_k(x^\star) - \|x^k - z^{k+1}\|^2 + C_1\|r^k\|^2 + C_2\|x^k - z^{k+1}\|^2 \\
        + C_3\|T(z^{k+1}) + n^k\|^2 - C_4\|\avgnoise^k\|^2 - C_5H_k + C_6 \frac{2\sigma^2}{L_0^2}
    \end{multline*}
    where
    \[
    \begin{aligned}
        C_1 &= 1 + 10\theta\gamma(1+\delta)\alpha_{k+1}\cz, & C_2 &= 10\theta\gamma(1+\delta)\omega_k\cz\lambda_k, \\
        C_3 &= 40\theta^2\gamma(1+\delta)^3\omega_k\cz, & C_4 &= 4\theta\gamma(1+\delta)^2\alpha_{k+1}, \\
        C_5 &= \frac{5}{2}\theta\gamma(1+\delta)\alpha_{k+1}(1+\alpha_{k+1}), & C_6 &= \omega_k^2 +4\theta\gamma (1+\delta)^2\alpha_{k+1}^2 + 10\theta\gamma(1+\delta)\omega_k^2\lambda_{k+1}.
    \end{aligned}
    \]
    By our bound on $\|r^k\|^2$ from \Cref{thm:bound_rk} and collecting like terms, we obtain
    \begin{align*}
        \E[V_{k+1}(x^\star)\mid\gF_k] &\leq V_k(x^\star) + \mathcolor{Blue4}{\left(C_1 - 1 + C_2 - \frac{11}{16}C_1 \frac{\omega_k\cz}{\gamma\lambda_k} \right) \|x^k - z^{k+1}\|^2} \\
        &\mspace{75mu} + \mathcolor{Red4}{\left(C_3 - \frac{2}{3}C_1\omega_k\cz\gamma\lambda_k\right) \|T(z^{k+1}) + n^k\|^2} \\
        &\mspace{75mu} + \mathcolor{Orange4}{\frac{5}{4} C_1 \omega_k\gamma \frac{\cz}{\lambda_k} \|e^k\|^2 - C_4\|\avgnoise^k\|^2 - C_5H_k} + C_6\frac{2\sigma^2}{L_0^2}.
    \end{align*}
    We first show that the terms corresponding to \textcolor{Orange4}{$\|e^k\|^2$, $\|\avgnoise^k\|^2$ and $H_k$} together become nonpositive, so they may be discarded. Indeed, \eqref{eq:ek_bound_young} and \eqref{eq:def_h} ensure that $\|e^k\|^2 \leq 2\|\avgnoise^k\|^2 + 2\lambda_kH_k$. Moreover, $10\theta\gamma(1+\delta)\cz \leq 10 \gamma \theta (1+\delta)^2 \leq \frac{10}{12} \gamma \leq \frac{5}{24} \leq \frac14$ by \eqref{eq:c_lam_bound} and \eqref{eq:param_bounds}, hence $C_1 \leq 1 + \frac{\alpha_{k+1}}{4}$ and
    \[
        \frac{5}{4} C_1\omega_k \gamma \frac{\cz}{\lambda_k} \|e^k\|^2 \leq \frac{5}{2} \left(1 + \frac{\alpha_{k+1}}{4}\right) \omega_k\gamma ((1+\delta)^2\|\avgnoise^k\|^2 + (1+\delta)H_k)
    \]
    where we also used \eqref{eq:c_lam_bound}.

    Therefore, utilizing the definitions of $C_4$ and $C_5$, and $\omega_k = \theta\alpha_{k+1}$, the group consisting of \textcolor{Orange4}{$\|e^k\|^2$, $\|\avgnoise^k\|^2$ and $H_k$} can be bounded by
    \[
        \omega_k\gamma(1+\delta)^2 \left( \frac{5}{2}\left(1 + \frac{\alpha_{k+1}}{4}\right) - 4\right) \|\avgnoise^k\|^2 + \frac{5}{2} \omega_k \gamma(1+\delta) \left( 1 + \frac{\alpha_{k+1}}{4} - (1+\alpha_{k+1})\right) H_k
    \]
    where the terms inside the parentheses are nonpositive due to $0 < \alpha_{k+1} \leq 1$. It follows that we may discard the terms corresponding to \textcolor{Orange4}{$\|e^k\|^2$, $\|\avgnoise^k\|^2$ and $H_k$}.
    
    For the coefficient corresponding to $\mathcolor{Blue4}{\|x^k - z^{k+1}\|^2}$, we may bound it as
    \[
        \left(C_1 - 1 + C_2 - \frac{11}{16}C_1 \frac{\omega_k\cz}{\gamma\lambda_k} \right) \leq - \frac{\omega_k\cz}{\gamma\lambda_k} \left(\frac{11}{16} - 10\gamma^2(1+\delta)\lambda_k(1+\theta\lambda_k)\right)
    \]
    by utilizing $C_1 \geq 1$ and $C_1 - 1 + C_2 = 10\gamma(1+\delta) \omega_k\cz(1+\theta\lambda_k)$. Similarly, we can bound the coefficient corresponding to $\mathcolor{Red4}{\|T(z^{k+1}) + n^k\|^2}$ via
    \[
        \left(C_3 - \frac{2}{3}C_1\omega_k\cz\gamma\lambda_k\right) \leq -\omega_k\cz \gamma\lambda_k\left(\frac{2}{3} - \frac{40\theta^2(1+\delta)^3}{\lambda_k} \right),
    \]
    In light of \eqref{eq:c_lam_bound} and \eqref{eq:param_bounds}, we have that $10\gamma^2(1+\delta)\lambda_k(1+\theta\lambda_k) \leq 10 (\gamma(1+\delta))^2 (1+ \theta(1+\delta)) \leq 10 \cdot \frac{1}{16} \cdot \frac{13}{12} = \frac{65}{96}$ and $\frac{40\theta^2(1+\delta)^3}{\lambda_k} \leq 40(\theta(1+\delta)^2)^2 \leq 40 \cdot \frac{1}{144} = \frac{5}{18}$. Since $\frac{11}{16} - \frac{65}{96} = \frac{1}{96}$ and $\frac{2}{3} - \frac{5}{18} = \frac{7}{18}$, and using \eqref{eq:c_lam_bound}, the desired coefficients corresponding to $\mathcolor{Blue4}{\|x^k - z^{k+1}\|^2}$ and $\mathcolor{Red4}{\|T (z^{k+1}) + n^k\|^2}$ follow.

    Lastly, we can write $C_6 = \omega_k^2 (1+10\theta\gamma(1+\delta)\lambda_{k+1}) + 4\theta\gamma(1+\delta)^2\alpha_{k+1}^2$. By \eqref{eq:c_lam_bound} and \eqref{eq:param_bounds}, we have $1 + 10\theta\gamma(1+\delta)\lambda_{k+1} \leq 1 + 10 \gamma \theta(1+\delta)^2 \leq 1 + \frac{10}{3}\gamma^2 \leq 1 + \frac{10}{3} \cdot \frac{1}{16} = \frac{29}{24} \leq \frac{5}{4}$. Similarly, $\omega_k^2 = \omega_k \theta\alpha_{k+1} \leq \frac{1}{3}\omega_k\gamma\alpha_{k+1}(1+\delta)^2$ and $4\theta\gamma(1+\delta)^2\alpha_{k+1}^2 = 4 \omega_k \gamma \alpha_{k+1}(1+\delta)^2$. Therefore, $C_6 \leq \left(\frac{5}{12} + 4\right) \omega_k\gamma\alpha_{k+1}(1+\delta)^2 = \frac{53}{12}\omega_k\gamma\alpha_{k+1}(1+\delta)^2$, finishing the proof.
\end{proof}

We now proceed to prove the results stated in the main text.
\subsection{Proof of \texorpdfstring{\Cref{thm:main}}{Theorem~\ref{thm:main}}}
\begin{proof}
    We start from \Cref{thm:Lyapunov_main}, take total expectation, use the tower property and sum from $k=0$ to $k=K-1$ to obtain
    \begin{multline*}
        \E[V_K(x^\star)] \leq \E[V_0(x^\star)] - \frac{\theta}{96\gamma(1+\delta)^2} \sum_{k=0}^{K-1}\alpha_{k+1} \E[\|x^k - z^{k+1}\|^2] \\ - \frac{7\theta\gamma}{18(1+\delta)^2} \sum_{k=0}^{K-1} \alpha_{k+1}\E[\|T(z^{k+1}) + n^k\|^2] + \frac{53}{12} \theta\gamma(1+\delta)^2 \frac{2\sigma^2}{L_0^2} \sum_{k=0}^{{K-1}}\alpha_{k+1}^2.
    \end{multline*}
    Since $V_K(x^\star) \geq 0$, and using $S_K = \sum_{k=0}^{K-1} \alpha_{k+1}$ and $Q_K = \sum_{k=0}^{K-1} \alpha_{k+1}^2$, it follows that
    \begin{equation} \label{eq:thm-main-end}
        \frac{1}{S_K} \sum_{k=0}^{K-1} \alpha_{k+1} \E[\|T(z^{k+1}) + n^k\|^2] \leq \frac{18}{7} \frac{(1+\delta)^2}{\theta\gamma S_K} \E[V_0(x^\star)] + \frac{159}{14}(1 + \delta)^4 \frac{2\sigma^2}{L_0^2} \frac{Q_K}{S_K}
    \end{equation}
    where we also dropped the nonpositive term corresponding to $\E[\|x^k - z^{k+1}\|^2]$. By \eqref{eq:def_lyap}, $\avgnoise^0 = \varepsilon^0$, \eqref{eq:unb_eps}, \eqref{eq:def_h}, $z^0 = x^0$, \eqref{eq:def_n} and \eqref{eq:c_lam_bound},
    the desired result follows with 
    \begin{equation}\label{eq:bound-EV0}
      \overline{V}_0(x^\star) = \|x^0 - x^\star\|^2 + 4\theta\gamma(1+\delta)^2 \tfrac{2\sigma^2}{L_0^2} + 40\theta^3\gamma(1+\delta)^4 \|T(x^0)\|^2.
    \end{equation}
\end{proof}

\subsection{Proof of \texorpdfstring{\Cref{cor:main_rate}}{Corollary~\ref{cor:main_rate}}}
\begin{proof}
    For $\alpha_{k+1} = K^{-1/2}$,
    we have that $S_K = \sum_{k=0}^{K-1}\alpha_{k+1} = \sqrt{K}$ and $Q_K = \sum_{k=0}^{K-1}\alpha_{k+1}^2 = 1$.
    Thus, by \eqref{eq:thm-main-end}, together with the bound on $\E[V_0(x^\star)]$ in \eqref{eq:bound-EV0},
    \begin{equation} \label{eq:res_rate}
        \frac{1}{K} \sum_{k=0}^{K-1} \E[\|T(z^{k+1}) + n^k\|^2] \leq \frac{18}{7}\frac{(1+\delta)^2}{\theta\gamma}\frac{1}{\sqrt{K}} \left(\overline V_0(x^\star) + \frac{53}{12}(1+\delta)^2 \theta\gamma \frac{2\sigma^2}{L_0^2}\right).
    \end{equation}
    Since $n^k \in N_C(z^{k+1})$ by \eqref{eq:def_n}, the desired claim is proven due to $\dist^2(0, (T + N_C)(z^{k+1})) \leq \|T(z^{k+1}) + n^k\|^2$.
\end{proof}

\subsection{Proof of \texorpdfstring{\Cref{thm:main_as}}{Theorem~\ref{thm:main_as}}}
Before moving to the convergence analysis we will require the supermartingale 
convergence theorem as adapted from \citet[Proposition 2]{bertsekas2011incremental}.
\begin{proposition} \label{thm:supermart}
    Let \((Y_k)_{k \in \mathbb N}, (Z_k)_{k \in \mathbb N}\) and \((W_k)_{k \in \mathbb N}\) be sequences of random variables and let \(\gP_k\) be a filtration. Suppose that:
    \begin{enumerate}
        \item The random variables are nonnegative and are \(\gP_k\)-measurable.
        \item For each \(k \geq 0\) we have
        \[
            \E[Y_{k+1} \mid \gP_k] \leq Y_k - Z_k + W_k.
        \]
        \item Almost surely, \(\sum_{k=0}^\infty W_k < +\infty\).
    \end{enumerate}
    Then, we have \(\sum_{k=0}^\infty Z_k < +\infty\) and the sequence \(Y_k\) converges to a nonnegative random variable \(Y\) almost surely.
\end{proposition}
We now move on to the proof of \Cref{thm:main_as}.
\begin{proof}
    Fix \( x^\star \in S\). We aim at using the standard almost supermartingale theorem to obtain almost sure convergence of \(\|x^k- x^\star\|^2\). To this end, we first show that $\avgnoise^k$ and $H_k$ vanish almost surely, and then show convergence of $\|x^k - x^\star\|^2$.

    From \Cref{thm:bound_noise} and $\sum_{k=0}^\infty \alpha_{k+1}^2<\infty$, we may apply \Cref{thm:supermart} with \(Y_k = \|\avgnoise^k\|^2\),  \(Z_k =\alpha_{k+1}\|\avgnoise^k\|^2\) and \(W_k = \alpha_{k+1}^2 \frac{2\sigma^2}{L_0^2}\) to obtain \(\sum_{k=0}^\infty \alpha_{k+1}\|\avgnoise^k\|^2 < +\infty\) almost surely. Furthermore, \(\|\avgnoise^k\|^2\) converges to some finite value almost surely. Since \(\sum_{k=0}^\infty \alpha_{k+1} = +\infty\), \(\|\avgnoise^k\|^2 \to 0\) almost surely.

    Next, from \Cref{thm:Lyapunov_main} and \Cref{thm:supermart} we get \((V_k(x^\star))_{k \in \mathbb N}\) converges to some finite limit a.s., and 
    \begin{equation} \label{eq:bounds_dk}
        \sum_{k=0}^\infty \alpha_{k+1}(\|T(z^{k+1})+n^k\|^2 + \|x^k-z^{k+1}\|^2) < +\infty,\qquad \text{a.s.}
    \end{equation}
    Since $V_k(x^\star) \geq \|x^k - x^\star\|^2$ by \eqref{eq:def_lyap}, and $V_k(x^\star)$ converges to some finite limit a.s., it follows that $(x^k)_{k\in\mathbb{N}}$ is bounded almost surely. For $H_k$, starting from \Cref{prop:bound_lyap_partial} we need to show that the positive error terms in the recursion are summable almost surely in order to apply \Cref{thm:supermart}. We have
    \begin{align*}
        \|r^k\|^2 \overset{\eqref{eq:def_r}, \eqref{eq:young}, \eqref{eq:c_lam_bound}}&{\leq} 2 \|x^k - z^{k+1}\|^2 + 2\omega_k^2(1+\delta)^2 \|T(z^{k+1}) + n^k\|^2.
    \end{align*}
    Therefore, \Cref{prop:bound_lyap_partial} becomes
    \begin{multline*}
        \E[H_{k+1}\mid\gF_k] \leq \left(1 - \tfrac{\alpha_{k+1}(1+\alpha_{k+1})}{4}\right) H_k + (\theta(1+\delta)^2 + 2(1+\delta))\alpha_{k+1} \|x^k - z^{k+1}\|^2 \\
        + 6 \theta^2 (1+\delta)^3 \alpha_{k+1} \|T(z^{k+1}) + n^k\|^2 + \theta^2(1+\delta) \frac{2\sigma^2}{L_0^2} \alpha_{k+1}^2
    \end{multline*}
    where we used \eqref{eq:c_lam_bound}, $\alpha_{k+1}\leq 1$ and $\omega_k=\theta\alpha_{k+1}$. Then, define \(Y_k = H_k\), \(Z_k = \tfrac{\alpha_{k+1}(1+\alpha_{k+1})}{4}H_k\) and 
    \[
        W_k = \alpha_{k+1}(C_1\|x^k - z^{k+1}\|^2+C_2\|T(z^{k+1}) + n^k\|^2) + \theta^2(1+\delta) \frac{2\sigma^2}{L_0^2} \alpha_{k+1}^2,
    \]
    for some positive deterministic constants \(C_1, C_2\). Note that $\sum_{k=0}^\infty W_k < +\infty$ almost surely in light of \eqref{eq:bounds_dk} and $\sum_{k=0}^\infty \alpha_{k+1}^2<\infty$, hence we may invoke \Cref{thm:supermart} to obtain that \((H_k)_{k \in \mathbb{N}}\) converges to a finite limit a.s.\ and \(\sum_{k=0}^\infty \alpha_{k+1}H_k < +\infty\) a.s., hence \(H_k \to 0\) a.s. Therefore, by \eqref{eq:def_lyap} and the fact that $V_k(x^\star)$ converges to some finite limit a.s., \(\|x^k- x^\star\|^2\) converges to some finite limit almost surely.

    From \eqref{eq:bounds_dk} and $\sum_{k=0}^\infty \alpha_{k+1} = +\infty$, it follows that
    \[
        \liminf_{k\to\infty}(\|T(z^{k+1})+n^k\|^2 + \|x^k-z^{k+1}\|^2) = 0\qquad \text{a.s.}
    \]
    Hence, almost surely, there exists a subsequence such that $\|T(z^{k_j + 1}) + n^{k_j}\| \to 0$ and $\|x^{k_j} - z^{k_j + 1}\|\to 0$. Since $(x^k)_{k\in\mathbb{N}}$ is bounded a.s., passing to a subsubsequence, if needed, we have $x^{k_j} \to \tilde{x}$. Since $z^{k_j + 1} \to \tilde{x}$ and $T$ is continuous, we have \(T(z^{k_j+1}) \to T(\tilde x)\) and \(n^{k_j} \to -T(\tilde x) \in N_C(\tilde x)\) a.s.\ by closedness of the graph of $N_C$. Therefore, \(\tilde x \in S\) due to $N_C(\tilde x)$ being a cone. The rest of the proof now follows by the standard arguments found in \citet[Proposition 9]{bertsekas2011incremental}.
\end{proof}

\subsection{Convergence in probability}

Throughout this section, we adopt the hypotheses of \Cref{cor:main_rate},
and fix a deterministic $x^\star\in S$.
Our main goal is to prove \cref{thm:convergence-in-prob}. 
The central challenge is that the established bounds on the residual with 
respect to the preconditioned operator \(T\) do not straightforwardly carry over 
to the original operator \(F\), as they are related by the random
factor \(L_0 + L_1 \| F(z^{J+1}) \|\), which need not be uniformly bounded.
To address this difficulty, we establish in \cref{lem:expectation-bound-randomit} 
a second-moment bound uniform in \(K\).
This ensures that the randomly selected iterate lies in a sufficiently large ball
with high probability. On this ball, continuity of \(F\) bounds the scaling factor.
Combining the residual bound inside the ball with the small probability of lying outside it yields the desired result.

In the unconstrained case, matters simplify significantly as the residuals with 
respect to \(T\) and \(F\) are related by a strictly increasing map.
The proof for \Cref{thm:convergence-in-prob-unconstrained} therefore follows 
straightforwardly from the bounds on \(T\).

\begin{lemma} \label{lem:avg-it-bound}
    We have the following inequality for $K\geq 1$:
    \[
        \frac{1}{K}\sum_{k=0}^{K-1} \E[\|x^k - z^{k+1}\|^2] \leq \frac{96\gamma(1+\delta)^2}{\theta \sqrt{K}}  \left(\overline{V}_0(x^\star) + \frac{53}{12} \theta\gamma(1+\delta)^2\frac{2\sigma^2}{L_0^2}\right)
    \]
\end{lemma}
\begin{proof}
  The proof follows analogously to the proof of \eqref{eq:res_rate}, where in the proof of \Cref{thm:main} we keep the term corresponding to $\E[\|x^k - z^{k+1}\|^2]$ and drop the nonpositive term corresponding to $\E[\|T(z^{k+1}) + n^k\|^2]$ instead.
\end{proof}

\begin{lemma} \label{lem:expectation-bound-randomit}
  Let random variable $J$ be uniformly distributed on $\{0, \dots, K-1\}$ and 
  independent from $\cF_K$. Then, 
  
  \begin{equation}\label{eq:lem-TzJ-bound}
    \E \big[
      \nrm{T(z^{J+1}) + n^J}^2 
    \big]
      \leq 
    \frac{18(1+\delta)^2}{7\theta\gamma\sqrt{K}} \left(\overline{V}_0(x^\star) + \frac{53}{12}(1+\delta)^2 \theta\gamma \frac{2\sigma^2}{L_0^2}\right)
  \end{equation}
  and, for any solution $x^\star \in S$,
  \begin{equation} \label{eq:lem-zJ-bound}
      \E[\|z^{J+1} - x^\star\|^2] \leq \left(2 + \frac{192\gamma(1+\delta)^2}{\theta\sqrt{K}}\right) \left(\overline{V}_0(x^\star) + \frac{53}{12}\theta\gamma(1+\delta)^2 \frac{2\sigma^2}{L_0^2}\right).
  \end{equation}
\end{lemma}
\begin{proof}
    Since $J$ is uniform on $\{0, \dots, K-1\}$ and independent of $\cF_K$,
    and $\nrm{T(z^{k+1}) + n^k}^2$ is $\cF_K$-measurable for all
    $0 \leq k \leq K-1$,
    \[ 
    \E [\nrm{T(z^{J+1}) + n^J}^2 \mid \cF_K] = \frac{1}{K}\sum_{k=0}^{K-1} \nrm{T(z^{k+1}) + n^k}^2.
    \]
    Thus, applying the tower property,
    \begin{align*}
        \E[\|T(z^{J+1}) + n^J\|^2] &= \E[\E[\|T(z^{J+1}) + n^J\|^2 \mid \gF_K]] = \frac{1}{K} \sum_{k=0}^{K-1} \E[\|T(z^{k+1}) + n^k\|^2 ].
    \end{align*}
    Invoking \eqref{eq:res_rate} immediately establishes \eqref{eq:lem-TzJ-bound}.
    By the same argument, \cref{lem:avg-it-bound} implies
    \begin{equation}\label{eq:bound-xJ-zJ1}
      \E [\nrm{x^{J}- z^{J+1}}^2]  \leq \frac{96\gamma(1+\delta)^2}{\theta\sqrt{K}} \left(\overline{V}_0(x^\star) + \frac{53}{12} \theta\gamma(1+\delta)^2\frac{2\sigma^2}{L_0^2}\right).
    \end{equation}
    
  To obtain \eqref{eq:lem-zJ-bound}, first take the total expectation in \cref{thm:Lyapunov_main} 
  and drop the nonpositive terms to obtain $ \E[V_{k+1}(x^\star)] \leq \E [V_k(x^\star)] + \frac{53}{12} \theta \gamma K^{-1}(1+\delta)^2 \frac{2\sigma^2}{L_0^2}$ for $0 \leq k \leq K - 1$. Hence, by induction, $\E[V_k(x^\star)] \leq \E[V_0(x^\star)] +\frac{53}{12}\frac{k\theta\gamma K^{-1}(1+\delta)^22\sigma^2}{L_0^2}$ for $0 \leq k\leq K$. Since $kK^{-1} \leq 1$, we have
  \begin{equation}\label{eq:bound-x-xstar}
    \E[\nrm{x^k-x^\star}^2]
    \leq \E[V_k(x^\star)]
    \leq
    \overline{V}_0(x^\star)
    +\frac{53}{12}\frac{\theta\gamma(1+\delta)^22\sigma^2}{L_0^2},
    \qquad 0\leq k\leq K.
  \end{equation}
  To convert this bound to \eqref{eq:lem-zJ-bound}, 
  add and subtract $x^J$ and use 
  Young's inequality \eqref{eq:young},
  \[
    \nrm{z^{J+1} - x^\star}^2 \leq 2 \nrm{x^J - x^\star}^2 + 2 \nrm{z^{J+1} - x^J}^2.
  \]
  Taking expectations and combining \eqref{eq:bound-x-xstar} with 
  \eqref{eq:bound-xJ-zJ1} yields the claimed result.
\end{proof}

\subsubsection{Proof of \texorpdfstring{\cref{thm:convergence-in-prob}}{Theorem~\ref{thm:convergence-in-prob}}} \label{sec:proof-convergence-in-prob}

\begin{proof}
  Starting from \cref{lem:expectation-bound-randomit} and applying 
  Markov's inequality, we obtain for any $\epsilon > 0$,
  \begin{equation}\label{eq:prob-bound-T}
    \Pr\{ \nrm{T(z^{J+1}) + n^J} > \epsilon \} \leq \frac{C_T}{\sqrt{K}\epsilon^2}, 
  \end{equation}
  with $C_T \coloneq \frac{18(1+\delta)^2}{7\theta \gamma} \Big(\overline{V}_0(x^\star) + \frac{53}{12}\frac{\theta\gamma (1+\delta)^22\sigma^2}{L_0^2}\Big)$.
  Our goal is to convert this
  to a bound on the residual with respect to $F$.
  For any $z \in \mathbb{R}^n$,
  \( 
    F(z) = (L_0 + L_1 \nrm{F(z)}) T(z),
  \) 
  and moreover, if $n \in N_C(z)$, then also $(L_0 + L_1 \nrm{F(z)}) n \in N_C(z)$,
  since $N_C(z)$ is a cone and $L_0 + L_1\nrm{F(z)} > 0$.
  Let $\tilde n^J = (L_0 + L_1 \nrm{F(z^{J+1})})n^J$, then 
  \begin{equation}\label{eq:relation-F-T}
    \nrm{F(z^{J+1}) + \tilde n^J} = (L_0 + L_1 \nrm{F(z^{J+1})})\nrm{T(z^{J+1}) + n^J}.
  \end{equation}
  
  Note that since the factor $L_0 + L_1 \nrm{F(z^{J+1})}$ is a random variable,
  which need not be uniformly bounded,
  \eqref{eq:lem-TzJ-bound}
  does not suffice to bound $\nrm{F(z^{J+1}) + \tilde n^J}$.

  To address this, let $R>0$ and define the event
  \(
    B_R = \{ \nrm{z^{J+1} - x^\star} \leq R \}.
  \)
Moreover, define
\[
  f_R
  \coloneqq
  \sup_{z\in C}
  \,\big\{ L_0+L_1\nrm{F(z)} \mid \nrm{z-x^\star}\leq R \big\}.
\]
Since $C$ is closed, the set
\(
  \left\{
    z\in C : \nrm{z-x^\star}\leq R
  \right\}
\)
is compact, so by continuity of $F$,
\(
  0 < f_R < \infty.
\)
Now, decomposing according to $B_R$,
\[
  \begin{aligned}
    \Pr\{ \nrm{F(z^{J+1}) + \tilde n^J} > \epsilon \} 
    &\overset{\eqref{eq:prob-decomposition}}{\leq}
    \Pr (\comp{B_R} )
    +\Pr ( \{ \nrm{F(z^{J+1}) + \tilde n^J} > \epsilon \} \cap B_R )\\
    &\overset{\eqref{eq:relation-F-T}}{=}  
    \Pr (\comp{B_R} )
    +\Pr (\{ (L_0 + L_1 \nrm{F(z^{J+1})}) \nrm{T(z^{J+1}) + n^J} > \epsilon \} \cap B_R ).
  \end{aligned}
\]
The first term is bounded by applying Markov's inequality to \eqref{eq:lem-zJ-bound},
and using the fact that $K \geq 1$,
\begin{equation}\label{eq:prob-BR-complement}
  \Pr(B_R^c)
  \leq
  \frac{C_z}{R^2},
  \; \text{ with } 
  C_z
  =
  \Big(
    2+\frac{192\gamma(1+\delta)^2}{\theta}
  \Big)
  \Big(
    \overline{V}_0(x^\star)
    +
    \frac{53}{12}\frac{\theta\gamma(1+\delta)^22\sigma^2}{L_0^2}
  \Big).
\end{equation}

For the second term, note that by definition,
\(
  B_R
  \subseteq
  \left\{
    L_0+L_1\nrm{F(z^{J+1})}\leq f_R
  \right\}.
\)
Therefore,
\begin{align}
  \nonumber
  \Pr \Big(
    \big\{ 
      (L_0 + L_1 \nrm{F(z^{J+1})}) \nrm{T(z^{J+1}) + n^J} > \epsilon 
    \big\} \cap B_R 
    \Big)
    &\leq 
  \Pr \{ f_R \nrm{T(z^{J+1}) + n^J} > \epsilon \} \\
  \nonumber
    &= 
    \Pr \big\{ 
    \nrm{T(z^{J+1}) + n^J} > \frac{\epsilon}{f_R} 
  \big\} 
  \\
  \label{eq:prob-bound-good-case}
    \overset{\eqref{eq:prob-bound-T}}&{\leq} \frac{C_T f_R^2}{\sqrt{K}\epsilon^2}.
\end{align}
Plugging in the bounds \eqref{eq:prob-BR-complement} and \eqref{eq:prob-bound-good-case}, we obtain
\begin{align}
  \Pr\left\{
    \nrm{F(z^{J+1})+\tilde n^J}>\epsilon
  \right\}
  &\leq
  \frac{C_z}{R^2}
  +
  \frac{C_T f_R^2}{\sqrt{K}\epsilon^2}.
  \label{eq:prob-bound-F}
\end{align}
Finally, 
fix $\eta\in(0,1)$. Choose $R_\eta>0$ such that
\(
  \frac{C_z}{R_\eta^2}\leq\frac{\eta}{2},
\)
and choose $C_\eta>0$ satisfying
\(
C_\eta\geq\frac{2}{\eta} C_Tf_{R_\eta}^2.
\)
These choices are independent of $K$.
Since $\tilde n^J\in N_C(z^{J+1})$, we have
\[
  \dist\bigl(0,(F+N_C)(z^{J+1})\bigr)
  \leq \nrm{F(z^{J+1})+\tilde n^J}.
\]
Applying \eqref{eq:prob-bound-F} with
$\epsilon=\sqrt{C_\eta}K^{-1/4}$ therefore yields
\[
\begin{aligned}
  \Pr\left\{
    \dist^2\bigl(0,(F+N_C)(z^{J+1})\bigr)
    >C_\eta K^{-1/2}
  \right\}
  &\leq
  \frac{C_z}{R_\eta^2}
  +\frac{C_Tf_{R_\eta}^2}{C_\eta}
  \leq \eta.
\end{aligned}
\]
\end{proof}

\subsubsection{Proof of \texorpdfstring{\Cref{thm:convergence-in-prob-unconstrained}}{Theorem \ref{thm:convergence-in-prob-unconstrained}}}
\begin{proof}
  Since $C=\mathbb{R}^n$, $N_C \equiv \{0\}$, $n^J = 0$ and applying Markov's inequality to 
  \cref{lem:expectation-bound-randomit} yields,
  for any $\epsilon > 0$,
  \begin{equation}\label{eq:prob-bound-T-unc}
    \Pr\{ \nrm{T(z^{J+1})} > \epsilon \} \leq \frac{C_T}{\sqrt{K}\epsilon^2}, 
  \end{equation}
  Since $r \mapsto \frac{r}{L_0 + L_1 r}$ is strictly increasing on $[0, \infty)$, 
  it holds for all $z \in \mathbb{R}^n$, $\epsilon > 0$ that
  \[
    \nrm{F(z)} > \epsilon \iff \nrm{T(z)} = \tfrac{\nrm{F(z)}}{L_0 + L_1 \nrm{F(z)}} > \tfrac{\epsilon}{L_0 + L_1 \epsilon}.
  \]
  Thus, the claim follows directly from \eqref{eq:prob-bound-T-unc}.
\end{proof}

\section{Experimental details and further results} \label{app:additional_experiments}

\subsection{Updates of compared methods}
SEG has the following standard update: sample \(\xi^k\) and compute \(G(x^k,\xi^k)\). Then update the predictor step:
\begin{equation} \tag{PS}
    z^{k} = x^k - \gamma G(x^k, \xi^k).
\end{equation}
Sample \(\zeta^k\) and compute \(G(z^{k},\zeta^k)\). Then update the corrector step:
\begin{equation} \tag{CS}
    x^{k+1} = x^k - \omega G(z^{k},\zeta^k).
\end{equation}
For SOG the update is simpler, since we only use one sample per iterate. Sample \(\xi^k\) and compute \(G(x^k,\xi^k)\). Then, perform the update:
\begin{equation*} \tag{{SOG}}
    x^{k+1} = x^k - \omega G(x^k,\xi^k) - \gamma (G(x^k,\xi^k) - G(x^{k-1},\xi^{k-1})).
\end{equation*}
We will also consider the stochastic past extragradient (SPEG) method, which has iterates
\begin{equation*} \tag{{SPEG}}
    z^{k+1} = x^k - \gamma G(z^k, \xi^k), \qquad x^{k+1} = x^k - \omega G(z^{k+1}, \xi^{k+1}).
\end{equation*}
Eliminating $x^k$ shows that the sequence $z^k$ follows the SOG recursion for $k\geq 1$, with the corresponding initialization.

For the Clipped SEG method, we follow the update of \cite[\S 4]{vankov2025generalized}:
sample \(\xi^k\), set \(a_k = \beta_k/\max\{1,\|G(x^k,\xi^k)\|\}\) and
compute the predictor step 
\[ 
  z^k = x^k - a_k G(x^k,\xi^k).
\] 
Then, sample \(\zeta^k\), and update 
\[ 
  x^{k+1} = x^k - a_k G(z^k,\zeta^k).
\]
Following \cite{vankov2025generalized},
we set $\beta_k = \frac{a}{a + k^q}$, where $a$ and $q$ are 
tuned according to the procedure below. In addition, we ran experiments taking 
$\beta_k = \beta$ constant. However, we consistently found it to underperform 
with respect to the diminishing step size variant.
To avoid clutter, we therefore opted to omit it from the reported results.

\subsection{Detailed setup and tuning procedures} \label{app:details}

\textbf{Tuning procedure for the experiments in \cref{sec:quartic_nonlinear,sec:cubic-game}}
For the proposed method, we consider $L_0 \in \{0.1, 1, 10\}$, and fix $L_1 = 1$.
\(\alpha_k = 0.1\). We tune the remaining parameters with a two-stage grid search over $98$ parameter configurations, each on three seeds.
Half of the runs are 
divided over a logarithmic grid over \([10^{-4}, 1]^2\) 
for $(\gamma, \omega)$, and $[1, 500] \times [0.1, 5]$ for $(a, q)$.
The remaining half is used for refinement a pass over a linear grid within the intervals defined by the
coarse winner's neighboring grid points.
Each candidate is evaluated using three sampling seeds and a maximum of
\(20000\) oracle calls per run, keeping the problem
instance and initial point fixed. The same seeds are reused across
candidates. We score each candidate by averaging the last \(100\)
recorded residuals within each run and then averaging
across the three runs.

\begin{table}[hbpt]
    \centering
    \caption{Selected hyperparameters for the experiments of \cref{sec:quartic_nonlinear}.}
    \label{tab:quartic_game_hyperparameters}
    \small
    \setlength{\tabcolsep}{4pt}
    \resizebox{\linewidth}{!}{%
    \begin{tabular}{@{}ccccccc@{}}
        \toprule
        $\kappa_B$
        & SEG $(\gamma,\omega)$
        & SOG $(\gamma,\omega)$
        & \multicolumn{3}{c}{SPEPG $(\gamma,\omega)$}
        & Cl. SEG $(a,q)$ \\
        \cmidrule(lr){4-6}
        & & & $L_0 = 0.1$ & $L_0 = 1$ & $L_0 = 10$ & \\
        \midrule
        $10$
        & $(0.0464,\,0.0464)$
        & $(0.0442,\,0.0243)$
        & $(0.215,\,0.0100)$
        & $(1,\,0.0464)$
        & $(1,\,0.364)$
        & $(177,\,1.02)$ \\

        $100$
        & $(0.0464,\,0.0464)$
        & $(0.0442,\,0.0243)$
        & $(0.364,\,0.00953)$
        & $(1,\,0.0442)$
        & $(1,\,0.364)$
        & $(63.0,\,0.863)$ \\
        \bottomrule
    \end{tabular}
    }
\end{table}

\begin{table}[hbpt]
    \centering
    \caption{Selected hyperparameters for the experiments of \cref{sec:cubic-game}.}
    \label{tab:cubic_game_hyperparameters}
    \small
    \setlength{\tabcolsep}{4pt}
    \resizebox{\linewidth}{!}{%
    \begin{tabular}{@{}lcccccc@{}}
        \toprule
        Oracle
        & SEG $(\gamma,\omega)$
        & SOG $(\gamma,\omega)$
        & \multicolumn{3}{c}{SPEPG $(\gamma,\omega)$}
        & Cl. SEG $(a,q)$ \\
        \cmidrule(lr){4-6}
        & & & $L_0 = 0.1$ & $L_0 = 1$ & $L_0 = 10$ & \\
        \midrule
        Finite sum
        & $(0.523,\,0.0785)$
        & $(0.364,\,0.0785)$
        & $(0.215,\,0.0100)$
        & $(1,\,0.113)$
        & $(1,\,0.523)$
        & $(1.91,\,0.450)$ \\

        Student-$t$ noise
        & $(0.523,\,0.0442)$
        & $(0.364,\,0.0317)$
        & $(0.205,\,0.00953)$
        & $(1,\,0.0464)$
        & $(1,\,0.215)$
        & $(9.33,\,0.698)$ \\
        \bottomrule
    \end{tabular}
    }
\end{table}

\textbf{Experimental details on the setup in \cref{sec:cubic-game}}
The problems are generated as
\(
A=U_A D_A U_A^\top,
C=U_C D_C U_C^\top,
\)
where \(U_A,U_C\) are orthogonal matrices obtained from QR decompositions.
\(D_A,D_C\) are diagonal matrices with strictly positive entries,
geometrically spaced between \(1\) and \(1/\kappa_A\), or \(1/\kappa_C\).
The code ensures that the condition-number parameters are at least \(1\).
\(B\) is generated similarly to the example in \cref{sec:quartic_nonlinear}.
We set \(\kappa(A) = \kappa(B) = \kappa(C) = 10\).

In the finite-sum case, the stochastic oracle is generated as
\[
    F_i(x,y)
    =
    \begin{pmatrix}
        \alpha_i \sqrt{x^\top A x}\,Ax + B_i y \\
        \delta_i \sqrt{y^\top C y}\,Cy - B_i^\top x
    \end{pmatrix},
    \qquad i = 1, \dots, N,
\]
where
\(
    \alpha_i = 1 + \sigma_{\alpha}\xi_i,
\)
\(
  \delta_i = 1 + \sigma_{\alpha} \eta_i,
\)
\(
    B_i = B + \sigma_B G_i.
\)
For \(i=1,\dots,N/2\), we independently draw
\(\xi_i,\eta_i \sim \mathcal U[-1,1]\) and matrices \(G_i\)
with independent entries \((G_i)_{jk}\sim\mathcal N(0,1/d)\).
The remaining components are constructed to be their negatives:
\(
    \xi_{i+N/2}=-\xi_i,
\)
\(
    \eta_{i+N/2}=-\eta_i,
\)
\(
    G_{i+N/2}=-G_i.
\)
We then apply a common random permutation to the indices.
Each oracle evaluation independently samples an index \(i\) uniformly
from \(\{1,\dots,N\}\) and returns \(F_i(x,y)\). 
In our experiments, we set \(N=300\), 
\(\sigma_{\mathrm{\alpha}}=0.3\), and \(\sigma_B=0.5\).

For the Student-t noise, we set \(\nu=3\), \(d=30\), and \(\sigma_B=0.5\).
Since \(\nu>1\), the perturbations have zero mean and the 
oracle remains unbiased. Note that in both noise models,
every oracle realization vanishes at the solution \((x,y)=(0,0)\).
For both experiments, we set \(d=30\).

\textbf{Tuning procedure for the WGAN GP experiment.} For our proposed method we choose \(L_0 = 0.1\) and \(L_1 = 1\), \(\alpha = 0.1\), \(\omega = \theta \alpha\). We tune the methods using as metric the mean of the SWD over the last five checkpoints taken at \(90\%, 92.5\%, 95\%, 97.5\%, \) and 100\% of the training process. The initial tuning phase consists of \(100\) runs (\(10 \times 10\) for \(\gamma\) and \(\omega\)) over \(3\) seeds with a logarithmic grid and it is followed by a refinement round with a linear grid for \(100\) more runs. The initial grids for SEG and SOG are \(\gamma \in [10^{-7}, 10^{-1}]\) and \(\omega \in [10^{-6}, 10^{-1}]\). For SPEPG, \(\gamma \in [10^{-3}, 5]\) and \(\theta \in [5 \times 10^{-4}, 5]\). For Clipped SEG, $a\in [5\times 10^{-4}, 5]$, $q\in [0.1, 10]$.

\begin{table}[hbpt]
    \centering
    \caption{Selected hyperparameters for the WGAN GP experiments.}
    \label{tab:wgangp_hyperparameters}
    \resizebox{\linewidth}{!}{\begin{tabular}{@{}ccccc@{}}
        \toprule
        $s_D$
        & SEG $(\gamma,\omega)$
        & SOG $(\gamma,\omega)$
        & SPEPG $(\gamma,\theta)$ & Cl.\ SEG $(a, q)$ \\
        \midrule
        $1$
        & $(1{\times}10^{-3},\,7.7{\times}10^{-3})$
        & $(4.7{\times}10^{-2},\,7.8{\times}10^{-3})$
        & $(2.17,\,3.96{\times}10^{-1})$ & $(0.41, 0.14)$\\

        $4$
        & $(2.1{\times}10^{-4},\,1.6{\times}10^{-4})$
        & $(1.2{\times}10^{-5},\,1.27{\times}10^{-3})$
        & $(7.53\times 10^{-1},\,7.54\times 10^{-1})$ & $(0.33, 0.83)$ \\
        \bottomrule
    \end{tabular}}
\end{table}

\subsection{Synthetic experiments}\label{app:synth}

Consider the following generalization of \cref{ex:sin}:
\begin{equation} \label{eq:ex_sin_nd}
    F(x, y) = \left(1 + \frac{1}{2} \sin\left(\frac{x^2 + y^2}{2}\right) \right) \begin{pmatrix}
        -y \\ x
    \end{pmatrix}
\end{equation}
where $x,y\in\mathbb{R}$. It is clear that $F(x,y) = 0$ if and only if $(x,y) = 0$, and we also have $\langle F(x,y), (x,y) \rangle = 0$ for all $(x,y)\in\mathbb{R}^2$, hence, for $C = \mathbb{R}^2$, \cref{def:minty} holds. Further, $F$ satisfies \cref{def:prec_lipschitz} with $L_0 = L_1 = 1$, whereas it is not $(L_0, L_1)$-Lipschitz for any finite $L_0$ and $L_1$ by similar arguments as in \cref{ex:sin}. Consider now the stochastic oracle $G(x,y, \varepsilon) = F(x,y) + \frac{\sigma}{\sqrt{2}}(\varepsilon_1, \varepsilon_2)^\top$ where $\varepsilon_1,\varepsilon_2$ are uniformly distributed from $\{-1, 1\}$. It is clear that $\E[G(x,y,\varepsilon)] = F(x,y)$ and $\E[\|G(x,y,\varepsilon) - F(x,y)\|^2] = \sigma^2$, hence \Cref{assum:bounded_variance} is satisfied. We emphasize that existing global convergence guarantees for SEG, SOG, and the method of \cite{vankov2025generalized} do not cover this problem.

We choose $\sigma=0.1$, so that $\gamma=0.22$ and $\theta=0.06$ satisfy \Cref{assum:param_choice}. Further, we let $\alpha_k = \frac{c}{c + k}$ where $c=500$, thus all assumptions of \Cref{thm:main_as} hold. We perform $10^6$ iterations of \cref{alg:polyak-sep} with $\omega_k=\theta\alpha_{k+1}$, starting from $(x^0,y^0)=(3,0)$.

As a second experiment, consider $9$ horizons $K_1, \ldots, K_9$, spaced logarithmically between $4\times 10^6$ and $30 \times 10^6$. For a horizon $K_i$, we run \cref{alg:polyak-sep} starting from $(x^0,y^0)=(3,0)$ in the setting of \cref{cor:main_rate} and compute the residual $R_{K_i} \coloneqq \frac{1}{K_i} \sum_{k=0}^{K_i-1} \|T(z^{k+1})\|^2$. We also compute and plot the upper bound from \eqref{eq:res_rate}, which equals approximately $2122 / \sqrt{K_i}$.

\begin{figure}[htb!]
  \centering
  \begin{minipage}[t]{0.48\textwidth}
    \vspace{0pt}
    \centering
    \input{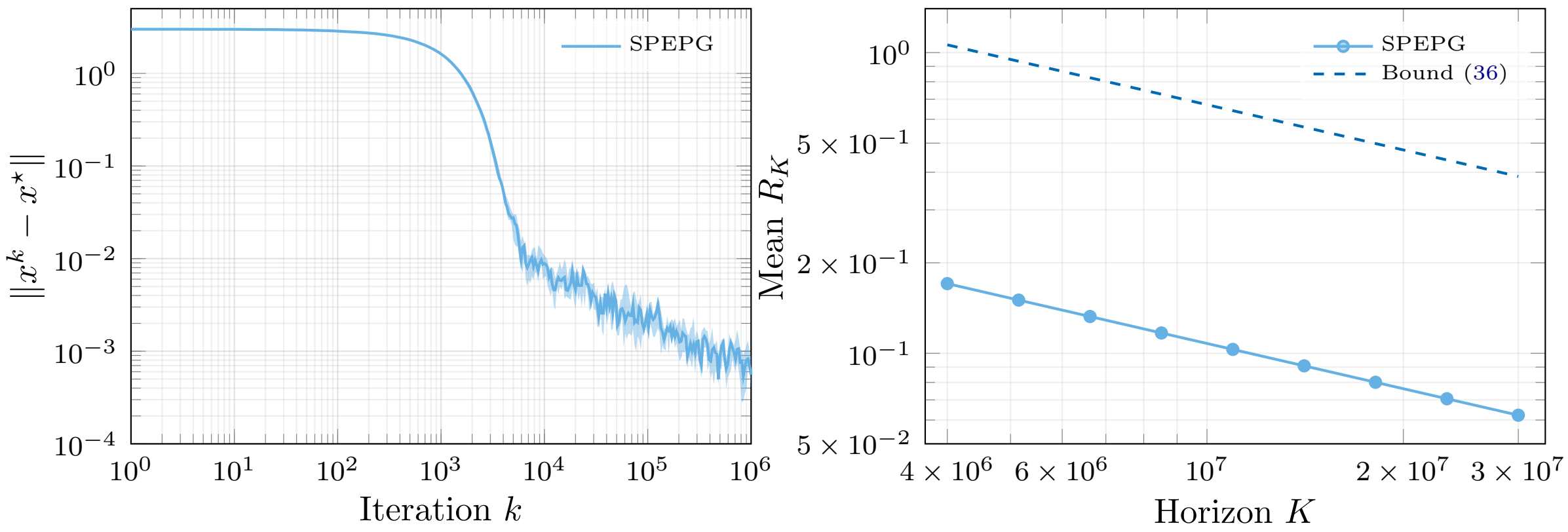}
  \end{minipage}
  \begin{minipage}[t]{0.48\textwidth}
    \vspace{0pt}
    \centering
    \definecolor{niceblue}{RGB}{0,125,210}
\colorlet{method1color}{niceblue!60!white}
\colorlet{method2color}{niceblue!85!black}

\pgfplotsset{compat=1.18}

\begin{tikzpicture}
\begin{axis}[
    set layers=standard,
    width=1.05\textwidth,
    height=0.8\textwidth,
    xlabel={Horizon $K$},
    ylabel={Mean $R_K$},
    try min ticks log=5,
    xmode=log,
    ymode=log,
    grid=both,
    major grid style={opacity=0.25},
    minor grid style={opacity=0.25},
    legend pos=north east,
    xlabel style={font=\small, yshift=3pt},
    ylabel style={font=\small, yshift=-5pt},
    xticklabel style={font=\scriptsize},
    yticklabel style={font=\scriptsize},
    legend style={
        font=\tiny,
        fill opacity=0.5,
        text opacity=1,
        draw=white,
        cells={anchor=west},
        row sep=-2pt,
    },
    enlarge x limits=false,
    xmin=3700000,
    xmax=33000000,
    ymin=0.05,
    ymax=1.4,
    xtick={
        4000000,
        6000000,
        10000000,
        20000000,
        30000000
    },
    xticklabels={
        {$4\times10^{6}$},
        {$6\times10^{6}$},
        {$10^{7}$},
        {$2\times10^{7}$},
        {$3\times10^{7}$}
    },
    minor xtick={
        5000000,
        7000000,
        8000000,
        9000000
    },
    ytick={0.05,0.1,0.2,0.5,1.0},
    yticklabels={
        {$5\times10^{-2}$},
        {$10^{-1}$},
        {$2\times10^{-1}$},
        {$5\times10^{-1}$},
        {$10^{0}$}
    },
    minor ytick={
        0.06,0.07,0.08,0.09,
        0.3,0.4,0.6,0.7,0.8,0.9
    },
]

\addplot[
    color=method1color,
    thick,
    mark=*,
    mark size=1.5pt,
    on layer=axis foreground,
]
coordinates {
    (4000000,0.17050226013011613)
    (5145679,0.15031779866487699)
    (6619502,0.13254895170517772)
    (8515457,0.11684415075187701)
    (10954451,0.10302841920550827)
    (14092021,0.090843149841033533)
    (18128252,0.08009334306481683)
    (23320540,0.070607523102879238)
    (30000000,0.062256364511351289)
};

\addlegendentry{SPEPG}

\addplot[
    color=method2color,
    thick,
    dashed,
    mark=none,
    on layer=axis foreground,
]
coordinates {
    (4000000,1.0611513394507537)
    (5145679,0.93559087430537868)
    (6619502,0.82488738768136338)
    (8515457,0.72728284247488628)
    (10954451,0.64122731437094493)
    (14092021,0.5653542761993926)
    (18128252,0.49845890307800234)
    (23320540,0.43947889208100716)
    (30000000,0.38747768369599966)
};

\addlegendentry{Bound \eqref{eq:res_rate}}

\end{axis}
\end{tikzpicture}
  \end{minipage}
    \caption{
      Convergence plots for the operator in \eqref{eq:ex_sin_nd} to validate \cref{thm:main_as} and \cref{cor:main_rate}.
      (\textit{left}) $10^6$ iterations with $\alpha_k = \frac{c}{c + k}$. The line and shaded areas represent the median and interquartile range over $10$ runs.
      (\textit{right}) Residuals at fixed horizons averaged over $10$ runs, and the theoretical bound.}
    \label{fig:exp_theory}
\end{figure}

The results are presented in \cref{fig:exp_theory}. We observe that the iterates converge to the solution with $\alpha_k = \frac{c}{c+k}$, which is consistent with \Cref{thm:main_as}. The behavior of the residual for the fixed horizons exhibits the predicted $K^{-1/2}$ decay. The constant factor of the theoretical bound could be further improved, though we have not done so to keep the presentation of the analysis accessible.

Finally, we consider the operator \(F(x) = x + \|x\|Jx\) where $J=\begin{psmallmatrix}
    0 & -1 \\ 1 & 0
\end{psmallmatrix}$ and construct the stochastic oracle as before. In this case the operator satisfies the Minty condition and in fact the \(p\)-quasi sharpness condition \cite[Assumption 2.4]{vankov2025generalized} with \(p=2\) and \(\mu=1\). Moreover,
\[
    \|\nabla F(x)\| \le 1 + 2\|x\| \le 1 + 2\|F(x)\|^{1/2}
\]
and thus it satisfies the generalized smoothness condition \cite[Assumption 2.1]{vankov2025generalized}. Therefore, the algorithms of \cite{vankov2025generalized} can also be applied to the root-finding problem. In fact, this problem setting is highly favorable to the methods of \cite{vankov2025generalized}, due to the $p$-quasi sharpness condition.

As for \cref{def:prec_lipschitz}, it is satisfied by \(L_0 = 2\) and \(L_1 = 0.23\) which are the ones we use in the algorithm. With \(\sigma = 0.1\), we get \(\gamma = 0.22\) and \(\theta = 0.067\). We compare the behavior of the methods with decreasing stepsizes such that almost sure convergence is guaranteed and tune the form of the stepsize sequence. More precisely, we tune the methods over 5 random initializations with 150 configurations: 50 stepsize schedules of the form \(\tfrac{\alpha}{k+1}\), with \(\alpha\) logarithmically spaced between \(0.25\) and \(20\), and 100 stepsize schedules of the form \(\tfrac{c}{c+k}\), with \(c\) logarithmically spaced between \(1\) and \(10^8\). We select the schedule minimizing the mean of the unsquared distance to the solution, averaged over the five runs and 41 equally spaced checkpoints in the final \(20\%\) of the tuning budget. In the end, the chosen stepsize sequences were \(\alpha_k = \tfrac{72.2}{72.2 + k}\) for SPEPG and \(\beta_{k} = \tfrac{1.14}{k+1}\) for Clipped SEG. The selected schedules satisfy the assumptions of the respective almost-sure convergence theorems.

\begin{figure}[hbpt]
    \centering
    \input{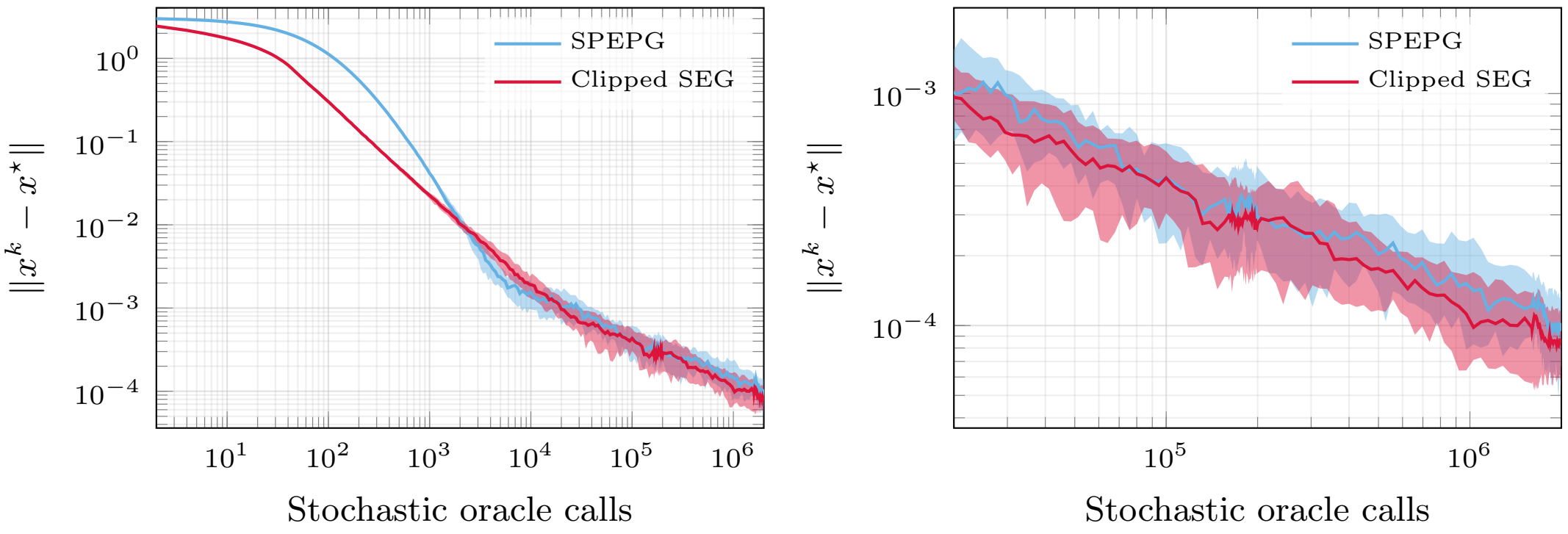}
    \caption{Comparison of SPEPG and Clipped SEG for the operator \(F(x) = x + \|x\|Jx\). The lines and shaded areas represent the median and interquartile range, respectively. (\emph{left}) Full horizon. (\emph{right}) After 20,000 oracle calls.}
    \label{fig:vankov_comp}
\end{figure}

We compare the two methods on 50 runs with an equal budget of \(2 \times 10^6\) stochastic oracle calls per method and plot the distance to the solution against stochastic oracle calls. The results can be found in \cref{fig:vankov_comp}. As expected, in this experiment Clipped SEG outperforms our proposed algorithm, but the two methods perform similarly after a thousand oracle calls.

\subsection{Neural optimal transport}

We consider the problem of estimating the optimal transport (OT) map \(T_\star\) between two probability measures \(P\) and \(Q\) with transport cost \(c(x,y) = h(x-y)\) for a strictly convex \(h \in \gC^1(\mathbb R^n)\), i.e.  the setting of \citet{gangbo1996geometry}. In this case, the map has the form \(T_\star(x) = x-\nabla h^*(\nabla g(x))\) where \(h^*\) is the convex conjugate of \(h\) and \(g\) is a \(c\)-concave \citep[Definition 1.1]{gangbo1996geometry} Kantorovich potential. We choose specifically \(h(x) = \sum_{i=1}^{16} \cosh(x_i)\) to go beyond the standard quadratic cost and make the problem potentially harder, and the \(c\)-concave potential \(g(x) = \sum_{i=1}^{16} \left(a_i \log \cosh ((Rx)_i) - \tfrac{\lambda}{4} (Rx)_i^4\right)\) where \(R\) is a fixed orthogonal matrix, \(a_i = 1-\mu_i\), \(\mu_i = \varepsilon(\mu_{\max}/\varepsilon)^{(i-1)/15}\), \(\varepsilon = 0.05\) and \(\mu_{\max} = 0.8\).

We generate the problem by setting \(P = \gN(0, I_{16})\), computing the exact form of the transport map and getting the target measure via the pushforward, \(Q = T_{\star_\#}P\). We consider the saddle formulation of the OT problem
\[
    \inf_T \sup_f \E_P[c(X,T(X)) - f(T(X))] + \E_Q[f(Y)],
\]
which is obtained by using the definition of \(c\)-conjugates. We parametrize the transport map by \(T_\phi\) and the potential by \(f_\psi\) with neural networks, similarly to \cite{korotin2022neural} and solve the associated minimax problem. We set \(T_\phi(x) = x + r_\phi(x)\) where \(r_\phi\) is a fully connected network with two \(\tanh\) hidden layers: \(16 \to 64 \to 64 \to 16\) while for \(f\) the layers are of the form: \(16 \to 64 \to 64 \to 1\). We estimate the transport map and use as a metric the MSE \(\tfrac{1}{N}\sum_{i=1}^N\|T_\phi(x^i) - T_\star(x^i)\|^2\), \(x_i \sim P\). 

We use a minibatch size of \(256\). For \(\lambda = 0\) we utilize a budget of \(10000\) minibatch oracle calls while for \(\lambda = 0.1\), a budget of \(20000\) minibatch oracle calls. Note that all the methods see an equal amount of data, but SEG and SPEPG perform half the iterations of SOG. We run the experiment ten times and present the results in \Cref{tab:neural_ot}. It can be seen that in terms of MSE SPEPG performs worse on the easier task with \(\lambda = 0\) but much better than the rest of the methods for \(\lambda = 0.1\). The performance in the \(\lambda = 0\) case seems to also depend on the specific instance in some cases performing better and in others much worse. SPEPG also has the best SWD from all the compared methods.

\textbf{Tuning.} For our proposed method we choose \(L_0 = 0.1\) and \(L_1 = 1\), \(\alpha = 0.1\), \(\omega = \theta \alpha\). We use three seeds for tuning the hyperparameters. We tune the methods using as metric the mean of the MSE over the last five checkpoints. The initial tuning consists of \(100\) runs (\(10 \times 10\) for \(\gamma\) and \(\omega\)) over \(3\) seeds with a logarithmic grid and it is followed by a refinement round with a linear grid for \(100\) more runs. The chosen parameters are presented in \Cref{tab:ot_hyperparameters}.

\label{app:experiments}
\begin{table}[hbpt]
    \centering
    \caption{Selected hyperparameters for the neural OT experiments.}
    \label{tab:ot_hyperparameters}
    \begin{tabular}{@{}cccc@{}}
        \toprule
        $\lambda$
        & SEG $(\gamma,\omega)$
        & SOG $(\gamma,\omega)$
        & SPEPG $(\gamma,\theta)$ \\
        \midrule
        $0$
        & $(4.64{\times}10^{-4},\,2.15{\times}10^{-1})$
        & $(3.64{\times}10^{-3},\,1.52{\times}10^{-1})$
        & $(3.48,\,9.57{\times}10^{-1})$ \\

        $0.1$
        & $(1.69{\times}10^{-2},\,1.20{\times}10^{-2})$
        & $(1.00{\times}10^{-2},\,5.76{\times}10^{-3})$
        & $(2.978,\,7.43{\times}10^{-1})$ \\
        \bottomrule
    \end{tabular}
\end{table}

\begin{table}[hbpt]
    \centering
    \caption{
        Final performance on the synthetic neural OT benchmark.
        ``Last-5 MSE'' denotes the mean map MSE over the final five evaluation
        checkpoints, while ``Final MSE'' denotes the map MSE at the end of
        training. Values are reported as mean $\pm$ standard deviation across
        runs. For $\lambda=0.1$, SEG was successful on $6$ of the
        $10$ runs, SOG was successful on $7$ of the $10$ runs and SPEPG was successful on all runs. The reported statistics are computed over the successful runs only.
    }
    \label{tab:ot_results}
    \begin{tabular}{@{}clcccc@{}}
        \toprule
        $\lambda$
        & Method
        & Last-5 MSE $\downarrow$
        & Final MSE $\downarrow$
        & Final SWD $\downarrow$ \\
        \midrule

        \multirow{3}{*}{$0$}
        & SEG
        & $0.04094 \pm 0.00427$
        & $0.03614 \pm 0.00412$
        & $0.02836 \pm 0.00530$ \\

        & SOG
        & $\mathbf{0.02605 \pm 0.00170}$
        & $\mathbf{0.02282 \pm 0.00193}$
        & $0.02309 \pm 0.00281$ \\

        & SPEPG
        & $0.06064 \pm 0.03642$
        & $0.07151 \pm 0.05210$
        & $\mathbf{0.02033 \pm 0.00222}$ \\

        \midrule

        \multirow{3}{*}{$0.1$}
        & SEG
        & $1.01134 \pm 0.34481$
        & $1.26617 \pm 1.05008$
        & $0.11843 \pm 0.11256$ \\

        & SOG
        & $0.76325 \pm 0.05849$
        & $0.70088 \pm 0.02763$
        & $0.04809 \pm 0.00251$ \\

        & SPEPG
        & $\mathbf{0.65410 \pm 0.12448}$
        & $\mathbf{0.65420 \pm 0.12635}$
        & $\mathbf{0.03975 \pm 0.00270}$ \\

        \bottomrule
    \end{tabular}
    \label{tab:neural_ot}
\end{table}

\subsection{GAN training on MNIST}

We compare the methods for training GANs \citep{goodfellow2020generative} to generate MNIST images \citep{deng2012mnist}. The MNIST training set is partitioned into 50000 training examples and 10000 validation examples using a fixed random split. We use a budget of 2048 minibatch oracle calls with batch size 528 and compare using Fr\'echet distance between Gaussian approximations to real and generated classifier features (mFID). For \(G_\psi\) the generator and \(D_\phi\) the discriminator's logit, \(l(t) \coloneqq \ln(1+\exp(t))\), the two players minimize their losses:
\[
    L_G(\psi, \phi) = \E_z[l(-D_\phi(G_\psi(z)))],
\]
\begin{multline*}
    L_D(\psi, \phi) = \E_u[l(-D_\phi(u))] + \E_z[l(D_\phi(G_\psi(z)))]
    \\
    + \tfrac{\lambda_{\text{out}}}{2}(\E_u[D_\phi(u)^2] + \E_z[D_\phi(G_\psi(z))^2]) + \tfrac{\lambda_{\text{wd}}}{2}\|\phi\|^2,
\end{multline*}
where we set \(\lambda_{\text{out}} = 10^{-3}\) and \(\lambda_{\text{wd}} = 10^{-6}\). The operator is simply defined as 
\[
    F(\psi, \phi) = \begin{pmatrix} \nabla_\psi  L_G(\psi, \phi)\\ \nabla_\phi  L_D(\psi, \phi) \end{pmatrix}.
\]

\begin{figure}[hbt]
    \centering
    \begin{subfigure}[t]{0.32\textwidth}
        \centering
        \includegraphics[width=\linewidth]{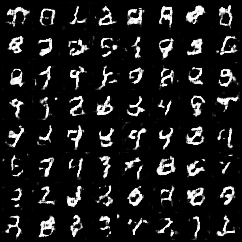}
        \caption{SEG}
        \label{fig:mnist_seg}
    \end{subfigure}
    \hfill
    \begin{subfigure}[t]{0.32\textwidth}
        \centering
        \includegraphics[width=\linewidth]{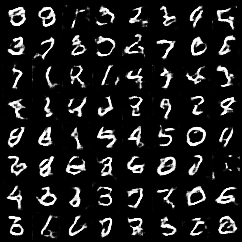}
        \caption{SPEG}
        \label{fig:mnist_speg}
    \end{subfigure}
    \hfill
    \begin{subfigure}[t]{0.32\textwidth}
        \centering
        \includegraphics[width=\linewidth]{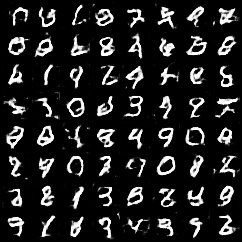}
        \caption{SPEPG}
        \label{fig:mnist_ppeg}
    \end{subfigure}

    \caption{
        Representative MNIST samples generated by SEG, SPEG, and SPEPG after
        2048 stochastic-oracle calls for one test seed. Each panel shows 64 samples from the final generator. All three methods produce recognizable digit-like images, although the qualitative sample quality differs across optimization methods. SEG exhibits more irregular and distorted digit shapes, while SPEG
        and especially SPEPG produce cleaner and more consistently recognizable
        samples.
    }
    \label{fig:mnist_samples}
\end{figure}

The generator and discriminator architectures can be found in the supplementary material. The results are presented in \Cref{tab:mnist} and a visualization of the outputs can be found in \Cref{fig:mnist_samples}. It can be seen that our proposed method achieves the lowest mFID.

\begin{table}[hbt]
    \centering
    \caption{
        Final MNIST generation performance after 2048 stochastic-oracle calls.
        Results are reported as mean $\pm$ standard deviation across held-out
        seeds. Lower mFID is better, while higher digit-classifier score is
        better. SEG failed numerically on 2 of the 5 seeds; consequently, its
        reported mean and standard deviation are computed only over the 3
        successful runs. SPEG and SPEPG completed all 5 runs.}
    \label{tab:mnist_results}
    \begin{tabular}{@{}lcc@{}}
        \toprule
        Method & mFID $\downarrow$ & Digit-classifier score $\uparrow$ \\
        \midrule
        SEG  & $266.78 \pm 10.39$ & $5.793 \pm 0.052$ \\
        SPEG & $82.77 \pm 11.05$  & $6.871 \pm 0.057$ \\
        SPEPG & $\mathbf{48.12 \pm 4.47}$ & $\mathbf{7.407 \pm 0.046}$ \\
        \bottomrule
    \end{tabular}
    \label{tab:mnist}
\end{table}

\textbf{Tuning.} For our proposed method we choose \(L_0 = 0.1\) and \(L_1 = 1\), \(\alpha = 0.1\). We use 3 seeds for tuning the hyperparameters and choose the best setting by averaging the mFID across the three seeds and five late training checkpoints. For this experiment due to the larger scale we followed a different tuning procedure: for all the methods we used a common range for \(\gamma \in [10^{-4}, 0.3]\) and kept the ratio \(\omega / \gamma\) in a range \([0.1, 4]\). In the first phase we use \(32\) configurations for each method obtained from a logarithmically spaced grid, which we then refine with \(32\) more configurations after we obtain the best parameters from the first run. The selected parameters can be found in 
\Cref{tab:gan_hyperparameters}.

\begin{table}[hbpt]
    \centering
    \caption{Selected hyperparameters for the GAN training on MNIST experiments.}
    \label{tab:gan_hyperparameters}
    \begin{tabular}{@{}ccc@{}}
        \toprule
         SEG $(\gamma,\omega)$
        & SPEG $(\gamma,\omega)$
        & SPEPG $(\gamma,\omega)$ \\
        \midrule
         $(7.7{\times}10^{-3},\,2.57{\times}10^{-2})$
        & $(8.1{\times}10^{-3},\,2.75{\times}10^{-2})$
        & $(3.12 {\times} 10^{-2},\,2.77{\times}10^{-1})$
        \\
        \bottomrule
    \end{tabular}
\end{table}
\end{document}